\pdfoutput=1
\documentclass[oneside]{amsart}

\usepackage[margin=1in]{geometry}

\usepackage[T1]{fontenc}

\usepackage[latin9]{inputenc}

\usepackage{bibentry}

\usepackage{amsthm}
\usepackage{amsmath}
\usepackage{amssymb}
\usepackage{nicefrac}

\usepackage{xcolor}

\usepackage{tikz}
\usetikzlibrary{decorations.pathmorphing,decorations.pathreplacing,decorations.shapes}
\usepackage{setspace}
\usepackage[all,cmtip,color]{xy}

\newcommand{\s}[1]{$\mathsf{{#1}}$}
\newcommand{\B}{\square}
\newcommand{\D}{\lozenge}
\newcommand{\bB}{\blacksquare}
\newcommand{\bD}{\blacklozenge}
\newcommand{\A}{\forall}
\newcommand{\E}{\exists}
\newcommand{\nR}{R\!\!\!\!\!\ /}
\newcommand{\nE}{E\!\!\!\!\!\ /}

\newcommand{\world}[1]{\underset{#1}{\bullet}}
\newcommand{\qmax}{\textup{q}\hspace{-.25em}\max}
\newcommand{\strict}[1]{\overset{\rightharpoonup}{#1}}
\newcommand{\xyR}[1]{
    \makeatletter
    \xydef@\xymatrixrowsep@{#1}
    \makeatother
}

\newcommand{\Sub}{\mathsf{Sub}(\varphi)}

\makeatletter
\def\labelbox#1{%
	\hbox{%
		\setbox\z@=\hbox{$\m@th\labelstyle{\,#1\,}$}%
		\setbox\tw@=\hbox{$\m@th\labelstyle\,$}%
		\dimen@=\ht\z@ \advance\dimen@ by \wd\tw@ \ht\z@=\dimen@
		\dimen@=\dp\z@ \advance\dimen@ by \wd\tw@ \dp\z@=\dimen@
		\box\z@
	}%
}
\makeatother

\newtheorem{prop}{Proposition}[section]
\newtheorem{lem}[prop]{Lemma}
\newtheorem*{lem*}{Lemma}
\newtheorem{theorem}[prop]{Theorem}
\newtheorem*{theorem*}{Theorem}
\theoremstyle{definition}
\newtheorem{defn}[prop]{Definition}
\newtheorem{remark}[prop]{Remark}

\newtheorem*{claim*}{Claim}

\newtheorem{cor}[prop]{Corollary}

\begin{document}

\title{Monadic intuitionistic and modal logics admitting provability interpretations}

\author{Guram Bezhanishvili}
\address{New Mexico State University}
\email{guram@nmsu.edu}

\author{Kristina Brantley}
\address{New Mexico State University}
\email{kleifest@nmsu.edu}

\author{Julia Ilin}
\address{University of Amsterdam}
\email{ilin.juli@gmail.com}

\maketitle

\begin{abstract}
The G\"odel translation provides an embedding of the intuitionistic logic $\mathsf{IPC}$ into the modal logic $\mathsf{Grz}$,
which then embeds into the modal logic $\mathsf{GL}$ via the splitting translation. Combined with Solovay's theorem that
$\mathsf{GL}$ is the modal logic of the provability predicate of Peano Arithmetic $\mathsf{PA}$, both $\mathsf{IPC}$ and $\mathsf{Grz}$
admit arithmetical interpretations. When attempting to `lift' these results to the monadic extensions $\mathsf{MIPC}$, $\mathsf{MGrz}$,
and $\mathsf{MGL}$ of these logics, the same techniques no longer work. Following a conjecture made by Esakia, we add an appropriate
version of Casari's formula to these monadic extensions (denoted by a `+'), obtaining that the G\"odel translation embeds
$\mathsf{M^{+}IPC}$ into $\mathsf{M^{+}Grz}$ and the splitting translation embeds $\mathsf{M^{+}Grz}$ into $\mathsf{MGL}$.
As proven by Japaridze, Solovay's result extends to the monadic system $\mathsf{MGL}$, which leads us to an arithmetical interpretation
of both $\mathsf{M^{+}IPC}$ and $\mathsf{M^{+}Grz}$.
\end{abstract}

\tableofcontents


\section{Introduction}

\subsection{Propositional case}

It is well known that the G\"odel translation embeds Intuitionistic Propositional Calculus $\mathsf{IPC}$ into the modal logic $\mathsf{S4}$.
We recall that the \emph{G\"odel translation} is defined as follows:
\begin{itemize}
\item $p^t=\B p$ for a propositional letter $p$;
\item $(\varphi\vee\psi)^t = \varphi^t\vee\psi^t$;
\item $(\varphi\wedge\psi)^t = \varphi^t\wedge\psi^t$;
\item $(\varphi\rightarrow\psi)^t = \B (\varphi^t\to\psi^t)$;
\item $(\neg\varphi)^t = \B (\neg\varphi^t)$.
\end{itemize}
McKinsey and Tarski \cite{McKinsey_Tarski_1948} proved that this translation is full and faithful; that is,
\[
\mathsf{IPC}\vdash\varphi\text{ iff } \mathsf{S4}\vdash \varphi^t.
\]

There are many other normal extensions of $\mathsf{S4}$, called \emph{modal companions} of $\mathsf{IPC}$, in which $\mathsf{IPC}$ is embedded
fully and faithfully. Esakia \cite{Esakia_1976} showed that the largest such companion is \emph{Grzegorczyk's logic} $\mathsf{Grz}$, which
is the normal extension of $\mathsf{S4}$ with the \emph{Grzegorczyk axiom}
\[
\mathsf{grz}=:\B(\B(p\rightarrow\B p)\rightarrow p)\rightarrow p.
\]
Thus, we have
\[
\mathsf{IPC}\vdash\varphi\text{ iff } \mathsf{Grz}\vdash \varphi^t.
\]

Goldblatt \cite{Goldblatt_1978}, Boolos \cite{Boolos_1980_2}, and Kuznetsov and Muravitsky \cite{Kuznetsov_1980} showed that the splitting
translation embeds $\mathsf{Grz}$ into the \emph{G\"odel-L\"ob logic} $\mathsf{GL}$ which is the normal extension of the least normal modal
logic $\mathsf{K}$ with the axiom
\[
\mathsf{gl}:=\B(\B p\to p)\to\B p.
\]
We recall that the splitting translation is defined by ``splitting boxes'' in formulas (see, e.g., \cite[p.~8]{Boolos_1993}); that is, for
a modal formula $\varphi$, let $\B^{+}\varphi$ be the abbreviation of the formula $\varphi\wedge\B\varphi$. Then the \emph{splitting translation}
is defined by letting $\varphi^s$ be the result of replacing all occurrences of $\B$ in $\varphi$ by $\B^{+}$. We then have
\[
\mathsf{Grz}\vdash\varphi \mbox{ iff } \mathsf{GL}\vdash\varphi^s.
\]

Combining these results yields
\[
\mathsf{IPC}\vdash\varphi \mbox{ iff } \mathsf{Grz}\vdash\varphi^t \mbox{ iff } \mathsf{GL}\vdash(\varphi^t)^s.
\]

By Solovay's theorem \cite{Solovay_1976}, $\mathsf{GL}$ can be thought of as the modal logic of the provability predicate in Peano Arithmetic
$\mathsf{PA}$. Thus, both $\mathsf{IPC}$ and $\mathsf{Grz}$ admit provability interpretations.

\subsection{Predicate case}

The G\"odel translation extends to the predicate case by setting
\begin{itemize}
\item $(\forall x\varphi)^t = \B\forall x (\varphi^t)$;
\item $(\exists x\varphi)^t = \exists x(\varphi^t)$.
\end{itemize}
Let $\mathsf{IQC}$ be the intuitionistic predicate calculus and $\mathsf{QS4}$ the predicate $\mathsf{S4}$. Then
\[
\mathsf{IQC}\vdash\varphi\text{ iff } \mathsf{QS4}\vdash \varphi^t,
\]
so the extension of the G\"odel translation to the predicate case remains full and faithful (see, e.g., \cite{Rasiowa_Sikorski_1963}). However, this is virtually the only positive result. Let $\mathsf{QGrz}$ be the predicate $\mathsf{Grz}$ and let $\mathsf{QGL}$ be the
predicate $\mathsf{GL}$. Montagna \cite{Montagna_1984} showed that Solovay's theorem no longer holds for $\mathsf{QGL}$. Moreover, the splitting translation does not embed $\mathsf{QGrz}$ fully and faithfully into $\mathsf{QGL}$ (see below), and as far as we know,
it remains an open problem whether the G\"odel translation embeds $\mathsf{IQC}$ fully and faithfully into $\mathsf{QGrz}$.

\subsection{Monadic case}

In view of the above, Esakia \cite{Esakia_1988} suggested to study these translations for the monadic (one-variable) fragments of
$\mathsf{IQC}$, $\mathsf{QGrz}$, and $\mathsf{QGL}$. The monadic fragment of $\mathsf{IQC}$ was introduced by Prior \cite{Prior_1955}
under the name of $\mathsf{MIPC}$. The monadic fragment of $\mathsf{QS4}$ was studied by Fischer-Servi \cite{Fischer-Servi_1977},
and the monadic fragments of $\mathsf{QGrz}$ and $\mathsf{QGL}$ by Esakia \cite{Esakia_1988}. We denote them by $\mathsf{MS4}$,
$\mathsf{MGrz}$, and $\mathsf{MGL}$, respectively.

Fischer-Servi \cite{Fischer-Servi_1977} proved that the G\"odel translation embeds $\mathsf{MIPC}$ into $\mathsf{MS4}$ fully and
faithfully. As we will see, the G\"odel translation also embeds $\mathsf{MIPC}$ fully and faithfully into $\mathsf{MGrz}$. Japaridze \cite{Dzhaparidze_1988,Dzhaparidze_1990} proved that Solovay's result extends to $\mathsf{MGL}$.
Therefore, to complete the picture, it would be sufficient to show that the splitting translation embeds $\mathsf{MGrz}$ into
$\mathsf{MGL}$ fully and faithfully. However, as was observed by Esakia, this is no longer true. To see this, we recall (see, e.g., \cite{Ono_1987}) that Casari's formula
\[
\mathsf{Cas: }\quad\forall x((P(x)\rightarrow\forall y P(y))\rightarrow\forall y P(y))\rightarrow\forall x P(x)
\]
is valid in an intuitionistic predicate Kripke frame provided the underlying poset is Noetherian. Consider the monadic version of Casari's formula
\[
\mathsf{MCas: }\quad\forall\left(\left(p\rightarrow\forall p\right)\rightarrow\forall p\right)\rightarrow \forall p.
\]
Using the same notation for the G\"odel and splitting translations in the monadic setting, we have that $\mathsf{MGrz}\not\vdash(\mathsf{MCas})^t$ but $\mathsf{MGL}\vdash ((\mathsf{MCas})^t)^s$. This yields that $\mathsf{MGrz}$ does not embed into $\mathsf{MGL}$ faithfully.

Let
\[
\mathsf{M^{+}IPC}=\mathsf{MIPC}+\mathsf{MCas}
\]
be the extension of \s{MIPC} by $\mathsf{MCas}$, and let
\[
\mathsf{M^{+}Grz}=\mathsf{MGrz}+(\mathsf{MCas})^t
\]
be the extension of $\mathsf{MGrz}$ by $(\mathsf{MCas})^t$. Esakia claimed that the translations
\[
\mathsf{IPC}\rightarrow\mathsf{Grz}\rightarrow\mathsf{GL}
\]
are lifted to
\[
\mathsf{M^{+}IPC}\rightarrow\mathsf{M^{+}Grz}\rightarrow\mathsf{MGL}.
\]
Verifying this claim will be our main goal.

\subsection{Main contribution and organization}

Our main result is the following theorem.

\begin{theorem*}
    $\mathsf{M^{+}IPC}\vdash\varphi$ iff $\mathsf{M^{+}Grz}\vdash \varphi^t$ iff $\mathsf{MGL}\vdash (\varphi^t)^s.$
\end{theorem*}

We will prove the theorem semantically. The most challenging part of our argument is in establishing the finite model property for $\mathsf{M^{+}IPC}$ and $\mathsf{M^{+}Grz}$ (see Sections~\ref{sec:5} and~\ref{section:fmpmodalcase}). It was established by Japaridze \cite{Dzhaparidze_1988} that \s{MGL} also has the finite model property. In fact, our technique of proving the finite model property for $\mathsf{M^{+}Grz}$ can be adapted to provide an alternative proof of Japaridze's result for \s{MGL}, but this is not needed for the above theorem.

The paper is organized as follows. Section~\ref{sec:preliminaries} provides a brief overview of monadic logics and their corresponding algebraic and relational semantics. Section~\ref{sec:3} discusses the G\"odel and splitting translations in the monadic setting. In Section~\ref{sec:4} we investigate how the addition of the adapted variations of Casari's formula affect the semantics. In Sections~\ref{sec:5} and~\ref{section:fmpmodalcase} we establish the finite model property for $\mathsf{M^{+}IPC}$ and $\mathsf{M^{+}Grz}$, respectively, using a modified selective filtration, which allows us to conclude the main result stated above.

We use \cite{Modal_Logic} as our standard reference for intuitionistic and modal propositional logic, and \cite{Many_Dim_Modal_Logics} as
our standard reference for intuitionistic modal logics and classical bi-modal systems.


\section{Monadic logics}
\label{sec:preliminaries}

In this section we recall the notion of monadic intuitionistic and modal logics and discuss their algebraic
and frame-based semantics.

\subsection{Monadic intuitionistic logic}

The \emph{monadic intuitionistic propositional calculus} $\mathsf{MIPC}$ was defined by Prior \cite{Prior_1955} and it
was shown by Bull \cite{Bull_1966} that $\mathsf{MIPC}$ axiomatizes the monadic fragment of the predicate intuitionistic logic.
To define $\mathsf{MIPC}$, let $\mathcal L$ be the language of propositional intuitionistic logic, and let $\mathcal{L}_{\forall\exists}$
be the enrichment of $\mathcal{L}$ with the \emph{quantifier modalities} $\A$ and $\E$.\footnote{$\B$ and $\D$ are also frequently used
in place of $\A$ and $\E$, respectively.}

\begin{defn}
$\mathsf{MIPC}$ is the smallest set of $\mathcal{L}_{\forall\exists}$-formulas containing
\begin{itemize}
\item all axioms of $\mathsf{IPC}$,
\item the $\mathsf{S4}$-axioms for $\A$,\footnote{$\A p\to p$, $\A p\to\A\A p$, and $\A(p\wedge q)\leftrightarrow(\A p\wedge\A q)$.}
\item the $\mathsf{S5}$-axioms for $\E$,\footnote{$p\to\E p$, $\E\E p\to\E p$, $\E(p\vee q)\leftrightarrow(\E p\vee\E q)$, and
$\E(\E p\wedge q)\leftrightarrow(\E p\wedge\E q)$.}
\item the connecting axioms $\E p\to\A\E p$ and $\E\A p\to\A p$,
\end{itemize}
and closed under the inference rules of substitution, modus ponens, and $\forall$-necessitation $\frac{\varphi}{\A\varphi}$.
\end{defn}

\begin{remark}
The non-symmetric feature of intuitionistic quantifiers is captured in the fact that while $\E$ is an $\sf S5$-modality, $\A$ is merely
an $\mathsf{S4}$-modality, and the $\A$-counterpart $\A(\A p\vee q)\leftrightarrow(\A p\vee\A q)$ of $\E(\E p\wedge q)\leftrightarrow(\E p\wedge\E q)$
is not provable in $\mathsf{MIPC}$.
\end{remark}

Algebraic semantics for $\mathsf{MIPC}$ is given by monadic Heyting algebras \cite{Monteiro&Varsavsky_1957,Bezhanishvili_1998_2}.

\begin{defn}
A \emph{monadic Heyting algebra} is a triple $(H,\A,\E)$ where
\begin{itemize}
\item $H$ is a Heyting algebra,
\item $\A:H\to H$ is an $\mathsf{S4}$-operator,\footnote{$\A a\le a$, $\A a\le\A\A a$, $\A(a\wedge b)=\A a\wedge\A b$, and $\A 1=1$.}
\item $\E:H\to H$ is an $\mathsf{S5}$-operator,\footnote{$a\le\E a$, $\E\E a\le\E a$, $\E(a\vee b)=\E a\vee\E b$, $\E 0=0$, and
$\E(\E a\wedge b)=\E a\wedge\E b$.}
\item $\E a\le\A\E a$ and $\E\A a\le\A a$.
\end{itemize}
\end{defn}

\begin{remark} \label{rem: H_0}
This in particular implies that the fixpoints of $\A$ and $\E$ are equal and form a Heyting subalgebra of $H$. In fact, every monadic Heyting
algebra can be represented as a pair $(H,H_0)$ where $H_0$ is a Heyting subalgebra of $H$ and the inclusion has both the right ($\A$) and
left $(\E)$ adjoint.
\end{remark}

As usual, propositional letters of $\mathcal{L}_{\forall\exists}$ are evaluated as elements of $H$, the connectives
as the corresponding operations of $H$, and the quantifier modalities as the corresponding modal operators of $H$.
The standard Lindenbaum-Tarski construction then yields:

\begin{theorem}\label{thm: algebraic completeness}
$\mathsf{MIPC}\vdash\varphi\Leftrightarrow\mathfrak{H}\vDash\varphi$ for each monadic Heyting algebra $\mathfrak{H}$.
\end{theorem}

Kripke semantics for $\mathsf{MIPC}$ is an extension of Kripke semantics for $\mathsf{IPC}$ \cite{Ono_1977,Fischer-Servi_1978_2,Esakia_1988}.

\begin{defn}
An $\mathsf{MIPC}$-\emph{frame} is a triple $\mathfrak{F}=(W,R,E)$ where $(W,R)$ is an $\mathsf{IPC}$-frame,\footnote{A nonempty partially
ordered set.} and $E$ is an equivalence relation on $W$ satisfying $(R\circ E)(x)\subseteq (E\circ R)(x)$ for all $x\in W$; that is,
if $xEy$ and $yRz$, then there is $w\in W$ such that $xRw$ and $wEz$.
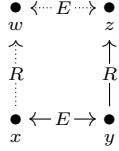
\[\xymatrix{
\world{w}\ar@{<.>}[r]|-{E} & \world{z} \\
\world{x}\ar@{<->}[r]|-{E}\ar@{.>}[u]|-{R} & \world{y}\ar[u]|-{R}
}\]
We refer to this condition as \emph{commutativity}. We will sometimes refer to $R$ as a `vertical relation', and to $E$ as a `horizontal relation',
as depicted in the diagram above.
\end{defn}

Valuations on $\mathsf{MIPC}$-frames are defined as for $\mathsf{IPC}$-frames; that is, a valuation on $\mathfrak F=(W,R,E)$ is an assignment
$\nu$ of propositional letters to $R$-upsets of $\mathfrak F$.\footnote{Recall that $U\subseteq W$ is an \emph{$R$-upset} if $u\in U$ and $uRv$
imply $v\in U$.} As usual, the truth relation in $\mathfrak F$ is defined by induction. The clauses for the connectives $\wedge,\vee,\to,\neg$ are the same
as for $\mathsf{IPC}$-frames:
\begin{alignat*}{3}
&w\vDash_\nu p \quad &\text{iff} & \quad w\in\nu(p); \\
&w\vDash_\nu\varphi\wedge\psi \quad  &\text{iff}& \quad  w\vDash_\nu\varphi \text{ and } w\vDash_\nu\psi;\\
&  w\vDash_\nu\varphi\vee\psi \quad  &\text{iff}&  \quad w\vDash_\nu\varphi \text{ or } w\vDash_\nu\psi;\\
&  w\vDash_\nu\varphi\to\psi \quad  &\text{iff}& \quad  (\text{for all } v)(wRv \text{ and } v\vDash_\nu\varphi \text{ implies } v\vDash_\nu\psi);\\
&  w\vDash_\nu\neg\varphi\quad &\text{iff} & \quad (\text{for all } v)(wRv \text{ implies } v\not\vDash_\nu\varphi).
\intertext{
To extend this to the truth relation for quantifier modalities, we first define a new relation $Q$ as the composition $R\circ E$ on $W$; that is,
$xQy$ iff there is $z\in W$ such that $xRz$ and $zEy$.
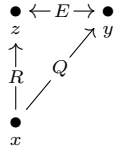
\[\xymatrix{
    \world{z}\ar@{<->}[r]|-{E} & \world{y} \\
    \world{x}\ar[u]|-{R}\ar[ur]|-{Q}
}\]
Then $Q$ is a quasi-order (reflexive and transitive) and $\A,\E$ are interpreted in $\mathfrak F$ as follows:}
& w\vDash_\nu\forall\varphi\quad&\text{iff}&\quad (\text{for all } v)(wQv \text{ implies } v\vDash_\nu\varphi);\\
& w\vDash_\nu\exists\varphi\quad&\text{iff}&\quad(\text{there exists } v)(wEv\text{ and } v\vDash_\nu\varphi).
\end{alignat*}
Sometimes we also write $(\mathfrak F, w) \vDash_\nu\varphi$ to emphasize the underlying frame $\mathfrak F$ or simply $w \vDash \varphi$
in case $\mathfrak F$ and $\nu$ are clear from the context.

There is a close connection between algebraic and relational semantics for $\mathsf{MIPC}$. To see this, let $\mathfrak F=(W,R,E)$ be
an $\mathsf{MIPC}$-frame. For $x\in W$, let
\[
Q(x)=\{y\in W\mid xQy\} \mbox{ and } E(x)=\{y\in W\mid xEy\}.
\]
Set $\mathfrak F^+=(\mathsf{Up}(\mathfrak F),\A,\E)$ where $\mathsf{Up}(\mathfrak F)$ is the Heyting algebra of $R$-upsets of $\mathfrak F$,
and for $U\in\mathsf{Up}(\mathfrak F)$,
\[
\A U=\{x\in W\mid Q(x)\subseteq U\} \mbox{ and } \E U=\{x\in W\mid E(x)\cap U\ne\varnothing\}.
\]
Then $\mathfrak F^+$ is a monadic Heyting algebra, and every monadic Heyting algebra is represented as a subalgebra of such. To see
this, for a monadic Heyting algebra $\mathfrak H=(H,\A,\E)$, let $W$ be the set of prime filters of $H$, let $R$ be the inclusion,
and let $E$ be defined by $\eta E \zeta$ iff $\eta\cap H_0=\zeta\cap H_0$, where we recall that $H_0$ is the fixpoint subalgebra of
$H$ (see Remark~\ref{rem: H_0}).
Then $\mathfrak H_+=(W,R,E)$ is an $\mathsf{MIPC}$-frame (where $\eta Q\zeta$ iff $\eta\cap H_0\subseteq\zeta\cap H_0$) and there is
an embedding $e:\mathfrak H\to\mathfrak (\mathfrak H_+)^+$
given by
\[
e(a)=\{\eta\in\mathfrak H_+\mid a\in\eta\}.
\]

In general, the embedding $e$ is not onto, so to recognize the $e$-image of $H$ in the Heyting algebra of upsets, we introduce the concept of
a descriptive $\mathsf{MIPC}$-frame. One way to do this is to introduce topology on an $\mathsf{MIPC}$-frame.

We recall that a topological space
is a \emph{Stone space} if it is compact Hausdorff and zero-dimensional\footnote{Clopen (closed and open) sets form a basis for the topology.}.
A relation $R$ on a Stone space $W$ is \emph{continuous} if (i)~$R(x)$ is closed for each $x\in W$ and (ii)~$U$ clopen implies $R^{-1}(U)$ is clopen,
where
\[
R^{-1}(U)=\{x\in W\mid xRu \mbox{ for some } u\in U\}.
\]

\begin{defn}
An $\mathsf{MIPC}$-frame $\mathfrak{F}=(W,R,E)$ is a \emph{descriptive $\mathsf{MIPC}$-frame} if
\begin{itemize}
\item $W$ is a Stone space,
\item $R$ and $Q$ are continuous relations,
\item $A$ clopen $R$-upset implies $E(A)$ is a clopen $R$-upset.
\end{itemize}
\end{defn}

\begin{remark}\label{rem: E}
This does not imply that $A$ clopen implies $E(A)$ is clopen; see \cite[p.~32]{Bezhanishvili_1999}. However, we do have that
$A$ closed implies $E(A)$ is closed; see \cite[Lem.~7]{Bezhanishvili_1999}.
\end{remark}

As follows from Esakia's representation of Heyting algebras \cite{Esakia_1974}, for a Heyting algebra $H$, there is a Stone topology on
the set $W$ of prime filters of $H$ generated by the basis
\[
\{e(a)\setminus e(b)\mid a,b\in H\},
\]
the inclusion relation $R$ on $W$ is continuous, and $e$ is a Heyting isomorphism from $H$ onto the Heyting algebra of clopen $R$-upsets of $W$.

By \cite[Thm.~13]{Bezhanishvili_1999}, if $\mathfrak H=(H,\A,\E)$ is a monadic Heyting algebra, then $(W,R,E)$ is a descriptive
$\mathsf{MIPC}$-frame, which we denote by $\mathfrak H_*$, and $e$ is an isomorphism from $\mathfrak H$ onto the monadic Heyting algebra
$(\mathfrak H_*)^*$ of clopen $R$-upsets of $\mathfrak H_*$. Thus, every monadic Heyting algebra can be thought of as the algebra of clopen
$R$-upsets of some descriptive $\mathsf{MIPC}$-frame. This representation together with Theorem~\ref{thm: algebraic completeness} yields:

\begin{theorem}\label{thm: frame completeness}
$\mathsf{MIPC}\vdash\varphi\Leftrightarrow\mathfrak{F}\vDash\varphi$ for each descriptive $\mathsf{MIPC}$-frame $\mathfrak{F}$.
\end{theorem}

If the descriptive $\mathsf{MIPC}$-frame is finite, then the topology is discrete, and hence finite descriptive $\mathsf{MIPC}$-frames are
simply finite $\mathsf{MIPC}$-frames. It is well known that $\mathsf{MIPC}$ has the finite model property:

\begin{theorem}
\label{thm:mipchasfmp}
$\mathsf{MIPC}\vdash\varphi\Leftrightarrow\mathfrak{F}\vDash\varphi$ for each finite $\mathsf{MIPC}$-frame $\mathfrak{F}$.
\end{theorem}

This was first proved by Bull \cite{Bull_1965} using algebraic semantics. However, Bull's proof contained a gap, which was later filled by Fischer-Servi \cite{Fischer-Servi_1978_1} and
Ono \cite{Ono_1977} independently of each other. For a more frame-theoretic proof, using the technique of selective filtration, see \cite[\S 10.3]{Many_Dim_Modal_Logics}.

We finish \S 2.1 by recalling an important property of descriptive $\mathsf{MIPC}$-frames, which will be useful later on.

\begin{defn}\label{defn: max}
Let $\mathfrak{F}=(W,R,E)$ be a descriptive $\mathsf{MIPC}$-frame and let $A\subseteq W$.
    \begin{enumerate}
        \item We say $x\in A$ is \emph{$R$-maximal} in $A$ if $xRy$ and $y\in A$ imply $x=y$.
        \item The \emph{$R$-maximum} of $A$ is the set of all $R$-maximal points of $A$, i.e.,
        \[\max A=\{x\in A \mid xRy\text{ and } y\in A \mbox{ imply } x=y\}.\]
    \end{enumerate}
\end{defn}

The next lemma states that every point in the $E$-saturation of clopen $A$ sees a point that is maximal in the $E$-saturation of $A$.
The proof follows from the result of Fine \cite{Fine_1974} and Esakia \cite{Esakia_1985} that can be phrased as follows: If $A$ is a
closed subset of a descriptive $\mathsf{IPC}$-frame, then for each $x\in A$ there is $y\in\max A$ such that $xRy$. Since $A$ clopen
implies that $E(A)$ is closed (see Remark~\ref{rem: E}), the proof is a consequence of the Fine-Esakia lemma.

\begin{lem}\label{MIPC points can see max}
Let $\mathfrak{F}=(W,R,E)$ be a descriptive $\mathsf{MIPC}$-frame. For each clopen $A$ and $x\in E(A)$, there is $y\in \max E(A)$
such that $xRy$.
\end{lem}

\subsection{Monadic modal logics}

Let $\mathcal{ML}$ be the basic propositional modal language (with one modality $\Box$). As usual, the least normal modal logic will be denoted
by $\mathsf{K}$, and normal modal logics are normal extensions of $\mathsf{K}$.

Let $\mathcal{ML_{\A}}$ be the bimodal language which enriches $\mathcal{ML}$ with the modality $\A$. We use the abbreviation $\E\varphi$ for
$\neg\forall\neg\varphi$.

\begin{defn}
\begin{enumerate}
\item[]
\item The \emph{monadic $\mathsf{K}$} is the least set of $\mathcal{ML_{\A}}$-formulas containing
\begin{itemize}
\item the $\mathsf{K}$-axiom for $\Box$,\footnote{$\B(p\to q)\to(\B p\to\B q)$.}
\item the $\mathsf{S5}$-axioms for $\A$,\footnote{$\A p\to p$, $\A p\to\A\A p$, $\neg\A p\to\A\neg\A p$, and $\A(p\to q)\to(\A p\to\A q)$.}
\item the bridge axiom $\B\A p\to\A\B p$,
\end{itemize}
and closed under $\A$-necessitation $\frac{\varphi}{\A\varphi}$ as well as under the usual rules of substitution, modus ponens, and
$\B$-necessitation. We denote the monadic $\mathsf{K}$ by $\mathsf{MK}$.
\item A \emph{normal extension} of $\mathsf{MK}$ is an extension of $\mathsf{MK}$ which is closed under both $\B$- and $\A$-necessitation.
We call normal extensions of $\mathsf{MK}$ \emph{normal monadic modal logics} or simply \emph{mm-logics}.
\item Let $\mathsf{L}$ be a normal modal logic (in $\mathcal{ML}$). The \emph{least monadic extension} $\mathsf{ML}$ of $\mathsf{L}$ is the
smallest mm-logic containing $\mathsf{MK}\cup\mathsf{L}$.
\end{enumerate}
\end{defn}

\begin{remark}
\begin{enumerate}
\item[]
\item Monadic modal logics are bimodal logics in the language with two modalities $\B,\A$, where $\A$ is an $\mathsf{S5}$-modality.
They correspond to expanding relativized products discussed in \cite[\S 9]{Many_Dim_Modal_Logics}.
\item The axiom $\A\B p\to\B\A p$, which is the converse of the bridge axiom, and is the monadic version of Barcan's formula, is not
provable in $\mathsf{MK}$.
\end{enumerate}
\end{remark}

Algebraic semantics for monadic modal logics is given by monadic modal algebras.

\begin{defn}
A \emph{monadic modal algebra} or simply an \emph{mm-algebra} is a triple $(B,\B,\A)$ where
\begin{itemize}
\item $(B,\B)$ is a modal algebra,\footnote{That is, $B$ is a boolean algebra and $\B:B\to B$ satisfies $\B 1=1$ and $\B(a\wedge b)=\B a\wedge\B b$.}
\item $(B,\A)$ is an $\mathsf{S5}$-algebra,\footnote{That is, $(B,\A)$ is a modal algebra satisfying $\A a\le a$, $\A a\le\A\A a$, and
$\neg\A a\le\A\neg\A a$.}
\item $\B\A a\le\A\B a$.
\end{itemize}
\end{defn}

\begin{remark}
As with monadic Heyting algebras, the $\A$-fixpoints of an mm-algebra $(B,\B,\A)$ form a subalgebra of the modal algebra $(B,\B)$,
and each mm-algebra is represented as a pair $(B,B_0)$ of modal algebras such that the embedding of $B_0$ into $B$ has a right adjoint ($\A$).
\end{remark}

Kripke semantics for mm-logics is given by augmented Kripke frames of Esakia \cite{Esakia_1988}.

\begin{defn}\label{defn:augFrm}
An \emph{augmented Kripke frame} is a triple $\mathfrak{F}=(W,R,E)$ where $(W,R)$ is a Kripke frame\footnote{$W$ is nonempty and $R$ is a binary
relation on $W$.} and $E$ is an equivalence relation on $W$ satisfying \emph{commutativity}, i.e., $(R\circ E)(x)\subseteq (E\circ R)(x)$ for all
$x\in W$; that is, if $xEy$ and $yRz$, then there is $w\in W$ such that $xRw$ and $wEz$.
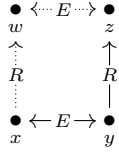
\[\xymatrix{
\world{w}\ar@{<.>}[r]|-{E} & \world{z} \\
\world{x}\ar@{<->}[r]|-{E}\ar@{.>}[u]|-{R} & \world{y}\ar[u]|-{R}
}\]
\end{defn}

As with $\mathsf{MIPC}$-frames, we may refer to $R$ as a `vertical relation,' and to $E$ as a `horizontal relation,' as depicted in the diagram above.

Valuations on augmented Kripke frames are defined analogously to Kripke frames; that is, a valuation $\nu$ on an augmented Kripke frame
$\mathfrak{F}=(W,R,E)$ assigns propositional letters to subsets of $W$.
The truth relation clauses for the connectives $\vee, \neg$, the modality $ \Box$, and its dual $\Diamond$ are defined as for Kripke frames:
\begin{alignat*}{3}
    &x \vDash_{\nu} p \quad &\text{ iff }&\quad x\in\nu(p); \\
    &x \vDash_{\nu}\psi\vee\chi \quad &\text{ iff } & \quad x \vDash_{\nu}\psi\text{ or } x \vDash_{\nu}\chi; \\
    &x \vDash_{\nu}\neg\psi \quad& \text{ iff } & \quad x \not\vDash_{\nu}\psi; \\
    & x \vDash_{\nu}\square\psi \quad & \text{ iff }&\quad (\text{for all } y\in W) (xRy\Rightarrow y\vDash_{\nu}\psi); \\
    &x \vDash_{\nu}\lozenge\psi \quad & \text{ iff }&\quad (\text{there exists } y\in W )(xRy\text{ and } y \vDash_{\nu}\psi).\\
\intertext{The modality $\A$ and its dual $\E$ are interpreted via the relation $E$ as follows:}
& x\vDash_{\nu}\forall\varphi\quad &\text{iff} &\quad (\text{for all } y\in W)(xEy\Rightarrow y\vDash_{\nu}\varphi)\\
& x\vDash_{\nu}\exists\varphi\quad & \text{iff} &\quad(\text{there exists } y\in W)(xEy\text{ and }y\vDash_{\nu}\varphi).
    \end{alignat*}
As in the case of $\mathsf{MIPC}$-frames, we also use the notation $(\mathfrak F, w) \vDash_\nu\varphi$ or $w \vDash \varphi$.

As in the case of $\mathsf{MIPC}$, there is a close connection between algebraic and relational semantics for mm-logics. For an augmented Kripke
frame $\mathfrak F=(W,R,E)$, set $\mathfrak F^+=(\wp(\mathfrak F),\B,\A)$ where $\wp(\mathfrak F)$ is the powerset of $\mathfrak F$, and for
$U\in\mathsf{Up}(\mathfrak F)$,
\[
\B U=\{x\in W\mid R(x)\subseteq U\} \mbox{ and } \A U=\{x\in W\mid E(x)\subseteq U\}.
\]
Then $\mathfrak F^+$ is an mm-algebra, and every mm-algebra is represented as a subalgebra of such. To see
this, for an mm-algebra $\mathfrak B=(B,\B,\A)$, let $W$ be the set of ultrafilters of $B$, and let $R$ and $E$
be defined by
\[
\eta R \zeta \mbox{ iff }\B a\in\eta \mbox{ implies } a\in\zeta \mbox{ and } \eta E \zeta \mbox{ iff } \eta\cap B_0=\zeta\cap B_0.
\]
Then $\mathfrak B_+=(W,R,E)$ is an augmented Kripke frame and there is an embedding $e:\mathfrak B\to\mathfrak (\mathfrak B_+)^+$ given by
\[
e(a)=\{\eta\in\mathfrak B_+\mid a\in\eta\}.
\]

In general, the embedding $e$ is not onto, so to recognize the $e$-image of $\mathfrak B$ in the powerset, we introduce the concept of
a descriptive augmented Kripke frame. As in the case of $\mathsf{MIPC}$, we do this by introducing topology on augmented Kripke frames.

\begin{defn}
An augmented Kripke frame $\mathfrak{F}=(W,R,E)$ is a \emph{descriptive augmented Kripke frame} if $W$ is a Stone space and $R$ and $E$
are continuous relations.
\end{defn}

As follows from the representation of modal algebras, for a modal algebra $B$, there is a Stone topology on
the set $W$ of ultrafilters of $B$ generated by the basis $\{e(a)\mid a\in B\}$, the relation $R$ on $W$ is continuous,
and $e$ is a modal isomorphism from $B$ onto the modal algebra of clopen subsets of $W$.

If $\mathfrak B=(B,\B,\A)$ is an mm-algebra, then $(W,R,E)$ is a descriptive augmented Kripke frame, which we denote by $\mathfrak B_*$,
and $e$ is an isomorphism from $\mathfrak B$ onto the mm-algebra $(\mathfrak B_*)^*$ of clopen subsets of $\mathfrak B_*$. Thus, every
mm-algebra can be thought of as the algebra of clopen subsets of some descriptive augmented Kripke frame.

\subsection{$\mathsf{MS4}$, $\mathsf{MGrz}$, and $\mathsf{MGL}$}

We next focus on the least monadic extension $\mathsf{MS4}$ of the modal logic $\mathsf{S4}$.

\begin{defn}
\begin{enumerate}
\item[]
\item An \emph{$\mathsf{MS4}$-algebra} is an mm-algebra $(B,\B,\A)$ such that $(B,\B)$ is an $\mathsf{S4}$-algebra.
\item An \emph{$\mathsf{MS4}$-frame} is an augmented Kripke frame $\mathfrak{F}=(W,R,E)$ such that $(W,R)$ is an $\mathsf{S4}$-frame.
\item A \emph{descriptive $\mathsf{MS4}$-frame} is a descriptive augmented Kripke frame $\mathfrak{F}=(W,R,E)$ such that $(W,R,E)$ is an
$\mathsf{MS4}$-frame.
\end{enumerate}
\end{defn}

As in the case of $\mathsf{MIPC}$, we have the following standard completeness results:

\begin{theorem}\label{thm: completeness 2}
\begin{enumerate}
\item[]
\item $\mathsf{MS4}\vdash\varphi\Leftrightarrow\mathfrak{B}\vDash\varphi$ for each $\mathsf{MS4}$-algebra $\mathfrak{B}$.
\item $\mathsf{MS4}\vdash\varphi\Leftrightarrow\mathfrak{F}\vDash\varphi$ for each descriptive $\mathsf{MS4}$-frame $\mathfrak{F}$.
\end{enumerate}
\end{theorem}

We also have that $\mathsf{MS4}$ has the finite model property. This can be proved by adopting the algebraic proof of the finite model property
of $\mathsf{MIPC}$ to the setting of $\mathsf{MS4}$ (see \cite{BC_2019}).

\begin{theorem}\label{MS4_Comp_FMP}
    $\mathsf{MS4}\vdash\varphi\Leftrightarrow\mathfrak{F}\vDash\varphi$ for each finite $\mathsf{MS4}$-frame $\mathfrak{F}$.
\end{theorem}

Let $\mathfrak{F}=(W,R,E)$ be a descriptive $\mathsf{MS4}$-frame and $A \subseteq W$. The \emph{$R$-maximal} points of $A$ and
the \emph{$R$-maximum} of $A$ are defined as in Definition~\ref{defn: max}. In the context of $\mathsf{MS4}$-frames, we also
need the notion of quasi-$R$-maximal points.

\begin{defn}
Let $\mathfrak{F}=(W,R,E)$ be a descriptive $\mathsf{MS4}$-frame and $A \subseteq W$.
    \begin{enumerate}
        \item We say $x\in A$ is \emph{quasi-R-maximal} in $A$ if $xRy$ and $y\in A$ imply $yRx$.
        \item The \emph{quasi-R-maximum} of $A$ is the set of all quasi-R-maximal points of $A$, i.e.,
        \[\qmax A=\{x\in A \mid xRy\text{ and }y\in A \text{ imply } yRx\}.\]
    \end{enumerate}
\end{defn}

Note that $\max A\subseteq\qmax A$ as $R$ is reflexive, but not conversely. The following lemma is a consequence of the
Fine-Esakia lemma \cite{Fine_1974,Esakia_1985} for descriptive $\mathsf{S4}$-frames.

\begin{lem}\label{lem:S4canseeqmax}
Let $\mathfrak{F}=(W,R,E)$ be a descriptive $\mathsf{MS4}$-frame.
For each closed $A \subseteq W$ we have $A\subseteq R^{-1}\qmax A$.
\end{lem}

\begin{defn}
\begin{enumerate}
\item[]
\item The \emph{monadic $\mathsf{Grz}$} is the least monadic extension $\mathsf{MGrz}$ of Grzegorczyk's logic $\mathsf{Grz}$.
\item An \emph{$\mathsf{MGrz}$-algebra} is an mm-algebra $(B,\B,\A)$ such that $(B,\B)$ is a $\mathsf{Grz}$-algebra.
\item An \emph{$\mathsf{MGrz}$-frame} is an augmented Kripke frame $\mathfrak{F}=(W,R,E)$ such that $(W,R)$ is a $\mathsf{Grz}$-frame.
\item A \emph{descriptive $\mathsf{MGrz}$-frame} is a descriptive $\mathsf{S4}$-frame $\mathfrak{F}=(W,R,E)$ validating
Grzegorczyk's axiom $\mathsf{grz}$.
\end{enumerate}
\end{defn}

Again, we have the following standard completeness results:

\begin{theorem}\label{thm: completeness 3}
\begin{enumerate}
\item[]
\item $\mathsf{MGrz}\vdash\varphi\Leftrightarrow\mathfrak{B}\vDash\varphi$ for each $\mathsf{MGrz}$-algebra $\mathfrak{B}$.
\item $\mathsf{MGrz}\vdash\varphi\Leftrightarrow\mathfrak{F}\vDash\varphi$ for each descriptive $\mathsf{MGrz}$-frame $\mathfrak{F}$.
\end{enumerate}
\end{theorem}

It is well known that an $\mathsf{S4}$-frame $\mathfrak{F}=(W,R)$ is a $\mathsf{Grz}$-frame iff $R$ is a \emph{Noetherian partial order};
that is, a partial order with no infinite ascending chains (of distinct points). Thus, if $\mathfrak{F}$ is finite, then $\mathfrak{F}$
is a $\mathsf{Grz}$-frame iff $R$ is a partial order.

It is a result of Esakia that a descriptive $\mathsf{S4}$-frame $\mathfrak{F}=(W,R)$ is a descriptive $\mathsf{Grz}$-frame iff for each clopen
$A \subseteq W$ the $R$-maximal and quasi-$R$-maximal points of $A$ coincide.
These results clearly hold for $\mathsf{MGrz}$ as well.

\begin{lem}[\cite{Esakia_1979}] \label{lem:grzchar}
\begin{enumerate}
\item[]
\item Let $\mathfrak{F}=(W,R,E)$ be a descriptive $\mathsf{MS4}$-frame. Then $\mathfrak{F}\vDash\mathbf{grz}$ iff for each clopen
$A$ we have $\qmax A=\max A$.
\item Let $\mathfrak{F}=(W,R,E)$ be a descriptive $\mathsf{MGrz}$-frame. For each clopen $A$ we have $A\subseteq R^{-1}\max A$.
\end{enumerate}
\end{lem}

\begin{defn}
\begin{enumerate}
\item[]
\item The \emph{monadic $\mathsf{GL}$} is the least monadic extension $\mathsf{MGL}$ of the G\"odel-L\"ob logic $\mathsf{GL}$.
\item An \emph{$\mathsf{MGL}$-algebra} is an mm-algebra $(B,\B,\A)$ such that $(B,\B)$ is a $\mathsf{GL}$-algebra.
\item An \emph{$\mathsf{MGL}$-frame} is an augmented Kripke frame $\mathfrak{F}=(W,R,E)$ such that $(W,R)$ is a $\mathsf{GL}$-frame.
\item A \emph{descriptive $\mathsf{MGL}$-frame} is a descriptive augmented Kripke frame $\mathfrak{F}=(W,R,E)$ validating $\mathsf{gl}$.
\end{enumerate}
\end{defn}

As before, we have the following standard completeness results:

\begin{theorem}\label{thm: completeness 4}
\begin{enumerate}
\item[]
\item $\mathsf{MGL}\vdash\varphi\Leftrightarrow\mathfrak{B}\vDash\varphi$ for each $\mathsf{MGL}$-algebra $\mathfrak{B}$.
\item $\mathsf{MGL}\vdash\varphi\Leftrightarrow\mathfrak{F}\vDash\varphi$ for each descriptive $\mathsf{MGL}$-frame $\mathfrak{F}$.
\end{enumerate}
\end{theorem}

It is well known that a Kripke frame $\mathfrak{F}=(W,R)$ is a $\mathsf{GL}$-frame iff $R$ is transitive and \emph{dually well founded}
(no infinite ascending chains). Call $R$ a \emph{strict partial order} if $R$ is irreflexive, antisymmetric, and transitive. If $W$ is finite,
then $\mathfrak{F}$ is a $\mathsf{GL}$-frame iff $R$ is a strict partial order.

A characterization of descriptive $\mathsf{GL}$-frames was originally established by Esakia and given in \cite{Abashidze_1981}.
It generalizes directly to descriptive $\mathsf{MGL}$-frames. For a transitive frame $\mathfrak{F}=(W,R)$ and $A\subseteq W$, define
the \emph{irreflexive maximum} of $A$ by
\[
\qquad\mu(A)=\{x\in A \mid R(x)\cap A=\varnothing\}.
\]

\begin{lem} [\cite{Abashidze_1981}] \label{lem:GL max props}
Let $\mathfrak{F}=(W,R,E)$ be a descriptive augmented Kripke frame. Then $\mathfrak{F}$ is a descriptive $\mathsf{MGL}$-frame iff
$\mathfrak F$ is transitive and $A\subseteq\mu(A)\cup R^{-1}\mu(A)$ for each clopen $A$.
\end{lem}

Thus, a descriptive augmented Kripke frame is a descriptive $\mathsf{MGL}$-frame iff it is transitive and each point in a clopen set
is either in the irreflexive maximum of the clopen or sees a point in the irreflexive maximum.
It was observed by Japaridze \cite{Dzhaparidze_1988,Dzhaparidze_1990} that $\mathsf{MGL}$ has the finite model property.

\begin{theorem}[Japaridze]\label{MGL compl}
$\mathsf{MGL}\vdash\varphi\Leftrightarrow\mathfrak{F}\vDash\varphi$ for all finite $\mathsf{MGL}$-frames $\mathfrak{F}$.
\end{theorem}


\section{The G\"odel and splitting translations in the monadic setting}
\label{sec:3}

In this section we discuss the G\"odel and splitting translations in the monadic setting. While the G\"odel translation embeds $\mathsf{MIPC}$ fully and faithfully
into $\mathsf{MGrz}$, the splitting translation from $\mathsf{MGrz}$ into $\mathsf{MGL}$ does not yield a faithful embedding.

\subsection{G\"odel translation}

The G\"odel translation extends to the monadic setting by defining
\begin{eqnarray*}
    (\A\varphi)^t & = & \B\A \varphi^t\\
    (\E\varphi)^t & = & \E \varphi^t.
\end{eqnarray*}

Using algebraic semantics, Fisher-Servi \cite{Fischer-Servi_1977,Fischer-Servi_1978_1} proved that this provides a full and faithful embedding of $\mathsf{MIPC}$ into $\mathsf{MS4}$. The proof also yields a full and faithful embedding of $\mathsf{MIPC}$ into $\mathsf{MGrz}$.
Below we give an alternate proof of this result, using relational semantics.
The proof extends a semantical proof that $\mathsf{IPC}\vdash\varphi$ iff $\mathsf{S4}\vdash\varphi^t$
as given, e.g., in \cite[pp.~96-97]{Modal_Logic}.

For notational simplicity, we abbreviate the formula $\square\forall \psi$ as $\blacksquare \psi$ and the formula $\D\E \psi$ as $\bD \psi$.
Observe that this keeps the duality between box and diamond since $\blacksquare \psi=\neg\blacklozenge\neg \psi$ as
$\blacksquare \psi=\square\forall \psi=\neg\lozenge\neg\neg\exists\neg \psi,$ which is provably equivalent to
$\neg\lozenge\exists\neg \psi=\neg\blacklozenge\neg \psi$.

\begin{remark}
The modalities $\blacksquare,\blacklozenge$ are $\mathsf{S4}$-modalities which can be modeled using the relation $Q=R\circ E$, i.e., we have
\begin{flalign*}
& w \vDash\bB\varphi\quad\text{iff}\quad (\text{for all } v)(wQv \text{ implies } v\vDash\varphi);\\
& w \vDash\blacklozenge\varphi\quad\text{iff}\quad (\text{there exists } v)(wQv\text{ and } v\vDash\varphi).
\end{flalign*}
\end{remark}
Using this notation, the $\A$-step in the G\"odel translation is
\begin{eqnarray*}
    (\A\varphi)^t & = & \blacksquare \varphi^t.
\end{eqnarray*}

For an
$\mathsf{MS4}$-frame $\mathfrak{F}=(W, R, E)$ define an equivalence relation $\sim$ on $\mathfrak{F}$ by
\[
x\sim y\text{ iff } xRy \text{ and } yRx.
\]
Let $[x]$ denote the equivalence class of $x$, and let $W_{\sim}=W/{\sim}$ be the set of all equivalence classes.
Define $R_{\sim}$ and $E_{\sim}$ on $W_{\sim}$ by
\begin{align*}
[x] R_{\sim} [y] \quad &\text{iff} \quad xRy; \quad \\
[x] E_{\sim} [y] \quad &\text{iff} \quad  xQy \text{ and } yQx.
\end{align*}
That $E_\sim$ is well defined follows from $R \circ Q \circ R \subseteq Q$ which is true by commutativity in $\mathfrak F$ and transitivity of $R$.
Let $\mathfrak{F}_{\sim}=(W_{\sim}, R_{\sim}, E_{\sim})$.
Set $Q_\sim = E_\sim \circ R_\sim$.

\begin{lem} \label{lem:propertiesQsim}
Let $\mathfrak{F}=(W, R, E)$ be an
$\mathsf{MS4}$-frame and $x, y \in W$.
\begin{enumerate}
\item $x E y$ implies $[x] E_\sim [y]$;
\item $x Q y$ iff $[x] Q_\sim [y]$.
\end{enumerate}
\end{lem}

\begin{proof}
(1) If $x E y$, then $x Q y$ and $y Q x$, so $[x] E_\sim [y]$ by definition of $E_\sim$.

(2) Suppose that $x Q y$. Then there is $y'$ with $x R y'$ and $y' E y$. Therefore, $[x] R [y']$ by definition of $R$ and $[y'] E_\sim [y]$ by (1). Thus, $[x] Q_\sim [y]$. Conversely, if $[x] Q_\sim [y]$, then there is $[y']$ with $[x] R_\sim [y']$ and $[y'] E_\sim [x]$. By the definitions of $R_\sim$ and $E_\sim$, we have $x R y'$ and $y' Q y$. Thus, $x Q y$.
\end{proof}

\begin{lem}
$\mathfrak{F}_{\sim}$ is an
$\mathsf{MIPC}$-frame.
\end{lem}

\begin{proof}
It is well known (and easy to verify) that $R_{\sim}$ is a partial order (see, e.g., \cite[p.~68]{Modal_Logic}).
Transitivity and reflexivity of $E_\sim$ easily follow from transitivity and reflexivity of $Q$, and $E_\sim$ is symmetric by definition.
To see that $\mathfrak{F}_{\sim}$ satisfies commutativity, let $[x], [y], [z] \in W_\sim$ with $[x] E_\sim [y]$ and $[y] R_\sim [z]$.
Then $xQy$ and $yRz$, so $xQz$. Therefore, there is $z'$ with $x R z'$ and $z' E z$.
From $xRz'$ it follows that $[x] R [z']$, and $z' E z$ implies $[z'] E_\sim [z]$ by Lemma \ref{lem:propertiesQsim}(1). Thus, $\mathfrak F_\sim$ satisfies commutativity.
\end{proof}

Given a valuation $\nu$ on $\mathfrak F$, define a valuation $\nu_{\sim}$ on $\mathfrak F_{\sim}$ by
\[
\nu_{\sim}(p)=\{[x]\mid x\in \nu(\B p)\}.
\]
Clearly $\nu_{\sim}(p)$ is an upset. We call $\mathfrak F_{\sim}$ the \emph{skeleton} of $\mathfrak F$ and $(\mathfrak F_{\sim}, \nu_\sim)$ the
\emph{skeleton} of $(\mathfrak F, \nu)$.

Conversely, given an
$\mathsf{MIPC}$-frame $\mathfrak F$, we regard it as an
$\mathsf{MS4}$-frame. In addition, if $\mathfrak F$ is finite, then we regard it as a finite $\mathsf{MGrz}$-frame.
If $\nu$ is a valuation on the $\mathsf{MIPC}$-frame $\mathfrak F$, then we regard it as a valuation on the $\mathsf{MGrz}$-frame $\mathfrak F$.

The following lemma describes how the above frame transformations behave with respect to the G\"odel translation. It is proved by induction on the complexity of $\varphi$.

\begin{lem}\label{lem: translations}
Let $\varphi$ be a formula of $\mathcal{L}_{\forall\exists}$.
\begin{enumerate}
\item For an
$\mathsf{MIPC}$-frame $\mathfrak{F}$ with a valuation $\nu$ and $x \in \mathfrak F$ we have
\[
(\mathfrak F, x) \vDash_{\nu}\varphi \Leftrightarrow  (\mathfrak F, x) \vDash_{\nu}\varphi^t.
\]

\item For an
$\mathsf{MS4}$-frame $\mathfrak F$ with a valuation $\nu$ and $x \in \mathfrak F$, we have
\[
(\mathfrak F, x) \vDash_{\nu}\varphi^t \Leftrightarrow (\mathfrak F_{\sim}, [x]) \models_{\nu_{\sim}}\varphi.
\]
\end{enumerate}
\end{lem}

\begin{proof}
If $\mathfrak F$ is an $\mathsf{MIPC}$-frame, then $\mathfrak F_\sim$ is isomorphic to $\mathfrak F$. Therefore, (1) follows from (2).
To prove (2), by \cite[Lem.~3.81]{Modal_Logic}, it is sufficient to only consider the case for the modalities $\forall$ and $\exists$.
Let $\varphi=\A\psi$. Then
        \begin{flalign*}
        \quad[x] \vDash \A\psi \quad &\Leftrightarrow\quad (\mbox{for all } [y]) ([x]Q_\sim [y] \Rightarrow [y]\vDash \psi) &\\
        &\Leftrightarrow\quad (\mbox{for all } [y]) ([x]Q_\sim[y] \Rightarrow y\vDash \psi^t) \quad\text{(Inductive Hypothesis)} &\\
        &\Leftrightarrow\quad (\mbox{for all } y) (xQy \Rightarrow y \vDash \psi^t) \quad\text{(Lemma \ref{lem:propertiesQsim}(2))} &\\
        &\Leftrightarrow\quad x\vDash \blacksquare \psi^t &\\
        &\Leftrightarrow\quad x\vDash (\forall\psi)^t. &
        \end{flalign*}

Next let $\varphi=\E\psi$. If $x \models (\E\psi)^t $, then there is $y$ with $x Ey $ and $y \models \psi^t $. Therefore,
$[y] \models \psi$ by the inductive hypothesis, and $[x] E [y]$ by Lemma \ref{lem:propertiesQsim}(1). Thus, $[x] \models \E\psi$.
Conversely, suppose that $[x] \models \E\psi$. Then there is $[y]$ with $[x] E_\sim [y]$ and $[y] \models \psi$. Therefore,
$y Q x$ by the definition of $E_\sim$. Thus, there is $x'$ with $y R x'$ and $x' E x$. By the definition of $R_\sim$, we have $[y] R_\sim [x']$.
So $[x'] \models \psi$ by the persistence in $\mathfrak F_\sim$. Consequently, $x' \models \psi^t$ by the inductive hypothesis, and hence
$x \models  \E \psi^t = (\E\psi)^t$.
\end{proof}

\begin{theorem}
\label{MIPC MS4 MGrz frames}
$ \mathsf{MIPC}\vdash\varphi \text{ iff } \mathsf{MS4}\vdash \varphi^t \text{ iff } \mathsf{MGrz}\vdash \varphi^t$.
\end{theorem}

\begin{proof}
Suppose that $\mathsf{MIPC}\not\vdash\varphi$. Since $\mathsf{MIPC}$ has the FMP (Theorem \ref{thm:mipchasfmp}), there exists a finite
$\mathsf{MIPC}$-frame $\mathfrak F$, a valuation $\nu$ on $\mathfrak F$, and $x \in \mathfrak F$ such that $x \not \models_\nu \varphi$.
By regarding $\mathfrak F$ as an $\mathsf{MGrz}$-frame, $x \not \models_\nu \varphi^t$ by Lemma~\ref{lem: translations}(1). Therefore,
$\mathsf{MGrz}\not\vdash \varphi^t$. Also, as $\mathsf{MS4}\subseteq\mathsf{MGrz}$, it follows that $\mathsf{MS4}\not\vdash \varphi^t$.

Conversely, if $\mathsf{MGrz}\not\vdash \varphi^t$, then $\mathsf{MS4}\not\vdash \varphi^t$. By the FMP for $\mathsf{MS4}$, there is a
finite $\mathsf{MS4}$-frame $\mathfrak F$, a valuation $\nu$ on $\mathfrak F$, and $x \in \mathfrak F$ such that
$(\mathfrak F, x) \not \models_{\nu} \varphi^t$. By Lemma~\ref{lem: translations}(2), $(\mathfrak F_{\sim}, [x]) \not \models_{\nu_{\sim}}\varphi$.
Thus, $\mathsf{MIPC}\not\vdash\varphi$.
\end{proof}

\subsection{Splitting translation}

Next we discuss the splitting translation in the monadic setting.
The key here is Esakia's observation that the splitting translation does not yield a faithful embedding of $\mathsf{MGrz}$ into $\mathsf{MGL}$. Since this result is unpublished, we give a proof of it.

\begin{defn}
\label{def:cleancluster}
Let $\mathfrak{F}=(W,R,E)$ be an augmented Kripke frame (modal or intuitionistic), and let $x\in W$.
    \begin{enumerate}
        \item An $E$\emph{-cluster} (or \emph{cluster}) is a subset of $W$ of the form $E(x)=\{w\in W:xEw\}$ (it is the equivalence
        class of $x\in W$ with respect to $E$).
        \item We say that the $E$-cluster $E(x)$ is \emph{dirty} if there are $u,v\in E(x)$ with $u\not=v$ and $uRv$.
        \item We say that the cluster is \emph{clean} otherwise; that is, $u,v\in E(x)$ and $uRv$ imply $u=v$.
    \end{enumerate}
\end{defn}

\begin{minipage}{0.45\textwidth}
    \[\xymatrix{
        & \world{v} \\
        \world{v}\ar@{<->}[r]|-{E} & \world{u}\ar[u]|-{R}
    }\]
	\vspace{0.5em}
    \begin{center}
        Dirty cluster
    \end{center}
\end{minipage}
\hfill
\begin{minipage}{0.45\textwidth}
    \[\xymatrix{
        \world{v} \\
        \world{u}\ar[u]|-{R}
        \save"1,1"."2,1"*\frm<15pt,32pt>{e}\restore
    }\]
    \vspace{0.5em}
    \begin{center}
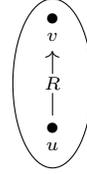
        Dirty cluster (alternate depiction - oval represents \emph{E}-cluster)
    \end{center}
\end{minipage}

\vspace{0.2in}

Descriptive $\mathsf{MGL}$-frames have the property that clusters in the irreflexive maximum of an $E$-saturated clopen are clean.

\begin{lem}\label{lem:mCasGL}
Let $\mathfrak{F}=(W,R,E)$ be a descriptive $\mathsf{MGL}$-frame. For clopen $A$ and $m\in \mu(E(A))$ we have that $E(m)$ is clean.
\end{lem}

\begin{proof}
Suppose there exist clopen $A$ and $m\in \mu (E(A))$ with $E(m)$ dirty. Then there are $x,y\in E(m)$ with $xRy$, $xEy$, and $x\not =y$.
By commutativity, there is $w$ such that $mRw$ and $wEy$, as shown below.
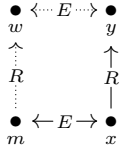
\[\xymatrix{
\world{w}\ar@{<.>}[r]|-E &\world{y}\\
\world{m}\ar@{.>}[u]|-{R}\ar@{<->}[r]|-{E} &\world{x}\ar[u]|-{R}
}
\]
Since $y\in E(A)$ we have $w\in E(A)$.
But this contradicts $R(m)\cap E(A)=\varnothing$. Thus, we cannot have a dirty cluster in $\mu (E(A))$.
\end{proof}

As a consequence of Lemma~\ref{lem:mCasGL}, we obtain:
	
\begin{lem}
\label{finite MGL frames}
Finite $\mathsf{MGL}$-frames are finite strict partial orders in which all clusters are clean.
\end{lem}

We next show that the splitting of the G\"odel translation of the monadic version of Casari's formula
\[
\mathsf{MCas: }\quad\forall\left(\left(p\rightarrow\forall p\right)\rightarrow\forall p\right)\rightarrow \forall p
\]
is provable in $\mathsf{MGL}$.

Since $\B\A p\leftrightarrow \B\A\B p$ is provable in $\mathsf{MS4}$, it is straightforward to check that $(\mathsf{MCas})^t$ is provably
equivalent to $\B\A\big(\B(\B p\rightarrow\B\A p)\rightarrow\B\A p\big)\rightarrow\B\A  p$. Using the notation $\blacksquare $ introduced above,
we have that $(\mathsf{MCas})^t$ is:
\[
\mathsf{M_\Box Cas}: \quad\bB\big(\B(\B p\rightarrow\bB p)\rightarrow\bB p\big)\rightarrow\bB p.
\]

Note that
$(\bB \psi)^s = (\B\A \psi)^s = \B^{+}\A \psi = \A \psi \wedge \B\A \psi = \A \psi \wedge \bB \psi$. So we can use $\bB^{+}\psi$ to abbreviate
$\A\psi\wedge\B\A\psi=\B^{+}\A\psi$, and so the splitting translation of $\mathsf{M_\Box Cas}$ is
\[
(\mathsf{M_\Box Cas})^s = \bB^{+}\big(\B^{+}(\B^{+}p\rightarrow\bB^{+} p)\rightarrow\bB^{+} p\big)\rightarrow\bB^{+} p.
\]

\begin{theorem}\label{MGL satisfies ST MCas}
$\mathsf{MGL}\vdash (\mathsf{M_\Box Cas})^s$.
\end{theorem}

\begin{proof}
Suppose $\mathfrak{F}=(W,R,E)$ is a descriptive $\mathsf{MGL}$-frame. We will prove that $\mathfrak{F}\vDash (\mathsf{M_\Box Cas})^s$.
Let $\nu$ be a valuation on $\mathfrak F$, $x \in \mathfrak F$, and $x\not\vDash_\nu\bB^{+}p$. We show that
$x\not\vDash\bB^{+}(\B^{+}(\B^{+}p\rightarrow\bB^{+}p)\rightarrow\bB^{+}p)$. Let $A=W\setminus \nu(\bB^{+}p)$. Then $x \in A$ and so by
Lemma \ref{lem:GL max props}, $x\in\mu(A)\cup R^{-1}\mu(A)$. If $x \in R^{-1}\mu(A)$, then there is $x' \in  \mu(A)$ with $x R x'$.
If $x \in \mu(A)$, we let $x' = x$. From $x' \in \mu(A)$ it follows that $x' \in A$, so $x' \not \models \bB^{+}p= \forall p \wedge \Box \forall p$.
We show that $x' \not \models \forall p$. If $x' \not \models  \Box \forall p$, then there is $y$ with $x'Ry$ and $y \not \models \A p$. Therefore,
$y \not \models \bB^{+} p$, so $y \in A$. But this contradicts $x' \in \mu (A)$. Thus, $x' \not \models \forall p$. So there is $w$ with $w E x'$
and $w \not \models p$. We show that $ w \not \models \B^{+}(\B^{+}p\rightarrow\bB^{+} p) \rightarrow \bB^{+} p$.

Since $w\not\vDash p$, we have $w\not\vDash p\wedge\B p$, so $w\not\vDash\B^{+} p$, and hence $w\vDash \B^{+}p\rightarrow\bB^{+}p$. Let
$w R z$.
By commutativity, there is $y$ such that $x'Ry$ and $yEz$. Since $x' \in \mu(A)$, we have $y \not \in A$. Therefore, $y\vDash\bB^{+}p$,
so $y\vDash\A p$, and hence $z\vDash\A p$.

\[
\xymatrix{
\vDash\bB^{+} p \hspace{-1.5em} & \world{y}\ar@{<.>}[r]|-{E} & \world{z} & \hspace{-1.75em}\vDash\forall p \\
\hspace{-1.5em} & \world{x'}\ar@{<->}[r]|-{E}\ar@{.>}[u]|-{R} & \world{w}\ar[u]|-{R} & \hspace{-2.25em}\not\vDash{p}
}
\]
In fact, if $zRt$, then $wRt$ by transitivity, and so by the same reasoning as above we have $t\vDash\A p$. It follows that $z\vDash \Box \A p$, and so $z\vDash \bB^{+}p$. Thus,
$z\vDash\B^{+}p\rightarrow\bB^{+}p$, and hence $w\vDash\B(\B^{+}p\rightarrow\bB^{+}p)$. This together with $w\vDash \B^{+}p\rightarrow\bB^{+}p$
yields $w\vDash\B^{+}(\B^{+}p\rightarrow\bB^{+} p)$. Since $w\not\vDash\bB^{+} p$, we obtain
$w \not \vDash\B^{+}(\B^{+}p\rightarrow\bB^{+} p) \rightarrow \bB^{+} p$.

If $x = x'$, then $x E w$, and so $x \not \vDash \forall (\B^{+}(\B^{+}p\rightarrow\bB^{+} p) \rightarrow \bB^{+} p)$. Otherwise, $x R x'$ and
$x' Ew$ imply $x Q w$, so $x \not \models  \bB (\B^{+}(\B^{+}p\rightarrow\bB^{+} p) \rightarrow \bB^{+} p)$. Thus, in either case,
$x \not \models  \bB^+ (\B^{+}(\B^{+}p\rightarrow\bB^{+} p) \rightarrow \bB^{+} p)$ as desired.
This yields $x\vDash (\mathsf{M_\Box Cas})^s$. Since $x$ was arbitrary,
$\mathfrak{F}\vDash (\mathsf{M_\Box Cas})^s$. Because $\mathfrak F$ is an arbitrary descriptive $\mathsf{MGL}$-frame, by Theorem~\ref{thm: completeness 4}(2), $\mathsf{MGL}\vdash (\mathsf{M_\Box Cas})^s$.
\end{proof}

\begin{theorem}\label{MGrz fails TMCas}
    $\mathsf{MGrz}\not\vdash \mathsf{M_\Box Cas}$.
\end{theorem}

\begin{proof}
Consider the $\mathsf{MGrz}$-frame $\mathfrak{F}=(W,R,E)$ where $W=\{x,y\}$, $R=\{(x,x),(y,y),(x,y)\}$, and
$E=W^{2}=\{(x,x),(y,y),(x,y),(y,x)\}$, as shown below.
\[\xymatrix{
\world{y} & \hspace{-1.5em}\vDash p\\
\world{x}\ar[u]|-{R} & \hspace{-1.5em}\not\vDash p
\save"1,1"."2,1"*\frm<15pt,32pt>{e}\restore
}\]  
The arrow represents the nontrivial $R$-relation and the circle represents that both points are in the same $E$-equivalence class.
It is easy to see that this is an $\mathsf{MGrz}$-frame. Let $\nu$ be a valuation on $\mathfrak{F}$ with $\nu(p)=\{y\}$.

First, we claim that $x\vDash\bB(\B(\B p\rightarrow\bB p)\rightarrow\bB p)$. To see this, note that both $x\not\vDash\bB p$ and $y\not\vDash\bB p$,
but since $y\vDash p$ and $y$ only sees itself (with respect to $R$), we have $y\vDash\B p$. Thus, $y\not\vDash\B p\rightarrow\bB p$, so
$x\not\vDash\B(\B p\rightarrow \bB p)$, and hence $x\vDash\B(\B p\rightarrow\bB p)\rightarrow\bB p$. Moreover, $y\not\vDash\B(\B p\rightarrow \bB p)$, so $y\vDash\B(\B p\rightarrow\bB p)\rightarrow\bB p$, and hence
$x\vDash\bB(\B(\B p\rightarrow\bB p)\rightarrow\bB p)$. However, $x\not\vDash\bB p$ as $xQx$ and $x\not\vDash p$. Thus,
$x\not\vDash\bB(\B(\B p\rightarrow\bB p)\rightarrow\bB p)\rightarrow\bB p$, hence $\mathfrak{F}\not\vDash_{\mathsf{MGrz}}\mathsf{M_\Box Cas}$,
and so $\mathsf{MGrz}\not\vdash \mathsf{M_\Box Cas}$.
\end{proof}

\begin{cor} (Esakia)
The splitting translation does not embed $\mathsf{MGrz}$ into $\mathsf{MGL}$ faithfully.
\end{cor}


\section{The logics $\mathsf{M^{+}IPC}$ and $\mathsf{M^+Grz}$}
\label{sec:4}

In the previous section we saw that the splitting translation does not embed $\mathsf{MGrz}$ into $\mathsf{MGL}$ faithfully.
In fact, while the G\"odel translation of $\mathsf{MCas}$ is not provable in $\mathsf{MGrz}$, the splitting translation of
the G\"odel translation of $\mathsf{MCas}$ is provable in $\mathsf{MGL}$. Esakia suggested to strengthen $\mathsf{MIPC}$ with
$\mathsf{MCas}$ and $\mathsf{MGrz}$ with the G\"odel translation of $\mathsf{MCas}$, and see whether this repairs the disbalance.
This is what we do in this section.

\subsection{$\mathsf{M^{+}IPC}$}

\begin{defn}
The logic $\mathsf{M^{+}IPC}$ is defined as the extension of $\mathsf{MIPC}$ by $\mathsf{MCas}$:
\[
\mathsf{M^{+}IPC}=\mathsf{MIPC}+\mathsf{MCas}.
\]
\end{defn}

Recall from Definition \ref{def:cleancluster} that a cluster of an $\mathsf{MIPC}$-frame is called clean if no distinct points in the
cluster are $R$-related. The following semantic characterization of $\mathsf{M^+IPC}$-frames was established by Esakia. For a proof see
\cite[Lem.~38]{Bezhanishvili_2000}. It states that a descriptive $\mathsf{MIPC}$-frame is a descriptive $\mathsf{M^+IPC}$-frame iff
the cluster of each point in the $R$-maximum of the $E$-saturation of clopen is clean.

\begin{lem}\cite[Lem.~38]{Bezhanishvili_2000}\label{cas condition}
Let $\mathfrak{F}=(W,R,E)$ be a descriptive $\mathsf{MIPC}$-frame. Then $\mathfrak{F}\vDash\mathsf{MCas}$ iff for each clopen $A$, if
$m\in\max E(A)$, then $E(m)$ is clean.
\end{lem}

\begin{remark}
The condition in \cite[Lem.~38]{Bezhanishvili_2000} is that $\mathfrak{F}\vDash\mathsf{MCas}$ iff for each clopen $A$ we have
$A\subseteq Q^{-1}(\max A\cap \max Q^{-1}A)$. But, as discussed after the proof of \cite[Lem.~38]{Bezhanishvili_2000},
this statement is equivalent to the statement in Lemma~\ref{cas condition}.
\end{remark}

As a consequence of Lemma \ref{cas condition}, we obtain:

\begin{lem}\label{finite M+IPC}
Finite $\mathsf{M^{+}IPC}$-frames are finite $\mathsf{MIPC}$-frames in which all clusters are clean.
\end{lem}

\subsection{$\mathsf{M^{+}Grz}$}

\begin{defn}
The logic $\mathsf{M^{+}Grz}$ is the extension of $\mathsf{MGrz}$ by $\mathsf{M_\Box Cas}$:
\[
\mathsf{M^{+}Grz}=\mathsf{MGrz}+\mathsf{M_\Box Cas}.
\]
\end{defn}

\begin{remark}
As we pointed out in the previous section, $\mathsf{M_\Box Cas}$ is provably equivalent to the G\"odel translation of $\mathsf{MCas}$.
\end{remark}

In order to obtain a semantic characterization of $\mathsf{M^+Grz}$, which is an analogue of Lemma \ref{cas condition},
we require the following lemma.

\begin{lem}
\label{lem:maxgrzclean}
Let $\mathfrak{F}=(W,R,E)$ be a descriptive $\mathsf{MGrz}$-frame, $A \subseteq W$ clopen, and $y \in \max E(A)$.
If $E(y)$ is clean, then:
\begin{enumerate}
\item $E(y) \subseteq \max E(A)$;
\item for all $z \in W$, from $yR z$ and $z R y$ it follows that $y = z$.
\end{enumerate}
\end{lem}

\begin{proof}
(1)  Let $z  \in E(y)$ and $w \in E(A)$ with $ z R w$. By commutativity, there is $w'$ with $y R w'$ and $w' E w$. Therefore, $w' \in E(A)$. Since $y \in \max E(A)$, we have $y=w'$. Thus, $z, w \in E(y)$ and $z R w$. As $E(y)$ is clean, $z = w$. This shows that $z \in  \max E(A)$.

(2) Suppose $yR z$ and $z R y$. From $y \in E (A)$ and $yRy$, we have $y \in R^{-1}E(A)$.
We show that $y \in \qmax R^{-1}E(A)$.
Let $yRw$ and $w\in R^{-1}E(A)$, so $wRu$ for some $u\in E(A)$. Then $yRu$ by transitivity,
and $y\in\max E(A)$ implies $y=u$, hence $wRy$, and so $y\in \qmax R^{-1}E(A)$. By Lemma \ref{lem:grzchar}(1), this means $y\in \max R^{-1}E(A)$.
Since $zR y$, we have $z \in R^{-1}E(A)$, so $yRz$ implies $z=y$.
\end{proof}

We now have the necessary machinery to prove a semantic characterization of $\mathsf{M^+Grz}$, which states that a descriptive $\mathsf{MGrz}$-frame is a descriptive $\mathsf{M^+Grz}$-frame iff the cluster of every point in the
maximum of the $E$-saturation of a clopen set is clean.

\begin{lem}\label{lem:mCasGrz}
Let $\mathfrak{F}=(W,R,E)$ be a descriptive $\mathsf{MGrz}$-frame. Then ${\mathfrak{F}\vDash \mathsf{M_\Box Cas}}$ iff
for each clopen $A$ and $m\in \max E(A)$ we have $E(m)$ is clean.
\end{lem}

\begin{proof}
First suppose $\mathfrak{F}\not\vDash \mathsf{M_\Box Cas}$. Then there is $x\in W$ such that
$x\not\vDash\bB(\B(\B p\rightarrow\bB p)\rightarrow\bB p)\rightarrow\bB p$, and hence
$x\vDash\bB(\B(\B p\rightarrow\bB p)\rightarrow\bB p)$ but $x\not\vDash\bB p$. Since $x\not\vDash\bB p$,
there is $x' \in W$ such that $xQx'$ and $x'\not\vDash p$. Let $A=\{w\in W \mid w\not\vDash p\}$. Then $x'\in A$,
and as $x'Ex'$, we have $x' \in E(A)$. Because $A$ is clopen, so is $E(A)$.
By Lemma~\ref{lem:grzchar}(2), there is $y \in\max E(A)$ with $x'R y$. If $E(y)$ is dirty, then we are done. So assume that $E(y)$ is clean.
We show that this leads to a contradiction.
Since $y \in E(A)$, there is $y' \in A$ with $y E y'$. By Lemma \ref{lem:maxgrzclean}(1), $y' \in \max E(A)$. Because $x Q y'$ and $x\vDash\bB(\B(\B p\rightarrow\bB p)\rightarrow\bB p)$, we have $y'\vDash\B(\B p\rightarrow\bB p)\rightarrow\bB p$. As $y' \in A$, we have $y'\not\vDash p$ and since $y'Qy'$, we have
$y'\not\vDash\bB p$, so we must have $y'\not\vDash\B(\B p\rightarrow\bB p)$. Thus, there is $z\in W$ such that $y'Rz$ and
$z\not\vDash\B p\rightarrow\bB p$, which means $z\vDash\B p$ but $z\not\vDash\bB p$. Because $z\not\vDash\bB p$, there exist $w',w\in W$ such that $zRw'Ew$ and $w\not\vDash p$ (see the diagram below).

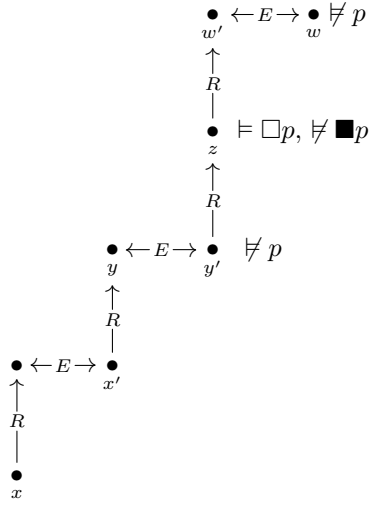
\[\xymatrix{
 &  & \world{w'}\ar@{<->}[r]|-{E} & \world{w}\not\vDash p \\
 &  & \world{z}\ar[u]|-{R} & \hspace{-2.5em}\vDash\B p, \,\not\vDash\bB p \\
 & \world{y}\ar@{<->}[r]|-{E} & \world{y'}\ar[u]|-{R} & \hspace{-5.5em}\not\vDash{p} \\
\bullet\ar@{<->}[r]|-{E} & \world{x'}\ar[u]|-{R} &  &  & \\
\world{x}\ar[u]|-{R} &  &  &
}\]

Now, since $w\not\vDash p$, we have $w\in A$ and hence $w'\in E(A)$. Thus, $y'Rw'$ and $w'\in E(A)$,
so by $R$-maximality of $y'$ in $E(A)$, we must have $y'=w'$.
But then $y'R z$ and $z R y'$, and so by Lemma~\ref{lem:maxgrzclean}(2), $y'=z$. This, however, is a contradiction since
$z\vDash\B p$, hence $z\vDash p$, whereas $y'\not\vDash p$.

For the converse, suppose that $A$ is clopen and $m\in\max E(A)$ with $E(m)$ dirty.
First observe that since $m$ is maximal in $E(A)$, from $m Q t$ it follows that $t \in E(m)$ for all $t \in E(A)$. Indeed, if $m Q t$ for
$t \in E(A)$, then there is $t'$ with $m R t'$ and $t' E t$. Since $t' \in E(A)$, we have $t' =m$ by maximality of $m$ in $E(A)$. Thus,
$t \in E(m)$.

Now, since $E(m)$ is dirty, there are $x , x'  \in E(m)$ with $x R x'$ and $x \neq x'$.
In particular, $x \not \in \max E(A)$. Since $E(A)$ is clopen, $\max E(A)$ is closed (see, e.g., \cite[Sec.~III.2]{Esakia_1985}).
Thus, we can find clopen $B$ such that $x\in B$ but $B\cap\max E(A)=\varnothing$, as shown below.

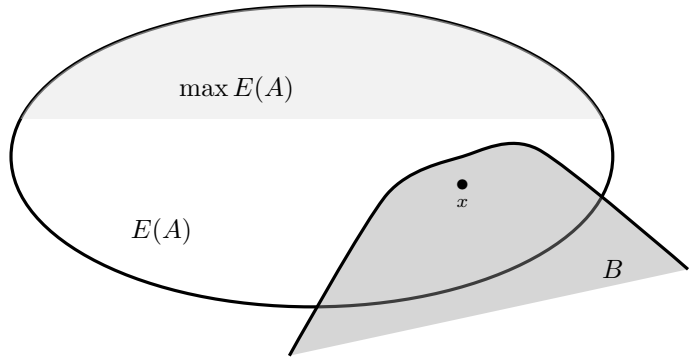
\begin{center}
\begin{minipage}{0.6\textwidth}
\vspace{1pt}
\begin{tikzpicture}
\begin{scope}
\draw [very thick] (4,1.5) ellipse (4cm and 2cm);
\clip (4,1.5) ellipse (4cm and 2cm);
\fill [gray!20, fill opacity=0.5] (0,2) rectangle (8,4);
\end{scope}
\draw [fill=gray!80, very thick, fill opacity=0.4] plot [smooth, tension=0.45] coordinates {(3.7, -1.15) (5, 1) (6, 1.5) (7, 1.6) (9, 0)};
\node (a) at (2,0.5) {$E(A)$};
\node (b) at (6,1) {$\underset{x}{\bullet}$};
\node (c) at (3,2.4) {$ \max E(A)$};
\node (d) at (8,0) {$B$};
\end{tikzpicture}
\end{minipage}
\end{center}

Choose a valuation $\nu$ with $\nu(p)=W\setminus(B\cap E(A))$. Note that
$\nu$ is well-defined as $B$ and $E(A)$ are clopen. We aim to show that
$x\vDash \bB(\B(\B p\rightarrow\bB p)\rightarrow\bB p)$ but $x\not\vDash\bB p$. Since $x \in B\cap E(A)$, we have $x \not \models p$.
This implies that $x \not \models \bB p$ because $x Q x$. To finish the argument it suffices to show that
$y \models  \B(\B p\rightarrow\bB p)\rightarrow\bB p$ for all $y$ with $x Q y$.
So let $x Q y$ and assume that $y\not \models \bB p$. Then there is $z$ with $yQ z$ and $z \not \models p$. Therefore, $z \in  B\cap E(A)$
and there is $z'$ with $y R z'$ and $z' E z$. Clearly $z' \in E(A)$. By Lemma~\ref{lem:grzchar}(2), there is $t \in\max E(A)$ with $z'R t$.

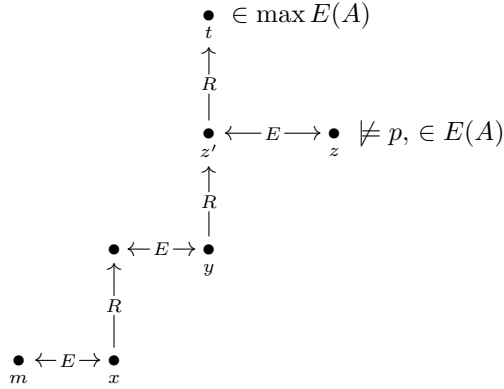
\[\xymatrix{
 & & \world{t}& \hspace{-2.5em}   \in \max E(A) & \\
 & & \world{z'}\ar[u]|-{R}\ar@{<->}[r]|-{E}& \world{z} & \hspace{-3.5em}  \not \models p, \, \in E(A)\\
 & \bullet\ar@{<->}[r]|-{E}& \world{y}\ar[u]|-{R} & &  \hspace{-5.5em}  \\
\world{m}\ar@{<->}[r]|-{E} & \world{x}\ar[u]|-{R} &  & &
}\]

Since $t \in \max E(A)$, we have $t \not \in B$, so $t \models p$ and if $t R v$ for $t \neq v$, then $v \not \in E(A)$ by maximality of $t$,
so $v \models p$. Thus, $t \models \Box p$. On the other hand, $x Q y$, $y R z' $, and $z' R t$ imply $m Q t$. As we saw above, this means $t \in E(m)$, and so $t Ex$. Since $x \not \models p$, we have $t \not \models \blacksquare p$. This implies that
$t \not \models \B p\rightarrow\bB p$, so $y \not \models \B(  \B p\rightarrow\bB p )$, and hence
$y \models  \B(\B p\rightarrow\bB p)\rightarrow\bB p$ as desired.
\end{proof}

As a consequence of Lemma \ref{lem:mCasGrz}, we obtain:

\begin{lem}\label{finite MGrz frames}
Finite $\mathsf{M^{+}Grz}$-frames are finite $\mathsf{MGrz}$-frames in which all clusters are clean.
\end{lem}

\subsection{The translations $\mathsf{M^{+}IPC}\to\mathsf{M^{+}Grz}\to\mathsf{MGL}$}

As we pointed out, the remaining part of the paper establishes the finite model property for the logics $\mathsf{M^+IPC}$ and $\mathsf{M^+Grz}$.
We finish this section by explaining how a proof of Esakia's claim is then obtained.


Let $R$ be a binary relation. We recall that the \emph{irreflexive reduction} of $R$, denoted $R^{i}$, is defined by
\[
aR^{i}b \text{ iff } aRb \text{ and } a\not=b;
\]
and the \emph{reflexive closure} of $R$, denoted $R^{r}$, is defined by
\[
aR^{r}b \text{ iff } aRb \text{ or } a=b.
\]

For an augmented Kripke frame $\mathfrak{F}=( W,R, E)$, let $\mathfrak{F}^i=( W,R^i, E)$ and $\mathfrak{F}^{r}=( W,R^{r}, E)$.
Following the terminology of \cite[pp. 98-99]{Modal_Logic}, we call $\mathfrak{F}^i$ the \emph{irreflexive reduction}
and $\mathfrak{F}^{r}$ the \emph{reflexive closure} of $\mathfrak{F}$.

\begin{lem}\label{lem: i & r}
\begin{enumerate}
\item[]
\item If $\mathfrak F$ is a finite $\mathsf{M^+Grz}$-frame, then $\mathfrak{F}^i$ is a finite $\mathsf{MGL}$-frame.
\item If $\mathfrak F$ is a finite $\mathsf{MGL}$-frame, then $\mathfrak{F}^r$ is a finite $\mathsf{M^+Grz}$-frame.
\end{enumerate}
\end{lem}

\begin{proof}
Since finite $\mathsf{M^+Grz}$-frames are finite partial orders with clean clusters (Lemma~\ref{finite MGrz frames}) and finite $\mathsf{MGL}$-frames are finite
strict partial orders with clean clusters (Lemma~\ref{finite MGL frames}), this is an immediate consequence of \cite[pp. 98-99]{Modal_Logic}.
\end{proof}

\begin{lem}\label{lem:withplustranslations}
Let $\varphi$ be a formula of $\mathcal{ML}_{\A}$.
\begin{enumerate}
\item For a finite $\mathsf{M^{+}Grz}$-frame $\mathfrak{F}$, a valuation $\nu$ on $\mathfrak F$, and $x \in \mathfrak F$ we have
\[
(\mathfrak F, x) \vDash_{\nu}\varphi \Leftrightarrow (\mathfrak F^i, x) \vDash_{\nu}\varphi^s.
\]
\item For a finite $\mathsf{MGL}$-frame $\mathfrak F$, a valuation $\nu$ on $\mathfrak F$, and $x \in \mathfrak F$ we have
\[
(\mathfrak F, x) \vDash_{\nu}\varphi^s \Leftrightarrow (\mathfrak F^r, x) \models_{\nu}\varphi.
\]
\end{enumerate}
\end{lem}

\begin{proof}
The proof is an immediate consequence of \cite[pp. 98-99]{Modal_Logic} since the quantifier modalities are not changed by the translation $(-)^s$, nor is the relation $E$ altered
going from $\mathfrak F$ to $\mathfrak F^i$ or $\mathfrak F^r$.
\end{proof}

Finally, we are ready to provide a proof of Esakia's claim.

\begin{theorem}
$\mathsf{M^{+}IPC}\vdash\varphi \text{ iff }\mathsf{M^{+}Grz}\vdash \varphi^t \text{ iff } \mathsf{MGL}\vdash (\varphi^t)^s$.
\end{theorem}

\begin{proof}
The first equivalence is proved exactly as Theorem \ref{MIPC MS4 MGrz frames} using the fact that finite $\mathsf{M^+IPC}$-frames and finite
$\mathsf{M^+Grz}$-frames coincide.

For the second equivalence, suppose $\mathsf{MGL}\not\vdash (\varphi^t)^s$. Since $\mathsf{MGL}$ has the FMP, there exist a finite
$\mathsf{MGL}$-frame $\mathfrak{F}$, a valuation $\nu$ on $\mathfrak F$, and $x \in \mathfrak F$ such that
$(\mathfrak F, x)\not\vDash_{\nu}(\varphi^t)^s$. By Lemma~\ref{lem:withplustranslations}(2), $(\mathfrak{F}^{r}, x) \not\vDash_{\nu}\varphi^t$,
and since $\mathfrak F^r$ is an $\mathsf{M^+Grz}$-frame by Lemma~\ref{lem: i & r}(2),
$\mathsf{M^{+}Grz}\not\vdash \varphi^t$. For the converse, suppose $\mathsf{M^{+}Grz}\not\vdash \varphi^t$. Since $\mathsf{M^{+}Grz}$ has the FMP,
there exist a finite $\mathsf{M^{+}Grz}$-frame $\mathfrak{F}$, a valuation $\nu$ on $\mathfrak F$, and $x \in \mathfrak F$ such that
$(\mathfrak{F},x)\not \models_{\nu}\varphi^t$.
By Lemma~\ref{lem:withplustranslations}(1), $(\mathfrak{F}^{i}, x) \not\vDash_{\nu}(\varphi^t)^s$, and since $\mathfrak F^i$ is an
$\mathsf{MGL}$-frame by Lemma~\ref{lem: i & r}(1), we conclude that
$\mathsf{MGL}\not\vdash (\varphi^t)^s$.
\end{proof}

We now have succeeded in lifting the original correspondences given by Goldblatt, Boolos, Kuznetsov and Muravitsky from the propositional
setting to the monadic setting, verifying Esakia's claim.
Combining this with Japaridze's result of arithmetical completeness for $\mathsf{MGL}$ yields
arithmetic interpretations of $\mathsf{M^{+}IPC}$ and $\mathsf{M^{+}Grz}$.


\section{The finite model property of $\mathsf{M^{+}IPC}$}
\label{sec:5}

This section is dedicated to the proof of the finite model property of $\mathsf{M^+IPC}$. We do this by modifying the selective filtration
technique originally developed by Grefe \cite{Gre98} to prove the finite model property of Fisher Servi's intuitionistic modal logic $\mathsf{FS}$. In \cite[\S 10.3]{Many_Dim_Modal_Logics} it was used to give an alternative proof of the finite model property of $\mathsf{MIPC}$.

We start by collecting some properties of descriptive $\mathsf{M^+IPC}$-frames that will be useful in what follows. The following lemma is the $\mathsf{M^+IPC}$-version of Lemma~\ref{lem:maxgrzclean}(1).

\begin{lem}
\label{lem:cleanclustersmipc}
Let $\mathfrak F= (W, R, E) $ be a descriptive $\mathsf{M^+IPC}$-frame, $A \subseteq W$ clopen, $y \in \max E(A)$, and $E(y)$ clean. Then $E(y) \subseteq \max E(A)$.
\end{lem}

\begin{proof}
If $E(y) \not\subseteq \max E(A)$, then there are distinct $t \in E(y)$ and $u \in E(A)$ with $t R u$.
By commutativity, there is $u'$ with $y R u'$ and $u' E u$. Therefore, $u' \in E(A)$, so by maximality of $y$ in $E(A)$ we have $y = u'$. This implies that $t E u$, contradicting that $E(y)$ is a clean cluster.
\end{proof}

We say a point $x$ is \emph{maximal with respect to} a formula $\psi$ if $x\not\vDash\psi$ and for each $y$ with $xRy$ and $x\neq y$ we
have $y\vDash\psi$ (that is, $x$ refutes $\psi $ and every point strictly above $x$ validates $\psi$).

\begin{lem}\label{lem:find maximal points new}
Let $\mathfrak F= (W, R, E) $ be a descriptive $\mathsf{M^+IPC}$-frame, $t \in W$, and $\nu$ a valuation on $\mathfrak F$.
\begin{enumerate}
\item Let $A \subseteq W$ be clopen. If $t \in E(A)$, then there is $x \in \max E(A)$ such that $tRx$ and $E(x)$ is clean.
\item If $t \not \models \forall \varphi$, then there is $x$ such that $t R x$, $x$ is maximal with respect to $\forall \varphi$,
and $E(x)$ is clean.
\item Let $A \subseteq W$ be clopen. If $t \in  A$, then there is $x \in A \cap \max E(A)$ such that $tQx$ and $E(x)$ is clean.
\item If $t \not \models \varphi$, then there is $x$ such that $t Q  x $, $x$ is maximal with respect to $\varphi$, and  $E(x)$ is clean.
\end{enumerate}
\end{lem}

\begin{proof}
(1) Let $t \in E(A)$.
By Lemma~\ref{MIPC points can see max}, there is $x \in \max E(A)$ such that $tRx$. By Lemma~\ref{cas condition}, $E(x)$ is clean.

(2) Suppose that $t \not \models \forall \varphi$. Let $A = \nu ( \forall \varphi)^c$. Then $A$ is clopen, $E(A) = A$, and $t \in E(A)$.
By (1), there is $x \in \max E(A)$ such that $tRx$ and $E(x)$ is clean. Since $E(A) =A$, it immediately follows that $x$ is maximal with
respect to $\forall\varphi$.

(3) Let $t \in A$. Then $t \in E(A)$. By (1), there is $x' \in \max E(A)$ such that $tRx'$ and $E(x')$ is clean. Since $x' \in E(A)$, there is $x \in A$ with $x' E x$. Therefore, $t Q x$, and because $E(x')$ is clean, we have that $x \in \max E(A)$ by Lemma \ref{lem:cleanclustersmipc}.

(4) Suppose that $t \not \models \varphi$. Let $A = \nu (\varphi)^c$. Then $A$ is clopen and $t \in A$. By (3), there is $x \in A \cap \max E(A)$ such that $t Q x$ and $E(x)$ is clean. Since $x \in A$, we also have $x \in \max A$.  But the latter means that $x$ is maximal with respect to $\varphi$. Thus, $x$ is as desired.
\end{proof}

\subsection{The construction}

We start with a formula $\varphi$, a descriptive $\mathsf{M^+IPC}$-frame $\mathfrak F=(W, R, E)$, and a valuation $\nu$ on $\mathfrak F$ such that
$\mathfrak F \not \models \varphi$. By modifying the construction in \cite[\S 10.3]{Many_Dim_Modal_Logics}, we will construct a sequence of finite
$\mathsf{M^{+}IPC}$-frames $\mathfrak{F}_{h}=(W_{h},R_{h},E_{h})$ such that $\mathfrak{F}_{h}\subseteq\mathfrak{F}_{h+1}$ for all $h<\omega$. For each point $t\in W_{h}$ that we select, we will be creating a copy of some original point in $W$. We give each added point a new name, say $t$, and let $\widehat{t}$ denote the original point in $W$ that $t$ was copied from and will behave similar to. Thus, it is possible to have two different points $x_{1}$ and $x_{2}$ in our new frame, where $\widehat{x_{1}}=\widehat{x_{2}}$. The main difference between our construction and the construction given in \cite[\S 10.3]{Many_Dim_Modal_Logics} will be seen in the $\rightarrow$-step,
which requires a more careful selection of new points.

To start the construction, let $\mathfrak{F}_{0}=(W_{0},R_{0},E_{0})$ where
\[
W_{0}=\{t_{0}\},\quad R_{0}=W_{0}^{2}, \quad E_{0}=W_{0}^{2},
\]
and $\widehat{t_{0}}$ is a point in $W$ such that $\widehat{t_0}$ is from a clean cluster
and is maximal with respect to $\varphi$. The existence of such $\widehat{t_0}$  follows from Lemma~\ref{lem:find maximal points new}(4). Moreover, let $W^{\forall H}_{-1}= \varnothing$.

Let $\Sub$ be the set of subformulas of $\varphi$, and let $(W', R', E')$ be any of our frames in the construction. To each $t \in W'$ we associate the following subsets of $\Sub$:
\begin{flalign*}
& &\Sigma^{\exists}(t) &= \{\exists\delta\in \Sub:\widehat{t}\vDash\exists\delta\text\}&\\
& &\Sigma^{\forall H}(t) &=\{\forall\beta\in \Sub:\widehat{t}\text{ is maximal wrt }\forall\beta\}&\\
& &\Sigma^{\forall V}(t) &= \{\forall\gamma\in \Sub:\widehat{t}\not\vDash\forall\gamma\text{ but is not maximal wrt }\forall\gamma\}&\\
& &\Sigma^{\rightarrow}(t) &= \{\alpha\rightarrow\sigma\in \Sub:\widehat{t}\not\vDash\alpha\rightarrow\sigma \text{ but is not maximal wrt } \alpha \rightarrow \sigma\}. &
\end{flalign*}
These are precisely the subformulas of $\varphi$ whose truth-value at $\widehat{t}$ is relevant for constructing our countermodel.

Suppose $\mathfrak{F}_{h-1}=(W_{h-1},R_{h-1},E_{h-1})$ has already been constructed so that $\mathfrak{F}_{h-1}$ is a finite
$\mathsf{M^+IPC}$-frame and $E(\widehat{w})$ is a clean cluster for each $w\in W_{h-1}$. We construct $\mathfrak{F}_{h}$ applying the four steps described below. They are designed to add the necessary witnesses required by the formulas in the sets $\Sigma^{\exists}(t)$, $\Sigma^{\forall H}(t)$, $\Sigma^{\forall V}(t)$, and $\Sigma^{\rightarrow}(t)$, respectively. In the $\exists$-step we ensure that for each formula in $\Sigma^{\exists}(t)$ the point $t$ has an $E$-successor that witnesses the existential statement. In the $\forall H$-step we ensure that for each formula in $\Sigma^{\forall H}(t)$ the point $t$ has an $E$-successor that witnesses the refutation of the universal statement. In the vertical steps  $\forall V$ and  $\rightarrow$ we make sure that $t$ has the necessary $R$-successors that are maximal with respect to the formulas in $\Sigma^{\forall V}(t)$ and  $\Sigma^{\rightarrow}(t)$, respectively.
In each step of the construction we add also points to witness commutativity.
Note that the first three of the following four steps are only done once per cluster. This is enough since all points of a cluster in $\mathfrak F$ agree on refuting an $\forall$- or $\exists$-formula and points from a clean cluster agree whether such a refutation is maximal.

Roughly speaking, points are added to the construction in the following order:
In the first round the cluster of the starting point $t_{0}$ is built by adding points for formulas in $\Sigma^{\exists}(t_0)$ and $\Sigma^{\forall H}(t_0)$. After this, no more points are added to this cluster. We call this the `bottom cluster' of our frame. The first round of the construction proceeds by adding vertical witnesses for each formula in $\Sigma^{\forall V}(t_0)$ and closing each such cluster by adding points for commutativity. The first round then finishes by adding for each point $t$ in the `bottom cluster' vertical witnesses for the formulas in $\Sigma^{\rightarrow}(t)$ and closing under commutativity. In the next round all these newly build clusters will possibly be enlarged in the horizontal steps and then new vertical clusters will be added in the $\A V$- and $\rightarrow$-steps.

\vspace{1em}

\emph{$\exists$-step $($Horizontal$)$}:
Let $W^\exists_{h} = W_{h-1}$,  $R^\exists_{h} = R_{h-1}$, and $E^\exists_{h} = E_{h-1}$. For each $E^{\exists}_{h}(t)\subseteq W^{\exists}_{h}\backslash W_{h-1}^{\A H}$,
if $\exists\delta\in\Sigma^{\exists}(t)$ but there is no $s\in W^\exists_{h}$ already such that $tE^{\exists}_{h}s$ and $\widehat{s}\vDash\delta$, then we
add a point $s$ to $W^\exists_{h}$ with $\widehat{s} \models \delta$ and $\widehat{t} E \widehat{s}$. Such a point $\widehat{s}$ exists in $W$ since
$\widehat{t} \models \exists \delta$. We then add the ordered pairs $(s,s)$ to $R^\exists_{h}$, the ordered pairs $(t,s)$ to $E^\exists_{h}$,
and generate the least equivalence relation.

\vspace{1em}

\emph{$\forall H$-step $($Horizontal$)$}:
Let $W^{\forall H}_{h}  = W^\exists_{h}$, $R^{\forall H}_{h} = R^\exists_{h}$, and $E^{\forall H}_{h} = E^\exists_{h} $.
For each  $E^{\forall H}_{h}(t)\subseteq W^{\forall H}_{h}\backslash W_{h-1}^{\A H}$,
if $\forall\beta\in\Sigma^{\forall H}(t)$ but there is no $s\in W_{h}^{\forall H}$ already such that $tE_{h}^{\forall H}s$ and
$\widehat{s}\not\vDash\beta$, then we add a point $s$ to $W_{h}^{\forall H}$ with $\widehat{s} \not \models \beta$ and
$\widehat{t} E \widehat{s}$. Such a point $\widehat{s}$ exists in $W$ since $\widehat{t}$ is maximal with respect to $\forall \beta$.
We then add the ordered pairs $(s,s)$ to $R^{\forall H}_{h}$, the ordered pairs $(t,s)$ to $E^{\forall H}_{h}$, and generate the
least equivalence relation.

\vspace{1em}

\emph{$\forall$V-step $($Vertical$)$}:
Let $W^{\forall V}_{h}  = W^{\forall H}_{h} $, $R^{\forall V}_{h} = R^{\forall H}_{h}$, and $E^{\forall V}_{h} = E^{\forall H}_{h}$.
For each $E^{\forall V}_{h}(t)\subseteq W^{\forall V}_{h}\backslash W_{h-1}^{\A H}$, consider $\forall\gamma\in\Sigma^{\forall V}(t)$.
Since $\widehat{t}\not\vDash\forall\gamma$, we can pick a point $\widehat{s} \in W$ as in Lemma~\ref{lem:find maximal points new}(2).
We add the point $s$ to $W_{h}^{\A V}$ and $(t,s)$ to $R_{h}^{\A V}$.

Since $W$ satisfies commutativity, for each $w\in E_{h}^{\forall V}(t)$, there is $z_{w}\in W$ such that $\widehat{w}Rz_{w}$ and $z_{w}E \widehat{s}$.
To ensure commutativity is satisfied in our new frame, we add the points $s_{w}$ to $W_{h}^{\A V}$ where $\widehat{s_{w}}=z_{w}$.
We then add $(w,s_{w})$ to $R_{h}^{\A V}$ and take the reflexive and transitive closure.
We also add $(s_{w},s)$ to $E_{h}^{\A V}$ and generate the least equivalence relation.

\vspace{1em}

\emph{$\rightarrow$-step $($Vertical$)$}:
Let $W^{\rightarrow}_{h} = W^{\forall V}_{h} $, $R^{\rightarrow}_{h} = R^{\forall V}_{h}$, and $E^{\rightarrow}_{h} = E^{\forall V}_{h}$.
For each $t\in W_{h}^{\forall H} \backslash W_{h-1}^{\A H}$ (hence including any points added in the horizontal steps above,
but not in the previous vertical step), consider all $\alpha\rightarrow\sigma\in\Sigma^{\rightarrow}(t)$ such that
there is no $s\in W^{\rightarrow}_{h}$ already such that $tR^{\rightarrow}_{h}s$ and $\widehat{s}\not\vDash\alpha\rightarrow\sigma$
maximally. Consider
\[
A=[W \setminus \nu(\alpha\rightarrow\sigma)]\cap\bigcap_{\psi\in\Sub}\{\nu(\psi):\widehat{t}\vDash\psi\}.
\]
Then $A$ is clopen and $\widehat{t} \in A$, so by Lemma \ref{lem:find maximal points new}(3) there is $z \in A$ with $z \in \max E(A)$, $\widehat{t} Q z$, and $E(z)$ clean. We add the point $s$ to $W_{h}^{ \rightarrow}$ where $\widehat{s}=z$ ($s$ is a distinct new copy of $z$) and $(t,s)$ to $R_{h}^{ \rightarrow}$.

\begin{remark}
It is at this step that we have altered the construction given in \cite[\S 10.3]{Many_Dim_Modal_Logics}, in which witnesses for implications
are added in the same manner as in the $\A V$-step. In our version, we took an original $Q$-relation and turned it into an $R$-relation. The reason for this is that we cannot guarantee the existence of an $R$-successor of $t$ that is maximal with respect to $\alpha \rightarrow \sigma$ and at the same time belongs to a clean cluster.
\end{remark}

Before wrapping up the step, we show two properties of the chosen points.

\begin{lem}\label{MIPC_z is max}
The point $\widehat{s}=z$, as chosen above, is maximal  with respect to $\alpha\rightarrow\sigma$.
\end{lem}

\begin{proof}
Suppose $zRu$ for some
$u\not\vDash\alpha\rightarrow\sigma$. Since $zRu$ and each $\nu(\psi)$ in $\{\nu(\psi):\widehat{t}\vDash\psi\}$ is an upset, we have $u\in A$. Because $z \in \max A$, we obtain $z =u$. Thus, $z$ is maximal with respect to $\alpha\rightarrow\sigma$.
\end{proof}

\begin{lem}\label{implication t neq x}
$E(\widehat{t})\neq E(\widehat{s})$.
\end{lem}

\begin{proof}
If $E(\widehat{t}) = E(\widehat{s})$, then $\widehat{t}\in\max E(A)$ by Lemma \ref{lem:cleanclustersmipc}. Since $\widehat{t} \in A$, we have $\widehat{t} \in \max A$. Therefore, the same argument as in the proof of the previous lemma yields that $\widehat{t}$ is maximal with respect to $\alpha \rightarrow \sigma$. This contradicts $\alpha\rightarrow\sigma\in\Sigma^{\rightarrow}(t)$.
\end{proof}

We wrap up the $\rightarrow$-step the same way as the $\forall V$-step. Since $W$ satisfies commutativity, for each $w\in E_{h}^{\forall V}(t)$ there is $z_w\in W$ with $\widehat{w}Rz_{w}$ and $z_w E \widehat{s}$.
We add the points $s_{w}$ to $W_{h}^{\A V}$ where $\widehat{s_{w}}=z_{w}$. We then add $(w,s_{w})$ to $R_{h}^{\rightarrow}$ and take the
reflexive and transitive closure.
We also add $(s_{w},s)$ to $E_{h}^{\rightarrow}$ and generate the least equivalence relation.

\vspace{2mm}

To end this stage of the construction, we let $\mathfrak{F}_{h}=(W_{h},R_{h},E_{h})$ where
\[
W_{h}=W_{h}^{\rightarrow}, \ \ R_{h}=R_{h}^{\rightarrow}, \ \ E_{h}=E_{h}^{\rightarrow}.
\]

\begin{lem}
\label{lem:shapeofFh}
$\mathfrak F_{h}$ is a finite $\mathsf{M^+IPC}$-frame.
\end{lem}

\begin{proof}
First we show that $R_{h}$ is a partial order. Since in the $\exists$- and $\forall H$-steps we only added reflexive arrows to $R_{h-1}$, the relation $R^{\forall H}_{h}$ is a partial order. By moving from $R^{\forall H}_{h}$ to $R_{h}$ we finished by taking the reflexive and transitive closure, hence $R_{h}$ is clearly reflexive and transitive. Antisymmetry of $R_{h}$ follows from the fact that every $R$-arrow added in the $\forall V$-step and $\rightarrow$-step is either reflexive or an arrow from a previously existing point into a freshly added point.

That $E_{h}$ is an equivalence relation is clear from the construction.
Moreover, the extra points added in the $\forall V$-step and $\rightarrow$-step make sure that commutativity is satisfied. In fact, the added points assure commutativity for immediate successors and by transitivity this implies commutativity for the whole frame.
Therefore, $\mathfrak F_{h}$ is an $\mathsf{MIPC}$-frame.

It follows from the construction that $\mathfrak F_{h}$ is finite. Thus, by Lemma~\ref{finite M+IPC}, it is left to show that $\mathfrak F_{h}$ has clean clusters.
Note that in the $\exists$-step and $\forall H$-step all freshly introduced $E_h$-relations are of the shape $(s, t)$ where either $s \in W_h$ and $t \in W^{\forall H}_{h} \setminus W_{h-1}$ or $s, t \in W^{\forall H}_{h} $. Since no non-reflexive $R_{h}$-arrows are introduced in these steps, no dirty cluster could have been built. We have already discussed the shape of the $R_{h}$ arrows introduced in the $\forall V$-step and $\rightarrow$-step.
This guarantees that no cluster in $W^{\forall H}_{h}$ is made dirty. The freshly introduced $E_{h}$-relations in these steps are of the shape $(s, t)$ where $s, t  \in W_{h} \setminus W^{\forall H}_{h}$. Since no non-reflexive $R_{h}$ relations exist between these points, we infer that all clusters are clean.
\end{proof}

\subsection{Auxiliary lemmas}

To prove that our construction terminates after finitely many steps, we require several auxiliary lemmas.

\begin{lem}\label{MIPC_relations}
Let $x,y\in W_h$.
\begin{enumerate}
\item If $xR_{h}y$ and $x\neq y$, then $\widehat{x}Q\widehat{y}$ and $E(\widehat{x})\neq E(\widehat{y})$.
\item If $xE_{h}y$, then $\widehat{x}E\widehat{y}$.
\item If $xQ_{h}y$, then $\widehat{x}Q\widehat{y}$.
\end{enumerate}
\end{lem}

\begin{proof}
(1) Observe that in the construction each non-trivial $R_{h}$-relation between immediate successors comes from either a non-trivial $R$-relation (as in the case of points added for commutativity or in the $\A V$-step) or a non-trivial $Q$-relation (as in the case of points added in the $\rightarrow$-step), in which case there is $w\in W$ with $\widehat{x}\neq w$, $\widehat{x}Rw$, and $wE\widehat{y}$. In that case we obviously have $\widehat{x} Q \widehat{y}$ and by Lemma~\ref{implication t neq x}, $E(\widehat{x})\neq E(\widehat{y})$ in $W$. Otherwise the relation $x R_h y$ was added by transitivity, so there is a chain $x = x_0 R_h x_1 R_h \dots R_h x_n = y$ of immediate $R_h$-successors to which the previous applies. In particular, by transitivity of $Q$ we have $\widehat{x_i} Q \widehat{y}$ for all $i \leq n$, so $\widehat{x} Q \widehat{y}$. Moreover, there is $z_1 \in W$ with $z_1 \neq x$, $\widehat{x} R z_1$ and $z_1 E \widehat{x_1}$. Since $\widehat{x_1} Q \widehat{y}$, by commutativity there is $y'$ in $W$ with $z_1 R y'$ and $y' E \widehat{y}$. If $\widehat{x} E \widehat{y}$, then also $\widehat{x} E y'$ and $\widehat{x} R y'$ by transitivity of $R$. Since $\widehat{x}$ is from a clean cluster, this implies $\widehat{x} = y'$. Therefore, $\widehat{x} R z_1 R\widehat{x}$, and so $\widehat{x} = z_1$ by antisymmetry of $R$. This is a contradiction since $\widehat{x} \neq z_1$. Thus, $E(\widehat{x})\neq E(\widehat{y})$.

(2) It is obvious that each $E_{h}$-relation in $W_{h}$ comes from a pre-existing $E$-relation in $W$.

(3) If  $xQ_{h}y$, then there is $z$ with $xR_{h}z$ and $zE_{h}y$. If $x = z$, then $x E_h y$, so $\widehat{x} E \widehat{y}$ by (2), and hence $\widehat{x} Q \widehat{y}$. If $x \neq z$, then $\widehat{x} Q \widehat{z}$ by (1). Also, $z E_{h} y$ implies $\widehat{z} E \widehat{y}$ by (2). Thus, $\widehat{x} Q \widehat{y}$.
\end{proof}

\begin{lem}[Persistence]\label{persistence}
If $u R_h w$, then $\widehat{u}\vDash\psi$ implies $\widehat{w}\vDash\psi$ for all $\psi\in\Sub$.
\end{lem}

\begin{proof}
Suppose $uR_{h}w$, $\psi\in\Sub$, and $\widehat{u}\vDash\psi$. It suffices to show the result for an immediate $R_{h}$-successor $w$ of $u$,
the general result then follows by induction. We consider how the $R_{h}$-arrow from $u$ to $w$ was added. By construction, either
$\widehat{u}R\widehat{w}$ or $w$ was added to witness some implication in $\Sigma^{\rightarrow}(u)$. If $\widehat{u}R\widehat{w}$, then
clearly $\widehat{u}\vDash\psi$ implies $\widehat{w}\vDash\psi$. If $w$ was added in a $\rightarrow$-step, then $w$ is specifically
chosen so that $\widehat{w}\in \nu(\gamma)$ for all $\gamma\in\Sub$ such that $\widehat{u}\vDash\gamma$. Thus, $\widehat{u}\vDash\psi$ implies
$\widehat{w}\vDash\psi$.
\end{proof}

\begin{lem}\label{MIPC_consistency of sigmas}
\begin{enumerate}
\item[]
\item If $t E_h u$, then $\Sigma^{\exists}(t) =\Sigma^{\exists}(u)$, $\Sigma^{\forall H}(t) =\Sigma^{\forall H}(u)$, and
$\Sigma^{\forall V }(t) = \Sigma^{\forall V}(u)$.
\item If $t R_h v$ and $\exists \gamma \in \Sigma^{\exists}(t) \cap \Sigma^{\exists}(v)$, then there are $u, w$ such that $t E_h u$,
$u R_h w$, $w E_h v$, $\widehat{u} \models \gamma$, and $\widehat{w} \models \gamma$.
\item If $t R_h v$ and $t\not=v$, then $\Sigma^{\forall H}(t) \cap\Sigma^{\forall H}(v) = \varnothing$.
\item Along an $R_{h}$-chain, each formula in $\{\forall \psi : \forall\psi\in\Sub\} \cup \{\E\psi : \E\psi\in\Sub\}$ can serve at most
once as a reason to enlarge a cluster in a horizontal step.
\item If $t R_h u$, then  $\Sigma^{\forall V }(u) \subseteq \Sigma^{\forall V}(t)$ and if $u$ was added as an immediate $R_h$-successor
to $t$ because of $\forall \alpha  \in \Sigma^{\forall V} (t)$, then $\Sigma^{\forall V }(u) \subset \Sigma^{\forall V}(t)$.
\item If $t R_h u$, then $\Sigma^{\rightarrow}(u) \subseteq \Sigma^{\rightarrow}(t)$ and if $u$ was added as an immediate $R_h$-successor
to $t$ because of $\alpha \rightarrow \beta \in \Sigma^{\rightarrow} (t)$, then $\Sigma^{\rightarrow}(u) \subset \Sigma^{\rightarrow}(t)$.
\end{enumerate}
\end{lem}

\begin{proof}
(1) Suppose $t E_h u$. Then $\widehat{t} E \widehat{u}$ by Lemma~\ref{MIPC_relations}(2). Therefore, $E(\widehat{t}) = E(\widehat{u})$ and $Q(\widehat{t}) = Q(\widehat{u})$. Thus, $\widehat{t} \vDash \exists \gamma$ iff $\widehat{u} \vDash\exists \gamma$, and $\widehat{t} \vDash \A \gamma$ iff $\widehat{u} \vDash \A \gamma$. Moreover, since $E(\widehat{t})$ is a clean cluster, $\widehat{t}$ is not maximal wrt $\A\gamma$ iff $\widehat{u}$ is not maximal wrt $\A\gamma$.
Consequently, $\Sigma^{\exists}(t) =\Sigma^{\exists}(u)$,
$\Sigma^{\forall H}(t) =\Sigma^{\forall H}(u)$, and $\Sigma^{\forall V }(t) = \Sigma^{\forall V}(u)$.

(2) Suppose $t R_h v$ and $\exists \gamma \in \Sigma^{\exists}(t) \cap \Sigma^{\exists}(v)$.
By the construction, there is $u$ with $t E_h u$ and $\widehat{u}\vDash\gamma$.
Since $\mathfrak F_h$ satisfies commutativity, there is $w$ with $u R_h w$ and $w E_h v$. By Lemma~\ref{persistence}, $\widehat{w}\vDash\gamma$.

(3) Suppose $t R_h v$ and $t\not=v$. Then $\widehat{t}Q\widehat{v}$ and $E(\widehat{t}) \neq E(\widehat{v})$ by Lemma~\ref{MIPC_relations}(1), so $\widehat{t}\neq\widehat{v}$. Thus, if $\forall \psi \in \Sigma^{\forall H}(t)$, then $\widehat{v}\vDash\A\psi$ by maximality of $\widehat{t}$, so $\A\psi\not\in\Sigma^{\A H}(v)$. Conversely, if $\forall \psi \in \Sigma^{\forall H}(v)$, then $\widehat{t}$ cannot be maximal with respect to $\A \psi$, so $\A\psi\not\in\Sigma^{\A H}(t)$.

(4) Let $\{v_i \mid i \in \mathbb N\}$ be an $R_h$-chain in $W_h$, i.e.~$v_i R_h v_{i+1}$ for all $i \in \mathbb N$. Suppose $\E\psi\in\Sub$.
Let $k$ be the least stage at which the formula $\E \psi$ has been used to
enlarge the cluster $E_{h}(v_{k})$ in a horizontal step. By (2), all $E_{h}(v_{l})$ for $l>k$ already contain a witness for $\psi$, so no cluster above will need to be enlarged in a horizontal step to witness the formula $\E\psi$. Now suppose $\A\psi\in\Sub$.
Let $l$ be a stage at which the formula $\A\psi$ has been used
to enlarge the cluster $E_{h}(v_{l})$ in a horizontal step. Then $\A\psi\in \Sigma^{\A H}(v_{l})$. By (3), $\A \psi\not\in\Sigma^{\A H}(v_{k})$ for $k\neq l$. Thus, $\A\psi$ is responsible for enlarging a cluster at most once in a horizontal step.

(5) We show the statement for immediate $R_{h}$-successors only, the general case follows by induction. Suppose $t R_h u$ and $\forall \psi \in \Sigma^{\forall V }(u)$. If $t=u$, then the result is clear. Suppose $t \neq u$. Since $tR_{h}u$, either $\widehat{t} R \widehat{u}$ in $W$ or $u$ was added as a successor of $t$ in some $\rightarrow$-step. If $\widehat{t} R \widehat{u}$, then $\forall \psi \in \Sigma^{\forall V }(t)$ by persistence (see Lemma~\ref{persistence}).
Suppose $u$ was added as an $R_h$-successor to $t$ as a witness to some implication. By the choice of $u$, we have $\widehat{u} \vDash\chi$ for all $\chi\in\Sub$ with $\widehat{t} \vDash \chi$. Therefore, if $\widehat{t}\vDash\A\psi$, then we would have $\widehat{u}\vDash\A\psi$, contradicting $\A\psi\in\Sigma^{\A V}(u)$. Thus, we must have $\widehat{t} \not \models\A\psi$. Moreover, since $\widehat{t}R\widehat{u}$, $\widehat{t}\neq\widehat{u}$, and $\widehat{u}\not\vDash\A\psi$, we have that
$\widehat{t}$ is not maximal with respect to $\A\psi$, so $\A\psi\in\Sigma^{\forall V} (t)$. Consequently, in either case we have
$\Sigma^{\forall V }(u) \subseteq \Sigma^{\forall V}(t)$.

Suppose that $u$ was added as an immediate $R_{h}$-successor to $t$ because of $\forall\alpha\in\Sigma^{\forall V}(t)$. Since
$\forall\alpha\in\Sigma^{\forall V}(t)$, we have $\widehat{t}\not\vDash\forall\alpha$ but $\widehat{t}$ is not maximal with respect to $\forall\alpha$. Since $u$ was added as an immediate $R_{h}$-successor of $t$ because of $\forall\alpha$, we specifically chose
$u$ so that $\widehat{u}\not\vDash\forall\alpha$ maximally, hence $\forall\alpha\not\in\Sigma^{\forall V}(u)$.

(6) We show the statement for immediate $R_h$-successors only, the general case follows by induction. Suppose $tR_{h}u$ and  $\alpha \rightarrow \beta \in \Sigma^{\rightarrow}(u)$. Then  $\widehat{u} \not\vDash\alpha \rightarrow \beta$ and $\widehat{u} \not \vDash\alpha$. If $t=u$, then the result is clear. Suppose $t \neq u$. Since $t R_h u$, either $\widehat{t} R \widehat{u}$ in $W$ or $u$ was added as a successor of $t$ in some $\rightarrow$-step. If $\widehat{t} R \widehat{u}$, then $\alpha \rightarrow \beta \in \Sigma^{\rightarrow}(t)$ by persistence (see Lemma~\ref{persistence}). Suppose $u$ was added as an $R_h$-successor to $t$ as a witness to some implication. By the choice of $u$, we have
$\widehat{u} \vDash\psi$ for all $\psi\in\Sub$ with $\widehat{t} \vDash \psi$. Therefore, we must have $\widehat{t} \not \models \alpha \rightarrow \beta$ and $\widehat{t} \not \models \alpha$, so $\Sigma^{\rightarrow}(u) \subseteq \Sigma^{\rightarrow}(t)$. Moreover, by construction, $\widehat{u}$ refutes
$\alpha \rightarrow \beta$ maximally (Lemma~\ref{MIPC_z is max}), and hence $\widehat{u} \vDash\alpha$. Thus, $\alpha \rightarrow \beta \not \in \Sigma^{\rightarrow}(u)$.
\end{proof}

\subsection{Termination of the construction}

With the aid of the auxiliary lemmas of the previous section, we will now prove that the end result of our construction is a finite frame. We will do this by looking at three important parameters of our frame: cluster size, $R$-branching, and $R$-depth.

\begin{defn} \label{defn: cluster size_branching_depth}
\begin{enumerate}
\item[]
\item A frame $\mathfrak{F}$ has \emph{bounded cluster size} if there exists $k\in\mathbb{N}$ such that $|E(t)|\leq k$ for all $t\in W$.
\item A frame $\mathfrak{F}$ has \emph{bounded R-branching} if there exists $m\in\mathbb{N}$ such that $t$ has at most $m$ distinct
immediate $R$-successors for all $t\in W$.
\item A frame $\mathfrak{F}$ has \emph{bounded R-depth} if there exists $n\in\mathbb{N}$ such that there is no $R$-chain in $\mathfrak{F}$
with more than $n$ distinct elements.
\end{enumerate}
\end{defn}

\begin{lem}\label{MIPC_finite}
Let $\mathfrak{F}=(W,R,E)$ be a partially ordered rooted augmented Kripke frame. If $\mathfrak{F}$ has bounded cluster size, bounded $R$-branching, and bounded $R$-depth, then $\mathfrak F$ is finite.
\end{lem}

\begin{proof}
Suppose $\mathfrak{F}=(W,R,E)$ is a partially ordered rooted augmented Kripke frame with bounded cluster size, $R$-branching, and $R$-depth. Consider the quotient $(W/{E},R_{E})$ whose worlds are the clusters $E(x)$ where $x\in W$ and $E(x)R_{E}E(y)$ iff $xQ y$. To see that $R_{E}$ is well defined, suppose $x Q y$, $x' \in E(x)$, and $y' \in E(y)$. Then $x'ExQyEy'$, so $x'Qy'$, and hence $R_E$ is well defined.

Because $Q$ is reflexive and transitive, so is $R_E$.
Since $R$ is a partial order and $\mathfrak F$ has bounded $R$-depth, from $x Q y $ and $y Q x$ it follows that $x E y$ by
\cite[Lem.~3(b)]{Bezhanishvili_1999}.
This shows that $R_E$ is anti-symmetric, and hence a partial order. Clearly $(W/{E},R_{E})$ is rooted since so is $\mathfrak F$.
Using commutativity in $\mathfrak F$ it is easy to verify that $(W/{E},R_{E})$ inherits bounded depth and bounded branching from $\mathfrak F$.
Since every rooted partial order with
the latter properties is finite, we have that $W/E$ is finite. Because $W$ has bounded cluster size, we conclude that $W$ is finite too.
\end{proof}

Let $m_{1},m_{2},m_{3}$ be the non-negative integers
\begin{flalign*}
&& m_{1} &= |\{\exists\psi:\exists\psi\in \Sub\}|&\\
&& m_{2} &= |\{\forall\psi:\forall\psi\in \Sub\}|&\\
&& m_{3} &= |\{\psi\rightarrow\chi:\psi\rightarrow\chi\in\Sub\}|.&
\end{flalign*}

\begin{lem}\label{MIPC_clusters}
$\mathfrak{F}_{h}=(W_{h},R_{h},E_{h})$ has cluster size bounded by $1+m_1+m_2$ for all $h<\omega$.
\end{lem}

\begin{proof}
Recall how the clusters of our frame are built. The `bottom cluster' of the starting point $t_{0}$ contains points added via the horizontal $\E$- and $\A H$-steps. After this, no more points are added to this cluster.

All other clusters are constructed as follows. First points of a new cluster are added via the vertical $\A V$- or $\rightarrow$-steps, and then the cluster is enlarged by the points added for commutativity. We refer to this stage as the `building phase' of the cluster. In the next round of the construction, the cluster is (possibly) enlarged via the two horizontal steps. After this, no more points are added to the cluster. In the horizontal steps, we enlarge the cluster for only two different reasons:
\[
\exists\gamma\in\Sigma^{\exists}(t) \mbox{ or } \forall\gamma\in\Sigma^{\forall H}(t).
\]
Thus, each enlargement of a cluster after its building phase is due to a formula in
\[
\{\forall \psi : \forall\psi\in\Sub\} \cup \{\E\psi : \E\psi\in\Sub\}.
\]
At the end of its building phase, the bottom cluster contains just one point. Observe that every cluster can be reached from the
bottom cluster by an $R_h$-chain. By Lemma~\ref{MIPC_consistency of sigmas}(4), every formula in
$\{\forall \psi : \forall\psi\in\Sub\} \cup \{\E\psi : \E\psi\in\Sub\}$ can serve at most once as a reason to enlarge a cluster
after its building phase along an $R_h$-chain. This entails that every cluster has size at most $1 +m_1 + m_2$.
\end{proof}

\begin{lem}\label{MIPC_branching}
$\mathfrak{F}_{h}=(W_{h},R_{h},E_{h})$ has $R_{h}$-branching bounded by $(1+m_{1}+m_{2})\cdot m_{3}+m_{2}$ for all $h<\omega$.
\end{lem}

\begin{proof}
Immediate $R_h$-successors are added in the $\forall V$-step and $\rightarrow$-step. First observe that since we are adding points to witness commutativity, every point in a cluster has the same number of immediate
$R_h$-successors by the end of a stage. Thus, it is enough to count the immediate successors of a point $t$ that we picked in the
$\forall V$-step.

To such a point $t$ we add immediate $R_h$-successors for three different reasons:
\begin{enumerate}
        \item $\A\gamma\in\Sigma^{\A V}(t)$,
        \item $\alpha \rightarrow \sigma \in\Sigma^{\rightarrow }(t)$, or
        \item $\alpha \rightarrow \sigma \in\Sigma^{\rightarrow }(y)$ for some $y \in E_h(t)$ with $y \neq t$.
    \end{enumerate}
The last reason covers the case where we add an $R_h$-successor to $t$ to witness commutativity. Note that all reasons occur at most once for
each formula in the respective sets. Therefore, reason (1) occurs at most $m_2$-times and reason (2) at most $m_3$-times. Finally, reason (3) occurs
at most $(m_{1}+m_{2})\cdot m_{3}$ times since by Lemma~\ref{MIPC_clusters} there are at most $m_1+ m_2$ points apart from $t$ in the cluster
of $t$. Thus, the $R_{h}$-branching of $\mathfrak F$ is bounded by
\[
m_2 + m_3 + (m_{1}+m_{2})\cdot m_{3} = m_{2} + (1+m_{1}+m_{2})\cdot m_{3}.
\]
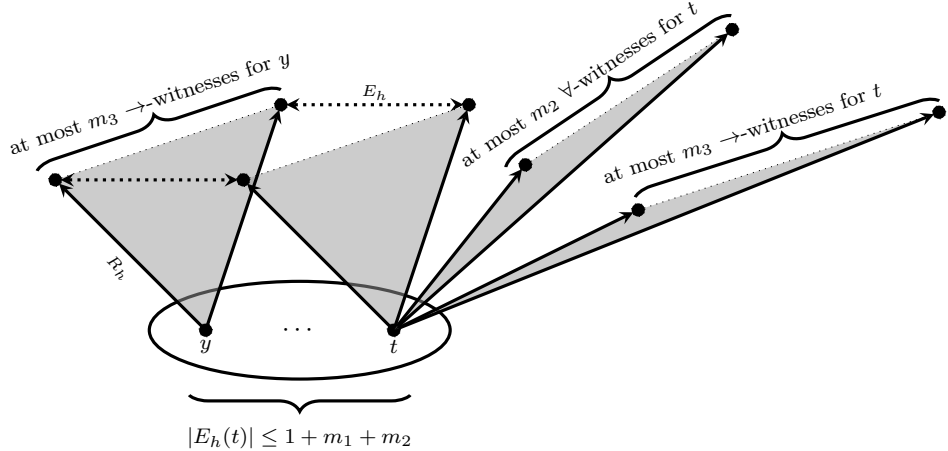
\begin{center}
\begin{minipage}{0.8\textwidth}
            \vspace{1pt}
            \begin{tikzpicture}
            [shorten <=2pt,shorten >=2pt,>=stealth]
            \draw [very thick] (5.25,0) ellipse (2cm and 0.65cm);
            \path[fill=gray!80, fill opacity=0.5] (4,0) -- (2,2) -- (5,3) -- cycle;
            \path[fill=gray!80, fill opacity=0.5] (6.5,0) -- (4.5,2) -- (7.5,3) -- cycle;
            \filldraw [dotted]
            (4,0) circle (2.5pt) node[align=center,   below] {\footnotesize{$y$}} --
            (2,2) circle (2.5pt) node[align=center, below] {}     --
            (5,3) circle (2.5pt) node[align=right,  below] {};
            \filldraw [dotted]
            (6.5,0) circle (2.5pt) node[align=center,   below] {\footnotesize{$t$}} --
            (4.5,2) circle (2.5pt) node[align=center, below] {}     --
            (7.5,3) circle (2.5pt) node[align=right,  below] {};
            \draw[decorate,decoration={brace,mirror,amplitude=8pt}, very thick] (3.75,-0.8) -- (6.75,-0.8) node [midway, below, yshift=-10pt]
            {\footnotesize{$\lvert E_h(t)\rvert\leq 1+m_{1}+m_{2}$}};
            \draw[->, very thick] (6.5,0) -- (4.5,2);
            \draw[->, very thick] (6.5,0) -- (7.5,3);
            \draw[->, very thick] (4,0) -- (2,2) node [midway, sloped, below] {\tiny{$R_h$}};
            \draw[->, very thick] (4,0) -- (5,3);
            \draw[dotted, <->, very thick] (2,2) -- (4.5,2);
            \draw[dotted, <->, very thick] (5,3) -- (7.5,3) node [midway,above, yshift=-2pt] {\tiny{$E_h$}};
            \draw[decorate,decoration={brace,amplitude=8pt}, very thick] (1.75,2.05) -- (5.05,3.15) node [midway, above, sloped, yshift=5pt]
            {\footnotesize{at most $m_3$ $\rightarrow$-witnesses for $y$}};
            \path[fill=gray!80, fill opacity=0.5] (6.5,0) -- (8.25, 2.2) -- (11,4) -- cycle;
            \draw[->, very thick] (6.5,0) -- (8.25,2.2);
            \draw[->, very thick] (6.5,0) -- (11,4);
            \filldraw [dotted]
            (6.5,0) --
            (8.25,2.2) circle (2.5pt) node[align=center, below] {}     --
            (11,4) circle (2.5pt) node[align=right,  below] {};
            \draw[decorate,decoration={brace, amplitude=8pt}, very thick] (8, 2.1) -- (11.05,4.15) node
            [pos=0.3, above, sloped, yshift=5pt, xshift=10pt] {\footnotesize{at most $m_2$ $\A$-witnesses for $t$}};
            \node (dots) at (5.25,0) {$\cdots$};
            \path[fill=gray!80, fill opacity=0.5] (6.5,0) -- (9.75,1.6) -- (13.75,2.9) -- cycle;
            \draw[->, very thick] (6.5,0) -- (9.75,1.6);
            \draw[->, very thick] (6.5,0) -- (13.75,2.9);
            \filldraw [dotted]
            (6.5,0) --
            (9.75,1.6)  circle (2.5pt) node[align=center, below] {}     --
           (13.75,2.9) circle (2.5pt) node[align=right,  below] {};
            \draw[decorate,decoration={brace, amplitude=8pt}, very thick] (9.7, 1.7) -- (13.8,3) node
            [pos=0.3, above, sloped, yshift=5pt, xshift=10pt] {\footnotesize{at most $m_3$ $\rightarrow$-witnesses for $t$}};
            \end{tikzpicture}
\end{minipage}
\end{center}
\end{proof}

\begin{lem}\label{MIPC_depth}
$\mathfrak{F}_{h}=(W_{h},R_{h},E_{h})$ has $R_{h}$-depth bounded by $(1+m_{1}+m_{2})\cdot(m_{2}+m_{3})$ for all $h<\omega$.
\end{lem}

\begin{proof}
The reason for adding an immediate successor to $t \in W_h$ via an $R_h$-relation is due to either a formula in $\Sigma^{\forall V}(t)$ or a formula in $\Sigma^{\rightarrow }(y)$ for some $y \in E_h(t)$ (as discussed in the proof of Lemma~\ref{MIPC_branching}). Let $s$ be a (not necessarily immediate) $R_h$-successor of $t$. Then $s$ could have been added via
direct formula witnessing, i.e.~there is an immediate predecessor $t'$ of $s$ with $t R_h t' R_h s$ and $s$ was added due to a formula in $\Sigma^{\forall V }(t')$ or $\Sigma^{\rightarrow }(t')$, or else $s$ was added to satisfy commutativity.

As we saw in Lemma~\ref{MIPC_consistency of sigmas}, moving up along an $R_{h}$-chain, the cardinality of the sets $\Sigma^{\forall V }(t)$
and $\Sigma^{\rightarrow }(t)$ does not increase, and it in fact decreases whenever an $R_h$-successor is added by direct formula witnessing. In particular, each point can have at most $m_{2}+m_{3}$ $R_h$-successors that have been added via direct formula witnessing and since in each cluster there are at most  $1+m_{1}+m_{2}$ points (Lemma~\ref{MIPC_clusters}), we have that the total $R_h$-depth cannot exceed $(1+m_{1}+m_{2})\cdot(m_{2}+m_{3})$.
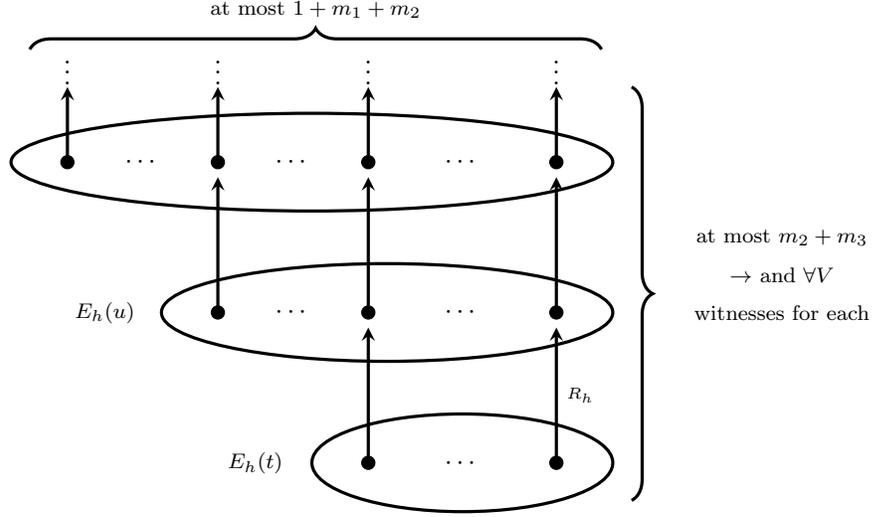
\begin{center}
\begin{minipage}{0.6\textwidth}
            \vspace{1pt}
            \begin{tikzpicture}
            [>=stealth]
            \node (Eht) at (2.5,0) {\footnotesize{$E_{h}(t)$}};
            \draw [very thick] (5.25,0) ellipse (2cm and 0.65cm);
            \filldraw (6.5,0) circle (2.5pt) node[align=center,   below] {};
            \node (f) at (5.25,0) {$\cdots$};
            \filldraw (4,0) circle (2.5pt) node[align=center,   below] {};
            \draw[->, very thick] (6.5,0) -- (6.5,1.8) node [midway,right] {\tiny{$R_{h}$}};
            \draw[->, very thick] (4,0) -- (4,1.8);

            \node (Ehu) at (0.5,2) {\footnotesize{$E_{h}(u)$}};
            \draw [very thick] (4.25,2) ellipse (3cm and 0.65cm);
            \filldraw (6.5,2) circle (2.5pt) node[align=center,   below] {};
            \node (d) at (5.25,2) {$\cdots$};
            \filldraw (4,2) circle (2.5pt) node[align=center,   below] {};
            \node (g) at (3,2) {$\cdots$};
            \filldraw (2,2) circle (2.5pt) node[align=center,   below] {};
            \draw[->, very thick] (6.5,2) -- (6.5,3.8);
            \draw[->, very thick] (4,2) -- (4,3.8);
            \draw[->, very thick] (2,2) -- (2,3.8);

            \draw [very thick] (3.25,4) ellipse (4cm and 0.65cm);
            \filldraw (6.5,4) circle (2.5pt) node[align=center,   below] {};
            \node (h) at (5.25,4) {$\cdots$};
            \filldraw (4,4) circle (2.5pt) node[align=center,   below] {};
            \node (i) at (3,4) {$\cdots$};
            \filldraw (2,4) circle (2.5pt) node[align=center,   below] {};
            \node (j) at (1,4) {$\cdots$};
            \filldraw (0,4) circle (2.5pt) node[align=center,   below] {};
            \draw[->, very thick] (6.5,4) -- (6.5,5);
            \draw[->, very thick] (4,4) -- (4,5);
            \draw[->, very thick] (2,4) -- (2,5);
            \draw[->, very thick] (0,4) -- (0,5);
            \node (k) at (6.5,5.3) {$\vdots$};
            \node (l) at (4,5.3) {$\vdots$};
            \node (m) at (2,5.3) {$\vdots$};
            \node (n) at (0,5.3) {$\vdots$};

            \draw[decorate,decoration={brace,amplitude=8pt}, very thick] (7.5,5) -- (7.5,-0.5);

            \node (p) at (9.5,3) {\footnotesize{at most $m_2+m_3$}};
            \node (q) at (9.5,2.5) {\footnotesize{$\rightarrow$ and $\A V$}};
            \node (r) at (9.5,2) {\footnotesize{witnesses for each}};

            \draw[decorate,decoration={brace,amplitude=8pt}, very thick] (-0.5,5.5) -- (7.1,5.5) node [midway, above, yshift=8pt]
            {\footnotesize{at most $1+m_1+m_2$}};

            \end{tikzpicture}
\end{minipage}
\end{center}
\end{proof}

\begin{lem}\label{lem:stop}
There is $h \in \mathbb N$ such that $\mathfrak F_{h'} = \mathfrak F_h$ for all $h' \geq h$.
\end{lem}

\begin{proof}
All points in the bottom cluster are added in round 1 and in each round we enlarge the $R_h$-length of a path by at most one. Thus, in stage $k$ of the construction, all $R_h$-chains are bounded by $k$.
The construction continues only until vertical witnesses are required. Since, by Lemma~\ref{MIPC_depth}, the $R_h$-depth of $\mathfrak F_k$ is bounded by $m = (1+m_{1}+m_{2})\cdot(m_{2}+m_{3})$, we have $\mathfrak F_{h'} = \mathfrak F_{m+1}$ for all $h' \geq m+1$.
\end{proof}

Set $\mathfrak{F}'=(W',R',E')$ where
\[
W'=W_{h},\quad R'={h}, \quad E'=E_{h},
\]
and $h$ is as in Lemma~\ref{lem:stop}. Then $\mathfrak F'$ is a finite $\mathsf{M^+IPC}$-frame by Lemma~\ref{lem:shapeofFh}.

\subsection{Truth lemma}

Define a valuation $\nu'$ on $W'$ by $\nu'(p)=\{t\in W':\widehat{t}\in \nu(p)\}$ for $p\in\Sub$ and $\nu'(q)=\varnothing$ for variables $q$ not occurring in $\varphi$. That $\nu'$ is well defined follows from Lemma~\ref{persistence}, which ensures that the
sets $\nu'(\psi)$ are in $\mathsf{Up}(\mathfrak{F'})$ for each $\psi\in\Sub$.

\begin{lem}[Truth Lemma]\label{MIPC_truth}
For all $t\in W'$ and $\psi\in\Sub$, we have $t\vDash'\psi$ iff $\widehat{t}\vDash\psi$.
\end{lem}

\begin{proof}
The proof is by induction on the complexity of $\psi$. The base cases $\psi=\bot$ and $\psi=p$ ($p$ a propositional variable) follow from
the definition, and the cases $\psi=\psi_{1}\wedge\psi_{2}$ and $\psi=\psi_{1}\vee\psi_{2}$ are easily verified. So we focus on the
cases $\psi=\psi_{1}\rightarrow\psi_{2}$ (and hence $\psi=\neg\psi_{1}=\psi_{1}\rightarrow\bot$), $\psi=\exists\psi_{1}$, and
$\psi=\forall\psi_{1}$.

\vspace{1em}

\noindent\underline{$\rightarrow$ case}:
Let $\psi=\psi_{1}\rightarrow\psi_{2}$ and $t\in W'$. Suppose $t\not\vDash'\psi_{1}\rightarrow\psi_{2}$. Then $tR's$
for some $s\in W'$ with $s\vDash'\psi_{1}$ and $s\not\vDash'\psi_{2}$. By the inductive hypothesis, $\widehat{s}\vDash\psi_{1}$ and $\widehat{s}\not\vDash\psi_{2}$. Thus, $\widehat{s} \not\vDash\psi_{1}\rightarrow\psi_{2}$. Since $tR's$, we have $\widehat{t} \not\vDash\psi_{1}\rightarrow\psi_{2}$ by persistence (Lemma \ref{persistence}).

Conversely, suppose $\widehat{t}\not\vDash\psi_{1}\rightarrow\psi_{2}$. If $\widehat{t}\vDash\psi_{1}$, then we have $\widehat{t}\vDash\psi_{1}$ but $\widehat{t}\not\vDash\psi_{2}$. By the inductive hypothesis, $t\vDash'\psi_{1}$ but $t\not\vDash'\psi_{2}$. By construction, $tR't$, so $t\not\vDash'\psi_{1}\rightarrow\psi_{2}$. If
$\widehat{t}\not\vDash\psi_{1}$, then in the $\rightarrow$-step of the stage immediately after $t$ is added to $W'$, we add $s$ to $W'$ and $tR's$ where $\widehat{s}\not\vDash\psi_{1}\rightarrow\psi_{2}$ maximally (Lemma \ref{MIPC_z is max}). Thus, $\widehat{s}\vDash\psi_{1}$ and $\widehat{s}\not\vDash\psi_{2}$, so by the inductive hypothesis, $s\vDash'\psi_{1}$ and $s\not\vDash'\psi_{2}$. Since $tR's$, we conclude that $t\not\vDash'\psi_{1}\rightarrow\psi_{2}$.

\vspace{1em}

\noindent\underline{$\E$ case}:
Let $\psi=\exists\psi_{1}$ and $t\in W'$. Suppose $t\vDash'\exists\psi_{1}$. Then $tE's$ for some $s\in W'$ with $s\vDash'\psi_{1}$. By the inductive hypothesis, $\widehat{s}\vDash\psi_{1}$, and $tE's$ implies $\widehat{t}E\widehat{s}$ by Lemma~\ref{MIPC_relations}(2). Thus, $\widehat{t}\vDash\exists\psi_{1}$.

Conversely, suppose $\widehat{t}\vDash\exists\psi_{1}$. Then
$\exists\psi_{1}\in\Sigma^{\exists}(t)$, so in the $\exists$-step of the next stage of the construction after $t$ is added, we add $s$ to
$W'$ and $(t,s)$ to $E'$ where $s$ is a copy of some $\widehat{s}\in W$ with $\widehat{t}E\widehat{s}$ and
$\widehat{s}\vDash\psi_{1}$. By the inductive hypothesis, $s\vDash'\psi_{1}$. Since $tE's$, we conclude that
$t\vDash'\exists\psi_{1}$.

\vspace{1em}

\noindent\underline{$\A$ case}:
Let $\psi=\forall\psi_{1}$ and $t\in W'$. Suppose $t\not\vDash'\forall\psi_{1}$. Then $tQ'w$ for some
$w\in W'$ with $w\not\vDash'\psi_{1}$. By the inductive hypothesis, $\widehat{w}\not\vDash\psi_{1}$, and $tQ'w$
implies $\widehat{t}Q\widehat{w}$ by Lemma~\ref{MIPC_relations}(3). Thus, $\widehat{t}\not\vDash\forall\psi_{1}$.

Conversely, suppose $\widehat{t}\not\vDash\forall\psi_{1}$.
If $\widehat{t}$ is maximal with respect to $\forall\psi_{1}$, then $\forall\psi_{1}\in\Sigma^{\forall H}(t)$, so at some point in the construction of the next stage after $t$ is added, we add $s$ to $W'$ and $(t,s)$ to $E'$ where $s$ is a copy of some $\widehat{s}\in W$ with
$\widehat{t}E\widehat{s}$ and $\widehat{s}\not\vDash\psi_{1}$. By the inductive hypothesis, $s\not\vDash'\psi_{1}$, so
$t\not\vDash'\psi_{1}$. If $\widehat{t}$ is not maximal, then we add $s$ to $W'$ and $(t,s)$ to $R'$ where $s$ is a copy of some $\widehat{s}\in W$ and $\widehat{s}$ is maximal with respect to $\forall\psi_{1}$. Therefore,
$\forall\psi_{1}\in\Sigma^{\forall H}(s)$, and in the next stage we add
$w$ to $W'$ and $(s,w)$ to $E'$ where
$\widehat{w}\in W$ and $\widehat{w}\not\vDash\psi_{1}$. But then $tQ'w$, hence $\widehat{t}Q\widehat{w}$ (see Lemma~\ref{MIPC_relations}),
and by the inductive hypothesis, $w\not\vDash'\psi_{1}$. Thus, $t\not\vDash'\forall\psi_{1}$.
\end{proof}

The FMP of $\mathsf{M^+IPC}$ is now an immediate consequence of the above.

\begin{theorem}
\label{MIPC_FMP}
$\mathsf{M^+IPC}$ has the finite model property.
\end{theorem}

\begin{proof}
Suppose $\mathsf{M^+IPC} \not \vdash \varphi$. By completeness of $\mathsf{M^+IPC} $ with respect to descriptive frames, there are a descriptive $\mathsf{M^+IPC} $-frame $\mathfrak F$ and a valuation $\nu$ on $\mathfrak F$ such that $(\mathfrak F, \nu) \not \models \varphi$. Let $\mathfrak F'$ be the finite $\mathsf{M^+IPC}$-frame constructed above. Since $t_0$ was chosen so that $\widehat{t_0}$ refutes $\varphi$ in $\mathfrak F$, by Lemma \ref{MIPC_truth}, $t_0$ refutes $\varphi$ in $\mathfrak F'$.
We thus found a finite $\mathsf{M^+IPC}$-frame refuting $\varphi$.
\end{proof}

Since $\mathsf{M^{+}IPC}$ is finitely axiomatizable and has the finite model property, as an immediate corollary to Theorem~\ref{MIPC_FMP},
we obtain decidability of $\mathsf{M^{+}IPC}$, meaning that there is an effective method for determining whether an arbitrary formula is a
theorem of $\mathsf{M^{+}IPC}$.

\begin{cor}
$\mathsf{M^{+}IPC}$ is decidable.
\end{cor}

\begin{remark}
\label{rem:Ono-Suzuki}
Another consequence of Theorem~\ref{MIPC_FMP} is that
$\mathsf{M^{+}IPC}$ is the monadic fragment of the intermediate predicate logic of Casari, which is obtained by adding to $\mathsf{IQC}$ the Casari formula $\mathsf{Cas}$. This can be seen by utilizing the Translation Theorem of Ono and Suzuki (see \cite[Thm.~3.5]{Ono_Suzuki_1988}).
\end{remark}

\section{The finite model property of $\mathsf{M^{+}Grz}$}
\label{section:fmpmodalcase}

In this section we prove that $\mathsf{M^{+}Grz}$ has the finite model property. Our proof, which consists of three steps, is a mixture of selective and standard filtration techniques. The main reasons why the same technique as for $\mathsf{M^{+}IPC}$ does not work is the lack of persistence in $\mathsf{M^{+}Grz}$-models and the fact that witnesses for $\A$-formulas cannot be chosen maximally wrt $Q$-relations. A rough structure of the proof is as follows.

Suppose $\mathsf{M^{+}Grz}\not\vdash\varphi$. Then there is a descriptive $\mathsf{M^{+}Grz}$-frame $\mathfrak{F}_{0}=(W_{0},R_{0},E_{0},P_{0})$ and a valuation $\nu_{0}$ on $W_{0}$ such that $\mathfrak{F}_{0}\not\vDash_{0}\varphi$. We build a finite $\mathsf{M^{+}Grz}$-frame from $\mathfrak{F}_{0}$ in three steps:

\begin{enumerate}
    \item First we select a (possibly infinite) partially ordered $\mathsf{MS4}$-frame $\mathfrak{F}_{1}=(W_1,R_1,E_1)$ from $\mathfrak{F}_{0}$, in which all clusters are clean and $\varphi$ is refuted. An important feature of this step is that $R_1$ is not simply the restriction of $R_0$ to $W_1$, but rather its strengthening. Its construction resembles the construction of $R$-relations from $Q$-relations in the $\rightarrow$-step of the $\mathsf{M^{+}IPC}$-construction.
    \item Next we construct a (possibly infinite) partially ordered $\mathsf{MS4}$-frame $\mathfrak{F}_{2}$ from $\mathfrak{F}_{1}$, in which all clusters are both clean and finite and $\varphi$ is refuted. In this step we use standard filtration to collapse $E_1$-clusters of $\mathfrak F_1$ so that each cluster contains only one point representing all points that satisfy the same formulas of $\mathsf{Sub}(\varphi)$.
    \item Finally, as in Step 1, we use selective filtration to construct a finite partially ordered $\mathsf{MS4}$-frame $\mathfrak{F}_{3}$ from $\mathfrak{F}_{2}$, in which all clusters are clean (hence $\mathfrak{F}_{3}$ is an $\mathsf{M^{+}Grz}$-frame) and $\varphi$ is refuted. This step resembles the $\mathsf{M^{+}IPC}$-construction, but in order for $\mathfrak F_3$ to inherit the bounded cluster size from $\mathfrak F_2$, we need to add only a single copy of an original point in $\mathfrak F_2$ to a cluster.
\end{enumerate}

\subsection{Step 1: Constructing $\mathfrak{F}_{1}$}

Let $\mathfrak{F}_{0}=(W_{0},R_{0},E_{0})$ be as above. For $x, y \in W_0$ let
\begin{equation*}
x \strict{Q_0} y \quad \text{iff} \quad \text{there is } w\in W_0 \text{ such that }  w\neq x, xR_{0}w, \text{ and } wE_{0}y.
\end{equation*}

We construct $\mathfrak{F}_{1}=(W_{1},R_{1},E_{1})$ as follows:

\begin{itemize}
    \item $W_{1}= \{x\in W_{0} \mid x\in\max_{R_0}E_{0}(A)\text{ for some clopen } A \text{ of } \mathfrak F_0\}$.
    \item $xR_{1}y\Leftrightarrow x=y$ or $x\strict{Q_{0}}y$ and $x\vDash_{0}\B\psi\Rightarrow y\vDash_{0}\B\psi$ for all $\B\psi\in\Sub$.
    \item $xE_{1}y\Leftrightarrow xE_{0}y$.
    \item We define a valuation $\nu_1$ on $\mathfrak F_1$ by  $\nu_{1}(p)=\{x : x\in \nu_{0}(p)\}$ for all $p\in\Sub$, and $\nu_{1}(q)=\varnothing$ for all other propositional variables $q$.
\end{itemize}

We first show that there is a point in $W_{1}$ which refutes $\varphi$ (in $\mathfrak F_{0}$).

\begin{lem}\label{start}
    There is $v\in W_{0}$ such that $v\not\vDash_{0}\varphi$ $(R_{0}$-maximally$)$ and $v\in\max_{R_0}E_{0}(\nu(\neg\varphi))$ $($hence $E_{0}(v)$ is clean and $v\in W_{1})$.
\end{lem}

\begin{proof}
    Since $\mathfrak{F}_{0}\not\vDash_{0}\varphi$, there is $t\in W_{0}$ such that $t\not\vDash_{0}\varphi$. Then $t\in \nu_{0}(\neg\varphi)$, so $t\in E_{0}(\nu_{0}(\neg\varphi))$. Because descriptive augmented Kripke frames satisfy $A\in P_{0}\Rightarrow E_{0}(A)\in P_{0}$, we have $E_{0}(\nu_{0}(\neg\varphi))\in P_{0}$. Thus, Lemma \ref{lem:S4canseeqmax} yields $u\in W_{0}$ with $tR_{0}u$ and $u\in \max_{R_{0}}E_{0}(\nu_{0}(\neg\varphi))$. Since $u\in E_{0}(\nu_{0}(\neg\varphi))$, there is $v\in W_{0}$ with $uE_{0}v$ and $v\not\vDash_{0}\varphi$. We now show that $v$ is our desired point. Because $u\in\max_{R_{0}}E_{0}(\nu_{0}(\neg\varphi))$, the cluster $E_{0}(u)=E_{0}(v)$ is clean (Lemma \ref{lem:mCasGrz}). We show that $v\in \max_{R_{0}}E_{0}(\nu_{0}(\neg\varphi))$. Suppose $vR_{0} w$ for some $w\in E_{0}(\nu_{0}(\neg\varphi))$. By commutativity, there is $u'$ such that $uR_{0}u'$ and $u'E_{0}w$. Then $u'\in E_{0}(\nu_{0}(\neg\varphi))$, so $u\in \max_{R_{0}}E_{0}(\nu_{0}(\neg\varphi))$ implies $u=u'$. Thus, $vR_{0}w$ and $vE_{0}w$, yielding that $v=w$ as $v$ is in a clean cluster. Now, since $v\in\max_{R_{0}}E_{0}(\nu_{0}(\neg\varphi))$ and $v\in \nu_{0}(\neg\varphi)$, it is easy to see that $v$ is $R_{0}$-maximal with respect to $\varphi$, hence is our desired point.
\end{proof}

We next highlight some fundamental properties of $\mathfrak{F_{1}}$.

\begin{lem}\label{F1 props}
    \begin{enumerate}
        \item[]
        \item $E_0(x) \subseteq W_0$ is a clean cluster in $\mathfrak F_0$ for all $x \in W_1$.
        \item\label{R1 Lem - E0 contained in W1} If $x\in W_{1}$, then $E_{0}(x)\subseteq W_{1}$.
        \item\label{Qbarcharacterization} $x\strict{Q_{0}}y$ iff  $xQ_{0}y$ but $x\nE_{0}y$ for all $x, y \in W_1$.
        \item The restriction of $\strict{Q_{0}}$ to $W_1$ is a strict partial order.
       \item $R_{1}$ is a partial order.
       \item $E_1$ is an equivalence relation.
       \item\label{R1 Lem - W1 commutative} $R_{1}$ and $E_{1}$ satisfy commutativity.
      \item\label{R1 Lem - E1(x) clean}
      $\mathfrak F_1$ has clean clusters.
       \item\label{R1 Lem - can find R0 max y in W1} For $x\in W_{1}$ and $\B\gamma\in\Sub$, if $x\not\vDash_{0}\B\gamma$, then there is $y\in W_{1}$ such that $xR_{1}y$, $y\in A \cap \max_{R_0}E_{0}(A)$, where
\[
        A=\nu_{0}(\neg\B\gamma)\cap\bigcap\{\nu_{0}(\B\psi) \mid \B\psi\in\Sub \text{ and } x\vDash_{0}\B\psi\},
\]
and $y\not\vDash_{0}\B\gamma$ $R_{0}$-maximally, hence $y \not\vDash_{0} \gamma$.
    \end{enumerate}
\end{lem}

\begin{proof}
(1) This is an immediate consequence of Lemma \ref{lem:mCasGrz}.

(2) Let $x \in W_1$ and $y \in E_0(x)$. Then $x \in \max_{R_0}E(A)$ for some clopen $A \subseteq W_0$. Therefore, $E_0(x)$ is clean by (1), and so $y \in \max_{R_0}E(A)$ by Lemma \ref{lem:maxgrzclean}(1). Thus, $y \in W_1$.

(3) The implication from right to left is obvious. For the converse, suppose that $x, y \in W_1$ and there is $w\in W_0$ such that $w\neq x, xR_{0}w, \text{ and } wE_{0}y$. Then clearly $xQ_{0}y$. Also, since $x$ is from a clean cluster, $x\nE_{0}w$. Thus, $x\nE_{0}y$.

(4) Irreflexivity of  $\strict{Q_{0}}$ on $W_1$ follows from the reflexivity of $E_0$ and (3). We show that $\strict{Q_{0}}$ is transitive on $W_1$. Suppose $x \strict{Q_{0}} y \strict{Q_{0}} z $ for $x, y, z \in W_1$. Then there are $y' \neq x$ and $z' \neq y$ with $x R_0 y'$, $y' E_0 y$ and $y R_0 z'$ and $z' E_0 z$. By commutativity, there is $z''$ with $y' R_0 z''$ and $z'' E_0 z$. Therefore, $x R_0 z''$ and $z'' E_0 z$. If we had $x = z''$, then we would obtain $x R_0 y' R_0 x$, and so $x = y'$ by Lemma \ref{lem:maxgrzclean}(2). The latter contradicts the choice of $y'$. Thus, $z'' \neq x$ and so $x \strict{Q_{0}} z $.

(5) $R_1$ is reflexive by definition. To see that $R_1$ is transitive, suppose $x, y, z \in W_1$ with $x R_1 y R_1 z$. Without loss of generality we may assume that $x,y,z$ are pairwise distinct. Then $x \strict{Q_{0}}y$ and  $y \strict{Q_{0}} z$, so $x \strict{Q_{0}}z$ by (3). Moreover, if $ x \models \Box \psi$ for $\Box \psi \in \mathsf{Sub}(\varphi)$, then since $x R_1 y R_1 z$, we have $ y \models \Box \psi$ and so $z \models \Box \psi$. Therefore, $R_1$ is transitive. Finally, if $x R_1 y R_1 x$ and $x \neq y$, then $x \strict{Q_{0}} y \strict{Q_{0}} x$. The latter implies $x \strict{Q_{0}} x$ by transitivity of $\strict{Q_{0}}$, which contradicts irreflexivity of $\strict{Q_{0}}$. Thus, $R_1$ is anti-symmetric.

(6) This is immediate since $E_1$ is an equivalence relation.

(7) Suppose that $xR_{1}y$ and $xE_{1}z$. Without loss of generality we may assume that $x\neq y$ and $x\neq z$. Then $x\strict{Q_{0}}y$, so there is $u\in W_{0}$ such that $x\neq u$, $xR_{0}u$, and $uE_{0}y$. By commutativity in $W_{0}$, there is $v$ such that $zR_{0}v$ and $vE_{0}u$. We show that $v$ is the required witness for commutativity in $W_{1}$. From $vE_{0}u$ and $uE_{0}y$ we have $vE_{0}y$, so $v\in W_{1}$ by (2). Because $x\neq u$, $xR_{0}u$, and $x$ is from a clean cluster, we have $x\nE_{0} u$ . Thus, $z\nE_{0} v$. In particular, $z\neq v$, and so $z\strict{Q_{0}}v$. Moreover, $zR_{0}v$ gives that if $z\vDash_{0}\B\gamma$, then $v\vDash_{0}\B\gamma$, so $zR_{1}v$. From $vE_{0}y$ we have $vE_{1}y$, yielding commutativity in $W_{1}$.

(8) Suppose there are $x,y\in W_{1}$ with $x\neq y$, $xE_{1}y$, and $xR_{1}y$. Since $xE_{1}y$, we have $xE_{0}y$, and because $xR_{1}y$ and $x\neq y$, we have $x\strict{Q_{0}}y$. Thus, there is $w\in W_{0}$ with $x\neq w$, $xR_{0}w$,  and $wE_{0}y$. From $xE_{0}y$ and $yE_{0}w$ we have $xE_{0}w$. By (1), $x$ is chosen from a clean cluster in $W_{0}$, so $xR_{0}w$ and $xE_{0}w$ imply $x=w$, a contradiction.

(9) Suppose $x\not\vDash_{0}\B\gamma$. Consider
\[
A=\nu_{0}(\neg\B\gamma)\cap\bigcap\{\nu_{0}(\B\psi) \mid \B\psi\in\Sub \mbox{ and } x\vDash_{0}\B\psi\}.
\]
Clearly $x\in A$, so $x\in E_{0}(A)$. We have $x\in\max_{R_{0}}E_{0}(A)$ or $x\not\in\max_{R_{0}}E_{0}(A)$.

\noindent\underline{Case 1: $x\in\max_{R_{0}}E_{0}(A)$}

If $x\in\max_{R_{0}}E_{0}(A)$, then from $xR_{0}w$ and $x\neq w$ it follows that $w\not\in E_{0}(A)$, so $w\not\in A$. But $xR_{0}w$ implies ${w\in \bigcap\{\nu_{0}(\B\psi) \mid \B\psi\in\Sub \mbox{ and } x\vDash_{0}\B\psi\}}$, so we must have $w\not\in \nu_{0}(\neg\B\gamma)$. Therefore, $w\vDash_{0}\B\gamma$. Since $x\not\vDash_{0}\B\gamma$ but $w\vDash_{0}\B\gamma$ for all $w\neq x$ with $xR_{0}w$, we must have $x\not\vDash_{0}\B\gamma$ $R_{0}$-maximally. Thus, $x\not\vDash_{0}\gamma$, and $x$ is our desired point.

\noindent\underline{Case 2: $x\not\in\max_{R_{0}}E_{0}(A)$}

If $x\not\in\max_{R_{0}}E_{0}(A)$, then Lemma~\ref{lem:grzchar}(2) yields $t\in \max_{R_{0}}E_{0}(A)$ such that $x\neq t$ and $xR_{0}t$. But then $tE_{0}y$ for some $y\in A$. Since $t\in \max_{R_{0}}E_{0}(A)$, we have $t\in W_{1}$, so $y\in W_{1}$ by (2).
From $x\neq t$ and $xR_{0}t$ it follows that $x\strict{Q_{0}}y$. By the choice of $y\in A$, if $x\vDash_{0}\B\psi$ then $y\vDash_{0}\B\psi$ for all $\Box \psi \in \Sub$, so $xR_{1}y$. Since $y\in A$, we have $y\not\vDash_{0}\B\gamma$. To see that $y\not\vDash_{0}\B\gamma$ $R_{0}$-maximally, suppose $yR_{0}z$ and $z\not\vDash_{0}\B\gamma$. If $x\vDash_{0}\B\psi$, then $y\vDash_{0}\B\psi$ (as $y\in A$), so $yR_{0}z$  implies $z\vDash_{0}\B\psi$. Thus, $z\in A$, hence $z\in E_{0}(A)$, and maximality of $y$ in $E_{0}(A)$ yields $y=z$. Consequently, $y$ is $R_{0}$-maximal with respect to $\B\gamma$.
\end{proof}

We conclude Step 1 by proving the truth lemma for $\mathfrak{F}_{1}$.

\begin{lem}[Truth Lemma]\label{truth F1}
For $x\in W_{1}$ and $\psi\in\Sub$,
\[
(\mathfrak{F}_{0},x)\vDash_{0} \psi \Leftrightarrow (\mathfrak{F}_{1},x)\vDash_{1}\psi.
\]
\end{lem}

\begin{proof}
The proof is by induction on the complexity of $\psi$. The base case $\psi=p$ is clear from the definition of $\nu_{1}$. The cases of $\psi=\psi_{1}\wedge\psi_{2}$ and $\psi=\neg\psi_{1}$ are straightforward, so we focus on the cases $\psi=\A\psi_{1}$ and $\psi=\B\psi_{1}$.

    Suppose $\psi=\A\psi_{1}$. If $x\not\vDash_{0}\A\psi_{1}$, then $xE_{0}y$ for some $y\not\vDash_{0}\psi_{1}$. By Lemma \ref{F1 props}(2), $y \in W_1$, so $y\not\vDash_{1}\psi_{1}$ by the inductive hypothesis. From $xE_{0}y$ we have $xE_{1}y$ by the definition of $E_1$. Thus, $x\not\vDash_{1}\A\psi_{1}$. The proof of the converse implication is immediate.

    Suppose $\psi=\B\psi_{1}$. If $x\not\vDash_{0}\B\psi_{1}$, then by Lemma~\ref{F1 props}(9), there is $y\in W_{1}$ such that $xR_{1}y$ and $y\not\vDash_{0}\psi_{1}$. By the inductive hypothesis, $y\not\vDash_{1}\psi_{1}$, hence $x\not\vDash_{1}\B\psi_{1}$.
     Conversely, if $x\not\vDash_{1}\B\psi_{1}$, then there is $y\in W_{1}$ such that $xR_{1}y$ and $y\not\vDash_{1}\psi_{1}$. By the inductive hypothesis, $y\not\vDash_{0}\psi_{1}$. If $x=y$, then $x\not\vDash_{0}\psi_{1}$, hence $x\not\vDash_{0}\B\psi_{1}$. If $x\neq y$, then as $xR_{1}y$, we have $x\strict{Q_{0}}y$ and $x\vDash_{0}\B\gamma$ implies $y\vDash_{0}\B\gamma$ for all $\B\gamma\in\Sub$. Since $y\not\vDash_{0}\psi_{1}$, we have $y\not\vDash_{0}\B\psi_{1}$. Thus, $x\not\vDash_{0}\B\psi_{1}$.
\end{proof}

\subsection{Step 2: Constructing $\mathfrak{F}_{2}$}

In this step we use the standard filtration technique to construct $\mathfrak F_2$ from $\mathfrak F_1$ by `collapsing' $E_{1}$-clusters into finitely many classes. Thus, $\mathfrak F_2$ will have finitely many clusters.

Define an equivalence relation $\sim$ on $W_{1}$ by
\[
x\sim y \Leftrightarrow \left(xE_{1}y \text{ and } x\vDash_{1}\gamma \Leftrightarrow y\vDash_{1}\gamma \text{ for all }\gamma\in\Sub \right).
\]

We construct $\mathfrak{F}_{2}=(W_{2},R_{2},E_{2})$ as follows:
\begin{itemize}
    \item $W_{2}= W_{1}/{\sim}=\{[x] : x\in W_{1}\}$ where $[x]$ denotes the $\sim$-equivalence class of $x$.
    \item For $[x],[y]\in W_{2}$, $[x]R_{2}[y]\Leftrightarrow [x]=[y]$ or $xR_{1}y$.
    \item For $[x],[y]\in W_{2}$, $[x]E_{2}[y]\Leftrightarrow xE_{1}y$.
    \item $\nu_{2}(p)=\{[x] : x\in \nu_{1}(p)\}$ for all $p\in\Sub$, and $\nu_{2}(q)=\varnothing$ for all other propositional variables $q$.
\end{itemize}

\begin{lem} The relations $E_2$ and $R_2$ are well defined, and so is the valuation $v_2$.
\end{lem}

\begin{proof}
It is easy to see that $E_2$ and $v_2$ are well defined. We show that $R_2$ is well defined. Let $x, y, x', y' \in W_1$ with $x\sim x'$, $y\sim y'$, and $[x]R_{2}[y]$.
Then $[x]=[y]$ or $xR_{1}y$. If $[x]=[y]$, we have $[x']=[x]=[y]=[y']$, and so $[x']R_{2}[y']$. 
If $xR_{1}y$, then $x=y$ or $x\strict{Q_{0}}y$ and $x\vDash_{0}\B\gamma$ implies $y\vDash_{0}\B\gamma$ for all $\B\gamma\in\Sub$. The former case implies $[x]=[y]$ which we have already considered. In the latter case,
from $x\strict{Q_{0}}y$ it follows that $xQ_{0}y$ and $x \nE_0 y$ by Lemma \ref{F1 props}(3).
Note that $x' \sim x$ implies $x' E_1 x$ and so $x ' E_0 x$. Similarly, $y ' E_0 y$. By transitivity of $Q_0$ we thus have $x' Q_0 y'$. Moreover, $x ' E_0 x$, $y ' E_0 y$, and $x \nE_0 y$ imply that $x' \nE_0 y'$. Thus, $x'\strict{Q_{0}}y'$ by Lemma~\ref{F1 props}(3).
If $\Box \gamma \in \Sub$ and $x'\vDash_{0}\B\gamma$, then
$x\vDash_{0}\B\gamma$ since $x' \sim x$. So $y\vDash_{0}\B\gamma$ by assumption. But then $y' \vDash_{0}\B\gamma$ since $y' \sim y$. This shows that $x' R_1 y'$, so $[x'] R_2 [y']$.
\end{proof}

In the following lemma we highlight some properties of $\mathfrak{F}_{2}$.

\begin{lem}\label{F2 props}
    \begin{enumerate}
        \item[]
        \item\label{R2 lem - R2 partial order} $R_{2}$ is a partial order.
        \item\label{R2 lem - E2 equiv} $E_{2}$ is an equivalence relation.
        \item\label{R2 lem - commutativity} $R_{2}$ and $E_{2}$ satisfy commutativity.
        \item\label{R2 lem - clean clusters} $\mathfrak{F}_{2}$ has clean clusters.
        \item\label{R2 lem - cluster bound} For $[x]\in W_{2}$, $|E_{2}([x])|\leq 2^n$, where $n=|\Sub |$.
        \item\label{R2 lem - can find y R1 max} For $[x]\in W_{2}$ and $\B\gamma\in\Sub$, if $x\not\vDash_{1}\B\gamma$, then there is $[y]\in W_{2}$ such that $[x]R_{2}[y]$ and $y\not\vDash_{1}\B\gamma$ $R_{1}$-maximally.
    \end{enumerate}
\end{lem}

\begin{proof}
(1) Reflexivity of $R_{2}$ is immediate from the definition, and transitivity and antisymmetry follow from transitivity and antisymmetry of $R_1$.

(2) This follows from $E_1$ being an equivalence relation.

(3) This follows from $R_{1}$ and $E_{1}$ satisfying commutativity.

(4) Suppose there are $[x]\neq [y]$ in $W_2$ with $[x]R_{2}[y]$ and  $[x]E_{2}[y]$. Then $x \neq y$, so by the definition of $R_{2}$ and $E_{2}$, we have $xR_{1}y$ and $xE_{1}y$ which yields a dirty cluster in $\mathfrak F_1$, contradicting Lemma \ref{F1 props}(8).

(5) This follows from the fact that there are at most $2^n$ $\sim$-equivalence classes in each cluster (see, e.g., \cite[Prop.~5.24]{Modal_Logic}).

(6) Suppose $x\not\vDash_{1}\B\gamma$. By Lemma \ref{truth F1}, $x\not\vDash_{0}\B\gamma$, so by Lemma \ref{F1 props}(9), there is $y\in W_{1}$ such that $xR_{1}y$, $y\in A\cap \max_{R_0}E_{0}(A)$, and $y\not\vDash_{0}\B\gamma$ $R_{0}$-maximally, where
\[
A=\nu_{0}(\neg\B\gamma)\cap\bigcap\{\nu_{0}(\B\psi) \mid \B\psi\in\Sub \mbox{ and } x\vDash_{0}\B\psi\}.
\]
Then $[x]R_{2}[y]$ and by Lemma~\ref{truth F1}, $y\not\vDash_{1}\B\gamma$. We show that $y$ is $R_{1}$-maximal with respect to $\B\gamma$. Suppose $yR_{1}z$ and $z\not\vDash_{1}\B\gamma$. By Lemma~\ref{truth F1}, $z\not\vDash_{0}\B\gamma$, and from $yR_{1}z$ it follows that $y=z$ or $y\strict{Q_{0}}z$ and $y\vDash_{0}\B\psi$ implies $z\vDash_{0}\B\psi$ for all $\B\psi\in\Sub$. Suppose the latter. Since $z\not\vDash_{0}\B\gamma$, we have $z\in \nu_{0}(\neg\B\gamma)$. If $x\vDash_{0}\B\psi$ for $\B\psi\in\Sub$, then $y\in A$ implies $y\vDash_{0}\B\psi$. So $yR_{1}z$ then gives $z\vDash_{0}\B\psi$. Therefore, $z\in \bigcap\{\nu_{0}(\B\psi) \mid \B\psi\in\Sub \mbox{ and } x\vDash_{0}\B\psi\}$, and hence $z\in A$. As $y\strict{Q_{0}}z$, there is $w\in W_{0}$ such that $y\neq w$, $yR_{0}w$, and $wE_{0}z$. Then $w\in E_{0}(A)$, and maximality of $y$ in $E_{0}(A)$ yields $y=w$, contradicting $y\neq w$. Thus, $y=z$, and so $y$ is $R_{1}$-maximal with respect to $\B\gamma$.
\end{proof}

We conclude Step 2 by showing the truth lemma for $\mathfrak{F}_{2}$.

\begin{lem}[Truth Lemma]\label{truth F2}
For $[x]\in W_{2}$ and $\psi\in\Sub$,
\[
(\mathfrak{F}_{1},x)\vDash_{1} \psi \Leftrightarrow (\mathfrak{F}_{2},[x])\vDash_{2}\psi.
\]
\end{lem}

\begin{proof}
The proof is by induction on the complexity of $\psi$. The base case $\psi=p$ follows from the definition of $\nu_{2}$. The cases of $\psi=\psi_{1}\wedge\psi_{2}$ and $\psi=\neg\psi_{1}$ are straightforward, and the $\forall$-case follows from the definition of $E_2$. Suppose that $\psi = \Box \psi_1$. If $x\not\vDash_{1}\B\psi_{1}$, then there is $y\in W_1$ with $xR_{1}y$ and $y\not\vDash_{1}\psi_{1}$. Therefore, $[x]R_{2}[y]$ and  $[y]\not\vDash_{2}\psi_{1}$ by the inductive hypothesis. Thus, $[x]\not\vDash_{2}\B\psi_{1}$. Conversely, if $[x]\not\vDash_{2}\B\psi_{1}$, then there is $y\in W_1$ with $[x] R_2 [y]$ and $[y]\not\vDash_{2}\psi_{1}$. By the inductive hypothesis, $y \not\vDash_{1} \psi_{1}$. If $[x] = [y]$, then $x \not \vDash_{1} \B \psi_{1}$ by definition of $\sim$. If $[x] \neq [y]$, then $x R_1 y$  and again $x \not \vDash_{1} \B \psi_{1}$.
\end{proof}

\subsection{Step 3: Constructing $\mathfrak{F}_{3}$}

We are ready for our final step, in which we construct $\mathfrak{F}_{3}=(W_{3},R_{3},E_{3})$ by selective filtration from $\mathfrak{F}_{2}$. This is done by constructing a sequence of finite partially ordered $\mathsf{MS4}$-frames $\mathfrak{F}_{3.h}=(W_{3.h},R_{3.h},E_{3.h})$ with clean clusters so that $\mathfrak{F}_{3.h}\subseteq \mathfrak{F}_{3.h+1}$ for all $h<\omega$. We then show that this construction eventually terminates.

Similar to the construction for $\mathsf{M^{+}IPC}$, for each point $[x]\in W_{2}$ that we select, we create a copy of the point, give it a new name, say $t$, and let $\widehat{t}=[x]$ denote the original point in $W_{2}$ that $t$ represents and will behave similar to.
However, we take a bit more care with the copies in this construction than in the construction for $\mathsf{M^{+}IPC}$. In particular, we will never create two copies of the same original point within one cluster. This will ensure that the cluster size in $\mathfrak{F}_{3}$ has the same bound as the cluster size in $\mathfrak{F}_{2}$. 

Before we begin the construction, we highlight an important property we will need for selecting our points.

\begin{lem}\label{R2 max point}
For $[x]\in W_{2}$ and $\B\gamma\in\Sub$, if $[x]\not\vDash_{2}\B\gamma$, then there is $[y]\in W_{2}$ such that $[x]R_{2}[y]$ and $[y]\not\vDash_{2}\B\gamma$ $R_{2}$-maximally.
\end{lem}

\begin{proof}
Suppose $[x]\not\vDash_{2}\B\gamma$. By Lemma~\ref{truth F2}, $x\not\vDash_{1}\B\gamma$, and by Lemma~\ref{F2 props}(6), there is $[y]\in W_{2}$ such that $[x]R_{2}[y]$ and $y\not\vDash_{1}\B\gamma$ $R_{1}$-maximally. Applying Lemma~\ref{truth F2} again yields  $[y]\not\vDash_{2}\B\gamma$. To see that $[y]$ is $R_{2}$-maximal with respect to $\B\gamma$, suppose $[y]R_{2}[z]$ and $[z]\not\vDash_{2}\B\gamma$. By definition of $R_{2}$, either $[y]=[z]$ or $yR_{1}z$. If $yR_{1}z$, then by $R_{1}$-maximality of $y$, we have $y=z$, so $[y]=[z]$, and hence $[y]$ must be $R_{2}$-maximal with respect to $\B\gamma$.
\end{proof}

Throughout the construction, for each $t\in W_{3.h}$, we associate the following sets of subformulas:
\begin{flalign*}
&&\Sigma^{\forall}(t) &= \{\forall\delta\in\Sub\mid\widehat{t}\not\vDash_{2}\forall\delta\text\}&\\
&&\Sigma^{\B}(t) &= \{\B\gamma\in\Sub\mid\widehat{t}\not\vDash_{2}\B\gamma, \widehat{t}\vDash_{2}\gamma\}.&
\end{flalign*}

We start with $\mathfrak{F}_{3.0}=(W_{3.0},R_{3.0},E_{3.0})$ where
\[W_{3.0}=\{t_{0}\},\quad R_{3.0}=W_{3.0}^{2}, \quad E_{3.0}=W_{3.0}^{2}\]
and $\widehat{t_{0}}=[x_0]\in W_{2}$ is a point with $[x_0]\not\vDash_{2}\varphi$. This will be a root of our frame and has $Q_{3}$-depth $1$. Let $W_{3.-1}=R_{3.-1}= E_{3.-1}= \varnothing$. Suppose $\mathfrak{F}_{3.h-1}=(W_{3.h-1},R_{3.h-1},E_{3.h-1})$ has already been constructed and is a partially ordered $\mathsf{MS4}$-frame with clean clusters. We construct $\mathfrak{F}_{3.h}$ by the following steps.

\vspace{1em}

\noindent\emph{Step $\forall$ $($Horizontal$)$}:
Let $W_{3.h}^{\A}=W_{3.h-1}$, $R_{3.h}^{\A}=R_{3.h-1}$, and $E_{3.h}^{\A}=E_{3.h-1}$.
For each cluster $E_{3.h}(t)\subseteq W_{3.h-1}\backslash W_{3.h-1}^{\A}$, consider $\forall\delta\in\Sigma^{\forall}(t)$. If there is no $s\in W_{3.h}^{\A}$ already such that $tE_{3.h}^{\A}s$ and $\widehat{s}\not\vDash_{2}\delta$, we add a witness to our new frame as follows. Since $\widehat{t}\not\vDash_{2}\forall\delta$, there exists $[x]\in W_{2}$ such that $\widehat{t}E_{2}[x]$ and $[x]\not\vDash_{2}\delta$. We add the point $s$ to $W_{3.h}^{\A}$ where $\widehat{s}=[x]$ ($s$ is a distinct new copy of $[x]$), the relations $(s,s)$ to $R_{3.h}^{\A}$, the relations $(t,s)$ to $E_{3.h}^{\A}$ and generate the least equivalence relation.

\vspace{1em}

\noindent\emph{Step $\B$ $($Vertical$)$}: 
Let $W_{3.h}^{\B}=W_{3.h}^{\A}$, $R_{3.h}^{\B}=R_{3.h}^{\A}$, and $E_{3.h}^{\B}=E_{3.h}^{\A}$. For $t\in W_{3.h}^{\forall}\backslash W_{3.h-1}^{\A}$ (hence including any points added in the horizontal step), consider $\B\gamma\in\Sigma^{\B}(t)$ where  $\widehat{t}\not\vDash_{2}\B\gamma$, but $\widehat{t}\vDash_{2}\gamma$ (thus, $t$ isn't witnessing the formula $\B\gamma$ itself), and there is no $s\in W_{3.h}^{\B}$ already such that $tR_{3.h}^{\B}s$ and $\widehat{s}\not\vDash_{2}\B\gamma$ $R_{2}$-maximally (such an $s$ could have been added in a previous stage to satisfy commutativity). For each such $\B\gamma$, since $\widehat{t}\not\vDash_{2}\B\gamma$ and $\widehat{t}=[w]$ for some $[w]\in W_{2}$, we have $[w]\not\vDash_{2}\B\gamma$. By Lemma~\ref{R2 max point}, there is $[x]\in W_{2}$ such that $[w]R_{2}[x]$ and $[x]$ is $R_{2}$-maximal with respect to $\B\gamma$.
We add the point $s$ to $W_{3.h}^{\B}$ where $\widehat{s}=[x]$, $(t,s)$ and $(s,s)$ to $R_{3.h}^{\B}$ and close under transitivity,
and add $(s,s)$ to $E_{3.h}^{\B}$.
To make sure commutativity is satisfied, for each $w\in E_{3.h}^{\B}(t)$, if there is already $s_{w}\in E_{3.h}^{\B}(s)$ such that $\widehat{w}R_{2}\widehat{s_{w}}$, we simply add the relation $(w, s_{w})$ to $R_{3.h}^{\B}$. If there is no such $s_{w}$, then by commutativity in $W_{2}$, there is $[x_{w}]\in W_2$ such that $\widehat{w}R_{2}[x_{w}]$ and $[x_{w}]E_{2}[x]$, so we add $s_{w}$ to $W_{3.h}^{\B}$, where $\widehat{s_{w}}=[x_{w}]$. We then add $(w,s_{w})$ to $R_{3.h}^{\B}$ and close it under reflexivity and transitivity, and add $(s_{w},s)$ to $E_{3.h}^{\B}$ and generate the smallest equivalence relation.

To end this stage of the construction, we let $\mathfrak{F}_{3.h}=(W_{3.h},R_{3.h},E_{3.h})$ where
 $$W_{3.h}=  W_{3.h}^{\B}, R_{3.h}= R_{3.h}^{\B} \text{ and } E_{3.h}=  E_{3.h}^{\B}.$$

\begin{lem}\label{finitegrzframe} $\mathfrak{F}_{3.h}$ is a finite partially ordered $\mathsf{MS4}$-frame with clean clusters.
\end{lem}

\begin{proof}
In the $\forall $-step we only added reflexive arrows to $R_{3.h}^{\A}$, so $R_{3.h}^{\A}$ is a partial order. In the $\Box$-step we close $R_{3.h}^{\Box}$ under reflexivity and transitivity each time we add a new arrow, so $R_{3.h}^{\Box}$ is reflexive and transitive. Moreover, we we only add $R_{3.h}^{\Box}$ arrows from points that were already present in $W_{3.h}^{\A}$ into points that are freshly added in the $\Box$-step of round $h$. Thus, $R_{3.h}^{\Box}$ is antisymmetric. 
That $E_{3.h}$ is an equivalence relation and that $\mathfrak{F}_{3.h}$ satisfies commutativity follow from the construction.
Finally, to see that $\mathfrak{F}_{3.h}$ has only clean clusters, note that in the $\forall$-step all freshly introduced $E_h$-relations are of the shape $(s, t)$ where $s$ or $t \in W^{\forall}_{3.h} \setminus W_{3.h-1}$. Since no non-reflexive $R_h$-arrows are introduced in this step, no dirty cluster could have been built. We have already discussed the shape of the $R_{h}$ arrows introduced in the $\B$-step.
This guarantees that no cluster in $W^{\forall}_{3.h}$ is made dirty. The freshly introduced $E_h$-relations in these steps are of the shape $(s, t)$ where $s, t  \in W^{\B}_{3.h} \setminus W^{\A}_{3.h}$. Since no non-reflexive $R_h$-relations exist between these points, we infer that all clusters are clean.
\end{proof}

The following lemma summarizes some useful properties of $\mathfrak{F}_{3}$. In the following let
\begin{equation*}
n=\lvert\Sub\rvert\ \text{ and } m= \lvert\{\B\psi:\B\psi\in \Sub\}\rvert.
\end{equation*}

\begin{lem}\label{mGrz_consistency of sigmas} Let $t,u\in W_{3.h}$.
    \begin{enumerate}
        \item\label{E-pres} If $tE_{3.h}u$, then $\widehat{t}E_{2}\widehat{u}$.
        \item\label{mGrz cons - A step once per cluster} If $t E_{3.h} u$, then $\Sigma^{\forall}(t) =\Sigma^{\forall}(u)$. $($This ensures that we only need to perform the $\A$-step once per cluster$)$.
        \item\label{mGrz cons - no two instances of the same pt in one cluster} If $t E_{3.h}u$ and $t\neq u$, then $\widehat{t}\neq\widehat{u}$. $($This ensures that one cluster does not contain two different copies of the same point, so our cluster size remains bounded$)$.
        \item\label{mGrz cons - strict Q2 implies strict Q0} If $[t]\strict{Q_{2}}[u]$, then $t\strict{Q_{0}}u$.
        \item\label{R-pres} If $tR_{3.h}u$, then $\widehat{t}R_{2}\widehat{u}$.
        \item\label{mGrz cons - strict Q3.h implies strict Q2} If $t \strict{Q_{3.h}} u$, then $\widehat{t}\strict{Q_{2}}\widehat{u}$. Thus, if $t\strict{Q_{3.h}}u$, then $\widehat{t}\neq\widehat{u}$.
        \item\label{mGrz cons - box witness bound} A formula $\B\gamma\in\Sub$ can be witnessed at most $2^n$ times in clusters along an $R_{3.h}$-chain. $($This shows that  $\B\gamma$ can be witnessed at most $2^n$ times per $Q_{3.h}$-chain.$)$
    \end{enumerate}
\end{lem}

\begin{proof}
(1) This follows from the construction.

(2) By (1), $t E_{3.h}u$ implies $\widehat{t} E_{2} \widehat{u}$, so $\widehat{t} \vDash_{2} \forall \gamma$ iff $\widehat{u} \vDash_{2}\forall \gamma$.

(3) Suppose $t E_{3.h}u$ and $\widehat{t}=\widehat{u}$, and without loss of generality assume that $t$ was added to the cluster before $u$, so either $u$ is added to witness some formula $\A\delta_{i}$ where $\widehat{u}\not\vDash_{2}\delta_{i}$, or $u$ is added as a commutativity witness for some point from the cluster below. However, by construction, $u$ would not have been added to witness a formula $\A\delta_{i}$, because if $\widehat{u}\not\vDash_{2}\delta_{i}$, then $\widehat{t}=\widehat{u}$ implies that $\widehat{t}\not\vDash_{2}\delta_{i}$, so $t$ is already a viable witness in the cluster for any such formula, contradicting the $\A$-step of the construction. Furthermore, $u$ would not be added as a commutativity witness for some point $w$ in the cluster immediately below, because then in $W_{2}$ we would have $\widehat{w}R_{2}\widehat{u}$, so $\widehat{w}R_{2}\widehat{t}$, and a new $R_{3.h}$-relation would have been added from $w$ to $t$ instead, contradicting the $\B$-step of the construction. Thus, we must have $\widehat{t}\neq\widehat{u}$.

(4) Suppose $[t]\strict{Q_{2}}[u]$. Then there is $[w]\in W_{2}$ with $[t]\neq[w]$, $[t]R_{2}[w]$, and $[w]E_{2}[u]$. From the definitions of $E_{1}$ and $E_{2}$, $[w]E_{2}[u]$ implies $wE_{0}u$. By definition of $R_{2}$, $[t]R_{2}[w]$ and $[t]\neq[w]$ imply $tR_{1}w$. Since $[t]\neq[w]$, we have $t\neq w$, so $t\strict{Q_{0}}w$ by the definition of $R_{1}$. Then there is $v\in W_{1}$ with $t\neq v$, $tR_{0}v$, and $vE_{0}w$. Since $vE_{0}w$, we have $vE_{0}u$. Thus, $t\neq v$, $tR_{0}v$, and $vE_{0}u$, and hence $t\strict{Q_{0}}u$.

(5) This follows from the construction.

(6) If $t \strict{Q_{3.h}} u$, then there is $w$ such that $t\neq w$, $tR_{3.h}w$, and $wE_{3.h}u$. By (5), $\widehat{t}R_{2}\widehat{w}$ and $\widehat{w}$ must come from a different cluster in $W_{2}$ than $\widehat{t}$, so $\widehat{t}\neq\widehat{w}$. We also have $\widehat{w}E_{2}\widehat{u}$ by (1), so $\widehat{t}\strict{Q_{2}}\widehat{u}$. Because $\mathfrak{F}_{2}$ has clean clusters, we must have $\widehat{t}\neq\widehat{u}$.

(7) Suppose that $x_1,...,x_{2^{n}+1}$ are all in different $E_{3.h}$-clusters along an $R_{3.h}$-chain (where $\widehat{x_1}=[w_{1}],...,\allowbreak\widehat{x}_{2^{n}+1}=[w_{2^{n}+1}]$), so $x_1\strict{Q_{3.h}}...\strict{Q_{3.h}}x_{2^{n}+1}$, and all have been added to witness a formula $\B\gamma\in\Sub$. Thus, $\widehat{x_{i}}\not\vDash_{2}\B\gamma$ $R_2$-maximally for $i=1,...,2^{n}+1$. Because there are only $2^{n}$ subsets of $\Sub$ (where $n=\lvert\Sub\rvert$), the pigeonhole principle implies that there are some $i$ and $j$ with $i\neq j$ (assume $i<j$) where $\widehat{x_i}$ and $\widehat{x_j}$ satisfy the same subformulas of $\varphi$.
By (6), $\widehat{x_{i}}\strict{Q_{2}}\widehat{x_{j}}$ and $\widehat{x_{i}}\neq\widehat{x_{j}}$. If $\widehat{x_{i}}R_{2}\widehat{x_{j}}$, then $R_{2}$-maximality of $\widehat{x_{i}}$ with respect to $\B\gamma$ implies $\widehat{x_{i}}=\widehat{x_{j}}$, contradicting $\widehat{x_{i}}\neq\widehat{x_{j}}$, so we must have $\widehat{x_{i}}\nR_{2}\widehat{x_{j}}$ and hence $[w_{i}]\nR_{2}[w_{j}]$. Since $\widehat{x_{i}}\strict{Q_{2}}\widehat{x_{j}}$ we have $[w_{i}]\strict{Q_{2}}[w_{j}]$. By (4), we then have $w_i\strict{Q_{0}}w_j$. Since $[w_i]$ and $[w_j]$ satisfy the same formulas in $\Sub$, we have $[w_i]\vDash_{2}\B\beta\Leftrightarrow [w_j]\vDash_{2}\B\beta$ for $\B\beta\in\Sub$. By Lemmas~\ref{truth F2} and~\ref{truth F1}, $w_{i}\vDash_{0}\B\beta\Leftrightarrow w_{j}\vDash_{0}\B\beta$. Thus, $w_{i}R_{1}w_{j}$ and hence $[w_i]R_{2}[w_j]$, a contradiction.
\end{proof}

We now prove that the end result of our construction is a finite frame, using the definitions of bounded cluster size, bounded $R$-branching, and bounded $R$-depth given in Definition \ref{defn: cluster size_branching_depth}.

\begin{lem}\label{mGrz_branching}
$\mathfrak{F}_{3.h}=(W_{3.h},R_{3.h},E_{3.h})$ has cluster size bounded by $2^{n}$ for all $h<\omega$.
\end{lem}

\begin{proof}
By Lemma~\ref{F2 props}(5), the cluster size in $\mathfrak F_2$ is bounded by $2^{n}$, and by Lemma~\ref{mGrz_consistency of sigmas}(3), we do not add copies of the same points to a cluster in $\mathfrak{F}_{3.h}$. Thus, cluster size in $\mathfrak F_{3.h}$ is bounded by $2^{n}$.
\end{proof}

\begin{lem}\label{mGrz_branching}
$\mathfrak{F}_{3.h}=(W_{3.h},R_{3.h},E_{3.h})$ has $R_{3.h}$-branching bounded by $2^{n}\cdot m$ for all $h<\omega$.
\end{lem}

\begin{proof}
It is sufficient to show that each $t\in W_{3.j-1}$, for $j\leq h$, has at most $2^{n}\cdot m$ immediate $R_{3.j}$-successors. By construction, we add at most $m$-many immediate $R_{3.j}$-successors to $t$ for formulas of the form $\B\psi\in\Sub$. Each $y\in E_{3.j}(t)$ also needs at most $m$-many immediate $R_{3.j}$-successors to witness $\B$-formulas. Since there are at most $2^{n}$-many such $y$ (including $t$ itself), we must add an immediate $R_{3.j}$ successor to $t$ for commutativity for at most $2^{n}\cdot m$ points. Thus, $t$ has at most $2^{n}\cdot m$ immediate $R_{3.j}$-successors.
    \begin{center}
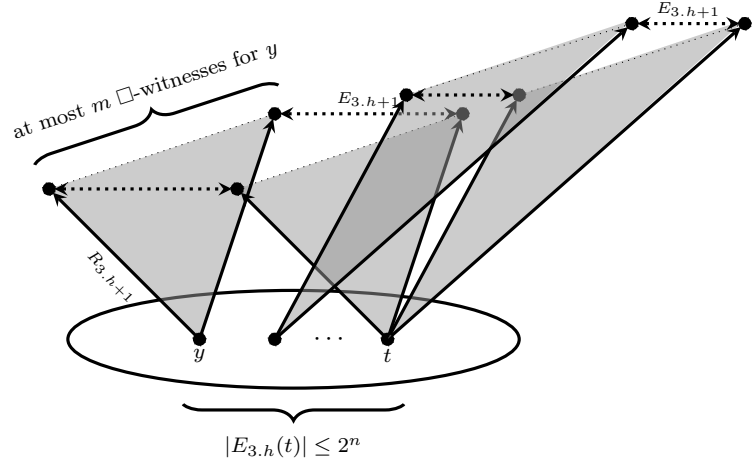
        \begin{minipage}{0.6\textwidth}
            \vspace{1pt}
            \begin{tikzpicture}
            [shorten <=2pt,shorten >=2pt,>=stealth]
            \draw [very thick] (5.25,0) ellipse (3cm and 0.65cm);
            \path[fill=gray!80, fill opacity=0.5] (4,0) -- (2,2) -- (5,3) -- cycle;
            \path[fill=gray!80, fill opacity=0.5] (6.5,0) -- (4.5,2) -- (7.5,3) -- cycle;
            \filldraw [dotted]
            (4,0) circle (2.5pt) node[align=center,   below] {\footnotesize{$y$}} --
            (2,2) circle (2.5pt) node[align=center, below] {}     --
            (5,3) circle (2.5pt) node[align=right,  below] {};
            \filldraw [dotted]
            (6.5,0) circle (2.5pt) node[align=center,   below] {\footnotesize{$t$}} --
            (4.5,2) circle (2.5pt) node[align=center, below] {}     --
            (7.5,3) circle (2.5pt) node[align=right,  below] {};
            \draw[decorate,decoration={brace,mirror,amplitude=8pt}, very thick] (3.75,-0.8) -- (6.75,-0.8) node [midway, below, yshift=-10pt] {\footnotesize{$\lvert E_{3.h}(t)\rvert\leq 2^{n}$}};
            \draw[->, very thick] (6.5,0) -- (4.5,2);
            \draw[->, very thick] (6.5,0) -- (7.5,3);
            \draw[->, very thick] (4,0) -- (2,2) node [midway, sloped, below] {\tiny{$R_{3.h+1}$}};
            \draw[->, very thick] (4,0) -- (5,3);
            \draw[dotted, <->, very thick] (2,2) -- (4.5,2);
            \draw[dotted, <->, very thick] (5,3) -- (7.5,3) node [midway,above, yshift=-2pt] {\tiny{$E_{3.h+1}$}};
            \draw[decorate,decoration={brace,amplitude=8pt}, very thick] (1.8,2.25) -- (5.1,3.33) node [midway, above, sloped, yshift=8pt] {\footnotesize{at most $m$ $\B$-witnesses for $y$}};
            \node (dots) at (5.75,0) {$\cdots$};

            \path[fill=gray!80, fill opacity=0.5] (6.5,0) -- (8.25, 3.25) -- (11.25,4.25) -- cycle;
            \draw[->, very thick] (6.5,0) -- (8.25,3.25);
            \draw[->, very thick] (6.5,0) -- (11.25,4.2);
            \filldraw [dotted]
            (6.5,0) --
            (8.25,3.25) circle (2.5pt) node[align=center, below] {}     --
            (11.25,4.2) circle (2.5pt) node[align=right,  below] {};

            \path[fill=gray!80, fill opacity=0.5] (5,0) -- (6.75, 3.25) -- (9.75,4.25) -- cycle;
            \draw[->, very thick] (5,0) -- (6.75,3.25);
            \draw[->, very thick] (5,0) -- (9.75,4.2);
            \filldraw [dotted]
            (5,0) circle (2.5pt) node[align=center, below] {} --
            (6.75,3.25) circle (2.5pt) node[align=center, below] {}     --
            (9.75,4.2) circle (2.5pt) node[align=right,  below] {};
            \draw[dotted, <->, very thick] (6.75,3.25) -- (8.25,3.25);
            \draw[dotted, <->, very thick] (9.75,4.2) -- (11.25,4.2) node [midway,above, yshift=-2pt] {\tiny{$E_{3.h+1}$}};

            \end{tikzpicture}
        \end{minipage}
    \end{center}
\end{proof}

\begin{lem}\label{mGrz_depth}
    $\mathfrak{F}_{3.h}=(W_{3.h},R_{3.h},E_{3.h})$ has $R_{3.h}$-depth bounded by $2^{n}\cdot m +1$ for all $h<\omega$.
\end{lem}

\begin{proof}
By construction, to make an immediate vertical move from some cluster $E_{3.h}(t)$ to another cluster $E_{3.h}(u)$ (with $t\neq u$), there must be some point $x\in E_{3.h}(t)$ and formula $\B\psi\in\Sigma^{\B}(x)$ requiring a witness $y$, where $y\in E_{3.h}(u)$, $xR_{3.h}y$, and $y$ is added in the $\B$-step of the construction. Starting from the bottom cluster $E_{3.h}(t_{0})$, by Lemma~\ref{mGrz_consistency of sigmas}(7), each of our $m$-many $\B$-formulas can be witnessed at most $2^{n}$ times in clusters along an $R_{3.h}$-chain. Thus, we add at most $2^n \cdot m$ elements to an $R_{3.h}$ chain originating from this cluster, with the total length of the chain (including the starting point) being at most $2^{n}\cdot m +1$.
    \begin{center}
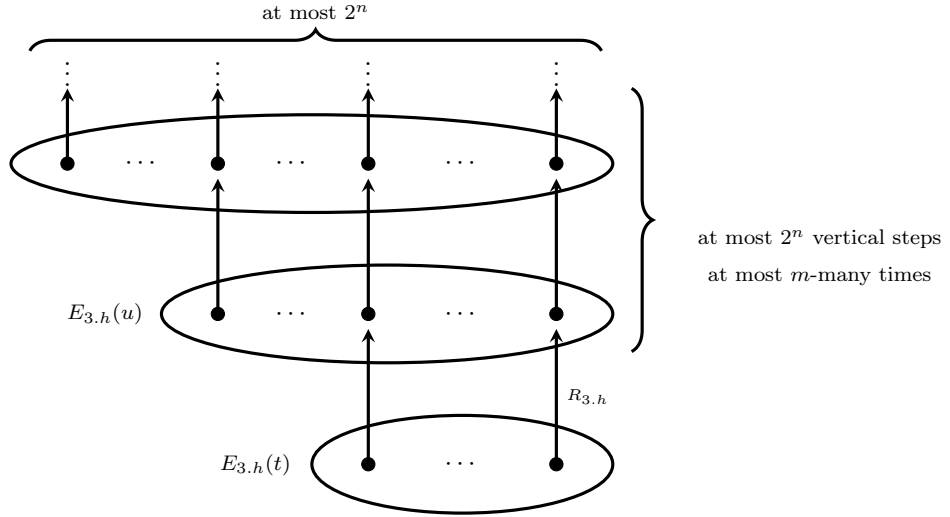
        \begin{minipage}{0.6\textwidth}
            \vspace{1pt}
            \begin{tikzpicture}
            [>=stealth]
            \node (Eht) at (2.5,0) {\footnotesize{$E_{3.h}(t)$}};
            \draw [very thick] (5.25,0) ellipse (2cm and 0.65cm);
            \filldraw (6.5,0) circle (2.5pt) node[align=center,   below] {};
            \node (f) at (5.25,0) {$\cdots$};
            \filldraw (4,0) circle (2.5pt) node[align=center,   below] {};
            \draw[->, very thick] (6.5,0) -- (6.5,1.8) node [midway,right] {\tiny{$R_{3.h}$}};
            \draw[->, very thick] (4,0) -- (4,1.8);

            \node (Ehu) at (0.5,2) {\footnotesize{$E_{3.h}(u)$}};
            \draw [very thick] (4.25,2) ellipse (3cm and 0.65cm);
            \filldraw (6.5,2) circle (2.5pt) node[align=center,   below] {};
            \node (d) at (5.25,2) {$\cdots$};
            \filldraw (4,2) circle (2.5pt) node[align=center,   below] {};
            \node (g) at (3,2) {$\cdots$};
            \filldraw (2,2) circle (2.5pt) node[align=center,   below] {};
            \draw[->, very thick] (6.5,2) -- (6.5,3.8);
            \draw[->, very thick] (4,2) -- (4,3.8);
            \draw[->, very thick] (2,2) -- (2,3.8);

            \draw [very thick] (3.25,4) ellipse (4cm and 0.65cm);
            \filldraw (6.5,4) circle (2.5pt) node[align=center,   below] {};
            \node (h) at (5.25,4) {$\cdots$};
            \filldraw (4,4) circle (2.5pt) node[align=center,   below] {};
            \node (i) at (3,4) {$\cdots$};
            \filldraw (2,4) circle (2.5pt) node[align=center,   below] {};
            \node (j) at (1,4) {$\cdots$};
            \filldraw (0,4) circle (2.5pt) node[align=center,   below] {};
            \draw[->, very thick] (6.5,4) -- (6.5,5);
            \draw[->, very thick] (4,4) -- (4,5);
            \draw[->, very thick] (2,4) -- (2,5);
            \draw[->, very thick] (0,4) -- (0,5);
            \node (k) at (6.5,5.3) {$\vdots$};
            \node (l) at (4,5.3) {$\vdots$};
            \node (m) at (2,5.3) {$\vdots$};
            \node (n) at (0,5.3) {$\vdots$};

            \draw[decorate,decoration={brace,amplitude=8pt}, very thick] (7.5,5) -- (7.5,1.5);

            \node (p) at (10,3) {\footnotesize{at most $2^{n}$ vertical steps}};
            \node (r) at (10,2.5) {\footnotesize{at most $m$-many times}};

            \draw[decorate,decoration={brace,amplitude=8pt}, very thick] (-0.5,5.5) -- (7.1,5.5) node [midway, above, yshift=8pt] {\footnotesize{at most $2^{n}$}};
            \end{tikzpicture}
        \end{minipage}
    \end{center}
\end{proof}

\begin{lem}\label{lem:stopgrz}
There is $h \in \omega$ such that $\mathfrak F_{3.h'} = \mathfrak F_{3.h}$ for all $h' \geq h$.
\end{lem}

\begin{proof}
As in the proof of Lemma~\ref{lem:stop}, we observe that in stage $k$ of the construction, all $R_{3.h}$-chains are bounded by $k$. Since, by Lemma~\ref{mGrz_depth}, the $R_{3.h}$-depth of $\mathfrak F_{3.h}$ is bounded by $2^{n}\cdot m +1$, we have $\mathfrak F_{3.h'} = \mathfrak F_{2^{n}\cdot m +1}$ for all $h' \geq 2^{n}\cdot m +1$.
\end{proof}

Set $\mathfrak{F_3}=(W_3,R_3,E_3)$ where
\[
W_3=W_{3.h},\quad R_3=R_{3.h}, \quad E_3=E_{3.h},
\]
and $h$ is as in Lemma~\ref{lem:stopgrz}.
As an immediate consequence of Lemma~\ref{finitegrzframe}, we obtain: 

\begin{lem}
$\mathfrak{F}_{3}=(W_{3},R_{3},E_{3})$ is a finite $\mathsf{M^{+}Grz}$-frame.
\end{lem}

Finally, we verify that our frame validates precisely the formulas we want it to. Define a valuation $\nu_{3}$ on $W_{3}$ by $\nu_{3}(p)=\{t\in W_{3}:\widehat{t}\in \nu_{2}(p)\}$ for $p\in\Sub$ and $\nu_{3}(q)=\varnothing$ for variables $q$ not occurring in $\varphi$.

\begin{lem}[Truth Lemma]\label{truth F3}
For $x\in W_{3}$ and $\psi\in\Sub$,
\[
(\mathfrak{F}_{2},\widehat{x})\vDash_{2} \psi \Leftrightarrow (\mathfrak{F}_{3},x)\vDash_{3}\psi.
\]
\end{lem}

\begin{proof}
The proof is by induction on the complexity of $\psi$ and again we only show the cases where $\psi=\A\psi_{1}$ or $\psi=\B\psi_{1}$.

Suppose $\psi=\A\psi_{1}$. If $\widehat{x}\not\vDash_{2}\A\psi_{1}$, then $\forall\psi_{1}\in\Sigma^{\forall}(x)$, so at some point in the construction of $\mathfrak{F}_{3}$ we add $s$ to $W_{3}$ and $(x,s)$ to $E_{3}$ where $\widehat{s}\not\vDash_{2}\psi_{1}$. By the inductive hypothesis, $s\not\vDash_{3}\psi_{1}$, hence $x\not\vDash_{3}\forall\psi_{1}$.
Conversely, if $x\not\vDash_{3}\forall\psi_{1}$, then there is $w\in W_{3}$ with $xE_{3}w$ and $w\not\vDash_{3}\psi_{1}$. By the inductive hypothesis, $\widehat{w}\not\vDash_{2}\psi_{1}$, and by Lemma~\ref{mGrz_consistency of sigmas}(1), $xE_{3}w$ implies $\widehat{x}E_{2}\widehat{w}$, so $\widehat{x}\not\vDash_{2}\forall\psi_{1}$.

Suppose $\psi=\B\psi_{1}$. If $\widehat{x}\not\vDash_{2}\B\psi_{1}$, then either $\widehat{x}\not\vDash_{2}\psi_{1}$ or $\widehat{x}\vDash_{2}\psi_{1}$. If $\widehat{x}\not\vDash_{2}\psi_{1}$, then by the inductive hypothesis we have $x\not\vDash_{3}\psi_{1}$, hence $x\not\vDash_{3}\B\psi_{1}$. If $\widehat{x}\vDash_{2}\psi_{1}$, then $\B\psi_{1}\in\Sigma^{\B}(x)$, so at some point in the construction of $\mathfrak{F}_{3}$ we add $s$ to $W_{3}$ and $(x,s)$ to $R_{3}$ where $\widehat{s}\not\vDash_{2}\psi_{1}$. By the inductive hypothesis, $s\not\vDash_{3}\psi_{1}$, hence $x\not\vDash_{3}\B\psi_{1}$. Conversely, if $x\not\vDash_{3}\B\psi_{1}$, then there is $w\in W_{3}$ with $xR_{3}w$ and $w\not\vDash_{3}\psi_{1}$. By the inductive hypothesis, $\widehat{w}\not\vDash_{2}\psi_{1}$, and by Lemma~\ref{mGrz_consistency of sigmas}(5), $xR_{3}w$ implies $\widehat{x}R_{2}\widehat{w}$, so $\widehat{x}\not\vDash_{2}\B\psi_{1}$.
\end{proof}

The three steps of our construction yield our desired result:

\begin{theorem}\label{FMP M+Grz}
$\mathsf{M^{+}Grz}$ has the finite model property.
\end{theorem}

As an immediate corollary to Theorem~\ref{FMP M+Grz}, we obtain:

\begin{cor}
$\mathsf{M^{+}Grz}$ is decidable.
\end{cor}

\begin{remark}
Another consequence of Theorem~\ref{FMP M+Grz} is that $\mathsf{M^{+}Grz}$ is the monadic fragment of the predicate modal logic obtained by adding to $\mathsf{QGrz}$ the G\"odel translation of Casari's formula $\mathsf{Cas}$ (cf.~Remark~\ref{rem:Ono-Suzuki}). 
\end{remark}

\bibliographystyle{plain}
\bibliography{References}

\end{document}